\input amstex
\input epsf
\documentstyle{amsppt}
\TagsOnRight
\def\epsfsize#1#2{.75\hsize}

\def\ro#1{\expandafter\def\csname#1\endcsname{\operatorname{#1}}} \ro{mesh}
\ro{dim} \ro{id} \ro{dist} \ro{diam} \ro{min} \ro{tr} \ro{im} \ro{N} \ro{Fr}
\ro{Int} \ro{lk} \ro{sh} \ro{st} \def\derlim{\lim\limits_\leftarrow{^1}\,}
\def\Cl#1{{\overline{#1}}} \def\col{{\ssize~\searrow~}} \redefine\o{\circ}
\let\alp\alpha \let\bta\beta \let\gma\gamma \let\dta\delta \let\eps\varepsilon
\let\lda\lambda \let\sma\sigma \let\phi\varphi \def\R{\Bbb R} \def\Z{\Bbb Z}
\let\emb\hookrightarrow \let\imm\looparrowright \let\x\times \let\but\setminus
\let\inv\leftrightarrow \let\tl\tilde \def\eq{\ind{eq}} \def\PL{\ind{PL}}
\def\ind#1{_{\text{\rm #1}}} \def\fig#1{\centerline{\epsffile{FIG#1.EPS}}}

\topmatter %%%%%%%%%%%%%%%%%%%%%%%%%%%%%%%%%%%%%%%%%%%%%%%%%%%%%%%%%%%%%%%%%%%%

\title On maps with unstable singularities \endtitle
\author Sergey A. Melikhov \endauthor

\abstract
If a continuous map $f\:X\to Q$ is approximable arbitrary closely by embeddings
$X\emb Q$, can some embedding be taken onto $f$ by a pseudo-isotopy?
This question, called Isotopic Realization Problem, was raised by
\v{S}\v{c}epin and Akhmet'ev.
We consider the case where $X$ is a compact $n$-polyhedron, $Q$ a PL
$m$-manifold and show that the answer is `generally no' for $(n,m)=(3,6)$;
$(1,3)$, and `yes' when:

1) $m>2n$, $(n,m)\neq (1,3)$;

2) $m>\frac{3(n+1)}2$ and $\Delta(f)=\{(x,y)\mid f(x)=f(y)\}$ has an
equivariant (with respect to the factor exchanging involution) mapping cylinder
neighborhood in $X\x X$;

3) $m>n+2$ and $f$ is the composition of a PL map and a TOP embedding.

\noindent
In doing this, we answer affirmatively (with a minor preservation) a question
of Kirby: does small smooth isotopy imply small smooth ambient isotopy in the
metastable range, verify a conjecture of Kearton--Lickorish: small PL
concordance implies small PL ambient isotopy in codimension $\ge 3$, and a
conjecture set of Repovs--Skopenkov.
\endabstract

\thanks Partially supported by the Russian Foundation for Basic Research
Grant No\. 99-01-00009 \endthanks
\address Moscow State University, Mech.--Math\. Dept., Division of Diff\.
Geometry; \newline Vorobjovy Gory, Moscow 119899, Russia \endaddress
\curraddr University of Florida, Department of Mathematics;\newline
358 Little Hall, PO Box 118105, Gainesville, FL 32611-8105, U.S. \endcurraddr
\email melikhov\@math.ufl.edu, sergey\@melikhov.mccme.ru \endemail
\endtopmatter

\document %%%%%%%%%%%%%%%%%%%%%%%%%%%%%%%%%%%%%%%%%%%%%%%%%%%%%%%%%%%%%%%%%%%%%

\head 1. Introduction. \endhead

A general mathematical problem is to decide whether a singular state of some
system is stable or unstable.
In terms of geometric topology it can be expressed as follows: given a
continuous map $f$ of a compactum $X$ into a manifold $Q$, can it be
$\eps$-approximated by an embedding $f_\eps\:X\emb Q$ for each $\eps>0$?
If this is the case, the map $f$ is called {\it realizable} \cite{Siek2}
or {\it discretely realizable} \cite{Akh1}.
If $f$ is a constant map, its realizability evidently coincides with
embeddability of $X$ into $\R^{\dim Q}$, meanwhile embeddability of a compactum
into $\R^m$ can be reduced to realizability of certain PL maps, cf\.
\cite{Siek2}, \cite{KW}, \cite{SS}, \cite{Akh1}.

As far as in some cases the $\eps$-approximation of $f$ can be made in infinite
number of inequivalent ways (i.e\. by embeddings, not joined by sufficiently
small ambient isotopies, -- e.g\. a map $S^n\sqcup S^n\to S^n\vee S^n\emb
\R^{2n+1}$, which embeds each sphere, can be approximated by links with
arbitrary linking number), it is natural to ask, whether an approximation of
$f$ can be viewed as a {\sl continuous process}, parametrized by real numbers?
The map $f$ is {\it isotopically realizable}, if there exists a homotopy
$H_t\:Q\to Q$, $t\in I=[0,1]$, such that $H_t$ is a homeomorphism for $t<1$
(such homotopy is called a {\it pseudo-isotopy}, cf\. \cite{Ke3}), $H_0=\id_Q$
and $H_1\o g=f$ for some embedding $g\:X\emb Q$.

\proclaim{Isotopic Realization Problem}
{\rm (E. V. \v{S}\v{c}epin, 1993; P. M. Akhmet'ev \cite{Akh1})}
When does discrete realizability imply isotopic realizability?
\endproclaim

To the best of the author's knowledge, the concept of isotopic realizability
was first considered by Blass and Holszty\'nski in 1971 \cite{Ho}.
It was independently introduced under the present name in a paper by
\v{S}\v{c}epin and \v{S}tan'ko \cite{SS} (subsequent to the work of
\v{S}\v{c}epin on uncountable inverse spectra and the earlier work of
\v{S}tan'ko on embedding dimension).
Although both papers \cite{Ho}, \cite{SS} also dealt with discrete
realizability, the relationship was not discussed there, and it appears that
until recently the IR Problem above has been virtually untouched.

Actually it traces back to the L. V. Keldy\v{s} Problem (1966) on
realizability of wildly embedded polyhedra by pseudo-isotopy of subpolyhedra
\cite{Ke1}, \cite{Ke2} (see also \cite{Ed2}, compare \cite{Mill4}).
That is, in the above $X$ should be replaced by a polyhedron, $Q$ by a PL
manifold, $f$ by an embedding and $g$ by a PL embedding, and, strictly
speaking, the pointwise equality $H_1\o g=f$ by the setwise $H_1(g(X))=f(X)$.
In a few succeeding years the Keldy\v{s} Problem was solved positively for wild
surfaces in $3$-manifolds \cite{Ke2} (see also \cite{Cr}) and for wild
$n$-polyhedra in PL $m$-manifolds, $m-n\ge 3$ \cite{Ed1} (cf\. Theorem 3.5a
below), and negatively for certain wild knots in $\R^3$ \cite{Ke4}, \cite{Sik}
(see Example 1.2 below).
On the other hand, it should be noticed that isotopic realizability as a
property of maps in the closure of the space of embeddings is similar to
tameness as a property of TOP embeddings in the closure of the space of
PL embeddings, moreover, in codimension $\ge 3$, the fact that all TOP
embeddings lie in the latter closure was used in proofs of equivalence of
tameness and the 1-LCC property \cite{BS1}, \cite{Ch2}, \cite{BS2} (see also
\cite{Ed1, 8.2}, \cite{Qu, 2.5.1}).

The concept of discrete realizability was studied widely (see brief
surveys in \cite{ReS1} and \cite{Akh3}).
For example, each self-map of the pseudo-arc is realizable (see \cite{Le}),
meanwhile for locally connected continua $X,Y$, $\dim X\leq 1$, all maps
$X\to Y\emb\R^2$ are realizable iff either $X$ is contained in triod and $Y$ is
in $S^1$, or $X$ is contained in the `letter q' and $Y$ is in $I$ \cite{Siek2}.
Realizability of a given map $X\to\R^2$ seems to be a harder question (see
\cite{Minc} for the PL case).
Any map of an $n$-dimensional compactum $X$ into $\R^m$ is realizable for
$m\geq 2n+1$ (cf\. \cite{Siek1}) and even for $m=2n$ if $\dim X\x X<2n$
\cite{DRS}, \cite{Sp}.
All maps $T^n\to T^n\emb_{\sssize\text{\rm standard}}~\R^{2n}$ are realizable
if $n>1$ \cite{KW}, meanwhile the maps
$S^n\to S^n\emb_{\sssize\text{\rm standard}}\R^{2n}$ are realizable whenever
$n\neq 1,2,3,7$, and are not, generally speaking, if $n=1,3$ or ~$7$
\cite{Akh1}, \cite{Akh3}.
Furthermore, for each $k$ one can find an $n$ such that all maps
$S^n\to S^n\emb_{\sssize\text{\rm st\.}}\R^{2n-k}$ are realizable \cite{Akh2},
\cite{Akh3}.
Surprisingly, in the space of maps
$S^2\to\R^3\emb_{\sssize\text{\rm st\.}}\R^4$ the subset of non-realizable maps
is dense \cite{AhRS}.

As for isotopic realizability, two principal results had been previously known.
From \v{C}ernavskij's Theorem on local contractibility of the
homeomorphism group of a closed manifold $M$ \cite{Ch3}, \cite{EK} it follows
that discrete realizability implies isotopic for self-maps of $M$ (thus for
$\dim M\neq 3$ both are equivalent to the property of being cell-like
\cite{Sieb2}, \cite{Qu}).
Secondly, Akhmet'ev showed in 1996 that all maps
$S^n\to S^n\emb_{\sssize\text{\rm st\.}}\R^{2n}$ are isotopically realizable
for $n=4k+1\ge 9$ \cite{Akh1}.

The real question, implicit in the above universal statement of the IR Problem
and originally motivating this deep result of Akhmet'ev as well as the present
paper was, does discrete realizability imply isotopic for maps of nice spaces
(say, of a manifold into Euclidean space) in high codimensions (say, greater
than $2$, in order to kill the fundamental group)?
Our main result is that it does not.
The counter-example (Example 1.9) is an explicit geometric construction with
a self-contained verification, but it was not until the rest of the paper had
been written when it naturally appeared.
The major part of the paper is devoted to the reduction of the IR Problem for
maps of nice spaces in the metastable range to a homotopy-theoretic question,
which, in turn, admits an answer in terms of vanishing of certain cohomological
obstructions.

Despite such an algebraization, it is still unknown\footnote{{\it Added in v.2:}
There has been some progress.
An example of a discretely but not isotopically realizable map
$S^n\to\R^{2n-1}\emb_{\sssize\text{\rm st.}}\R^{2n}$ for any even $n\ge 4$ is
included in \cite{Me2}, where it is also shown that if $n\ne 5,6$ and neither
$n$ nor $n+1$ is a power of $2$, then any map
$S^n\to S^n\emb_{\sssize\text{\rm st.}}\R^{2n}$ is isotopically realizable.}
whether the above condition $n=4k+1$ is really necessary, and even whether
discrete realizability implies isotopic for all maps
$S^n\to\R^{2n-1}\emb_{\sssize\text{\rm st.}}\R^{2n}$.
A new technique seems to be necessary here, which may be also useful in
attacking the following problem \cite{AhRS}: Suppose $f\:M\to\R^n$ is a generic
smooth map and $i\:\R^n\emb\R^{n+k}$ the standard inclusion, does discrete
realizability of $i\o f$ imply that $f$ can be factored into the composition of
an embedding $M\emb\R^{n+k}$ and the projection $\R^{n+k}\to\R^n$?
(The latter clearly implies isotopic realizability of $i\o f$.)

\subhead A. Low-codimensional examples\endsubhead

In the general setting it is easy to construct discretely realizable maps which
are not isotopically realizable.

\example{Example 1.1} Let $S$ be the countable union of $n$-spheres
$S^n_1,S^n_2,\dots$, compactified by a point $p$, and let
$f\:S\sqcup q\to\R^{n+1}$ be a map, throwing the points $p,q$ onto the origin
and each $S^n_k$ homeomorphically onto the standard sphere of radius
$\frac{1}{k}$ centered at the origin.
Clearly, $f$ is realizable but not isotopically.
\endexample

There are also somewhat less straightforward examples.

\example{Example 1.1$'$}
Let $P$ denote the pseudo-arc, $p\:P\sqcup P\to P$ the trivial double cover
and $i\:P\emb\R^2$ any embedding yielded by the Bing definition of the
pseudo-arc \cite{Bi1}, \cite{Le} where all links are round disks in the plane.
Clearly, the composition $i\o p$ is discretely realizable, however in \S2 we
show that it is not isotopically realizable.
\endexample

Perhaps it is worth determining, which compacta admit such natural maps into
Euclidean space, realizable discretely but not isotopically, in particular,
whether the standard embedding of the $p$-adic solenoid into $\R^3$ (cf\.
Example 1.9) precomposed with the trivial double cover is isotopically
realizable.

However, in this paper we treat such cases as pathological, and to eliminate
them we restrict the spaces under consideration in the IR Problem.
From now {\sl we assume the domain $X$ to be a compact $n$-polyhedron and
the target space $Q$ a PL $m$-manifold (without boundary).}
In this setting, the following example was known.

\example{Example 1.2} Let $f\:I\emb\R^3$ be the Wilder arc (i.e\. one of the
two wild arcs shown on Fig\. 1 below) or, more generally, a non-trivial Wilder
arc in the sense of \cite{FH}.
Up to an ambient isotopy, we can assume that $f$ consists of infinitely many
tame knots $f|_{[a_i,a_{i+1}]}\:[a_i,a_{i+1}]\emb\R^2\x[a_i,a_{i+1}]$ (each
of them can be chosen of arbitrary non-trivial isotopy class), where
$a_i=\frac{1}2-\frac{1}{2^i}$, $i=1,2,\dots$, and of a straight line segment
$f|_{[1/2,1]}$.
It was noticed by Keldy\v{s} \cite{Ke4} and Sikkema \cite{Sik} that $f$ cannot
be obtained by a pseudo-isotopy of a tame arc.

For convenience of the reader (since in \cite{Sik} the reduction of Theorem 1
to Theorem 2 was omitted, while the argument in \cite{Ke4} seems to be too
complicated to prove this particular statement), we outline a proof.
Indeed, suppose on the contrary that $g\:I\emb\R^3$ is a PL arc and
$H_t\:\R^3\to\R^3$ a pseudo-isotopy such that $H_0=\id$ and $H_1\o g=f$.
The arc $g$ is the restriction of a PL knot $\bar g\:S^1\emb\R^3$ (it is
supposed that $I\i S^1$), and without loss of generality
$H_t\o\bar g(S^1\but I)$ is sufficiently far from $f([\frac{1}4,\frac{3}4])$
for all $t\in I$ (for $H_t$ and $g$ can be assumed as close as desired to the
identity and to $f$, respectively).
But then for each $n$ there exists an $\eps>0$ such that $H_{1-\eps}\o\bar g$
can be decomposed into at least $n$ knots, which contradicts the uniqueness
of decomposition of $\bar g$ into prime knots (cf\. \cite{Fox}). \qed
\endexample

We call a discretely realizable map $f\:X\to Q$ {\it continuously realizable},
if $\forall\eps>0$ $\exists\dta>0$ such that each embedding $g_\dta\:X\emb Q$,
$\dta$-close to $f$, can be taken onto $f$ by an $\eps$-pseudo-isotopy.
Of course, every continuously realizable map is isotopically realizable, but
not vice versa, as Example 1.2 shows.
We will see below that maps, realizable discretely but not continuously, are
often easier to find and to classify than ones realizable discretely but not
isotopically.
That is why in what follows we keep in mind, along with the IR Problem, the
following {\sl Pre-limit IR Problem:} When does realizability imply continuous
realizability?

\example{Example 1.3} The map $f\:I\sqcup I\to I\lor I\emb\R^3$, whose image
is shown on Fig\. ~1, is not isotopically realizable.
As in 1.2, the proof rests on the Schubert Theorem of uniqueness of
decomposition into prime knots.

\fig 1

The following argument was inspired by an idea due to P. Akhmet'ev.
First let us define an invariant of PL links.
Given a PL link $l\:S^1_1\sqcup S^1_2\to\R^3$ with vanishing linking number,
we consider a decomposition of $l|_{S^1_1}$ into the connected sum of prime
knots $k_1,k_2,\dots,k_p\:S^1_1\to\R^3$.
We call $k_i$ {\it inessential}, if the homotopy class of $l|_{S^1_2}$ in
$\R^3\but l(S^1_1)$, regarded as a conjugate class in
$\pi_1(\R^3\but l(S^1_1))$, lies in the kernel of the homomorphism
$\pi_1(\R^3\but l(S^1_1))\to\pi_1(\R^3\but k_i(S^1_1))$, yielded by
intoduction of the commutativity relations killing the rest prime knots
$k_1,\dots,k_{i-1},k_{i+1},\dots,k_p$.
We define $\alp(l)$ to be the number of essential prime knots among
$k_1,\dots,k_p$.

Now suppose that there exists a (possibly wild) embedding
$g\:I_1\sqcup I_2\emb\R^3$ and a pseudo-isotopy $H_t\:\R^3\to\R^3$ such that
$H_0=\id$ and $H_1\o g=f$.
Extend $g$ by adding two arcs to obtain a link of two (possibly wild)
knots $\bar g\:S^1_1\sqcup S^1_2\emb\R^3$ with zero linking number.
It can be asssumed that $H_t\o\bar g(S^1_i\but I_i)$, $i=1,2$, is sufficiently
far from $f([\frac{1}4,\frac{3}4]\sqcup [\frac{1}4,\frac{3}4])$ for all
$t\in I$.
It is easy to see that for each positive integer $n$ there exists an
$\eps>0$ such that for every PL link $l$, sufficiently close to
$H_{1-\eps}\o\bar g$, the invariant $\alp(l)$ is greater than $n$.
On the other hand, for reasons of compactness (see \S2 for details) $\alp(l)$
is bounded for PL links $l$, sufficiently close to $\bar g$.
\endexample

\remark{Remark} It is not clear whether $\alp(l)$ necessarily stabilizes as
$l\to\bar g$.
This is evidently true if some polyhedral neighborhood of $\bar g(S^1_1)$ in
$\R^3\but\bar g(S^1_2)$ is homeomorphic to the solid torus.
In general it seems (cf\. Example 1.5) that, whatever happens to the quantity
of the essential prime knots, their isotopy types need not stabilize.
\endremark

\example{Example 1.3$'$} If we replace the Wilder arcs in the previous example
by the wild arcs from \cite{Miln, p\. 303}, we will obtain an isotopically
realizable map (compare to \cite{Ke4, Example 1}).
Initial steps of a pseudo-isotopy are indicated on Fig\. 2.
\endexample

\bigskip\bigskip
\fig 2

The above argument in Example 1.3 works equally well for plenty of maps,
similar to the one on Fig\. 1, in particular when one of the two wild arcs on
Fig\. ~1 is replaced by a straight line segment.
One may focus his attention on a class of maps for which it does not work,
namely, the maps $I\sqcup I\to I\lor I\emb\R^3$ whose restriction on each
component is a tame arc.
(For example, the map $I\sqcup I\to\R^3$, obtained from Fig\. ~1 by replacing
each elementary link of two trefoils with the Whitehead link.)
The question, whether these maps are isotopically realizable, seems to be
important, because the contrary would show that, in the range, the
phenomenon of a map, realizable but not isotopically, is not just a
ramification of the phenomenon of a wild embedding, but is somewhat completely
different.
More generally:

\proclaim{Question I} Does there exist a discretely realizable but not
isotopically realizable map which is a locally flat topological immersion?
\footnote{{\it Added in v.2:} A series of such maps $S^n\to\R^{2n}$, $n\ge 3$,
has been constructed \cite{Me2}. However for maps of an $1$-manifold into
$\R^3$ the problem remains open.}
\endproclaim

The positive answer would follow from the positive answer to a general problem
in the link theory \cite{MM}.
Speaking informally, is there a natural theory of `links modulo knots' with a
well-defined operation of connected sum admitting accumulation of complexity
(that is, for some `link modulo knot' $\lambda$ and any $\lambda'$ and any
positive integer $n$ there exists a positive integer $N$ such that for any
$\lambda''$, the connected sum $(\sharp_N\lambda)\sharp\lambda''$ is not
equivalent to $(\sharp_n\lambda)\sharp\lambda'$)?
See \cite{MM} and \cite{MR} for precise statement and some partial results
concerning the latter question, which turns out to be somewhat related to the
long-standing problem of equivalence of the bounded Engel condition and
nilpotence in the class of finitely generated groups.

Another question, arising from the above examples: should a map, realizable
discretely but not isotopically, necessarily be of infinite complexity (in some
sense), or does there exist, say, a PL map such that the better we try to
approximate it, the more `knotted' embedding we should use?
We make this more precise as follows:

\proclaim{Question II} Does there exist a PL map which is PL discretely
realizable but not PL isotopically realizable?
\endproclaim

The definitions of {\it PL (discrete, isotopic, continuous) realizability} can
be obtained by stating the definitions above in the PL category.
Although the answer is unknown in general (see part {\bf B} for the
codimension $\ge 3$ case), we suggest the following negative answer to the
pre-limit version of the latter problem:

\example{Example 1.4}
The PL unknot $f_0\:S^1\emb\R^3$ is not PL continuously realizable, for there
exist PL knots $f_{1/k}\:S^1\emb\R^3$ (see Fig\. 3), arbitrarily close to
$f_0$, which cannot be taken onto $f_0$ by a small PL pseudo-isotopy.
(One can drop `small' by the price of replacing the knots $f_{1/k}$ with the
links $f'_{1/k}$, obtained in the similar way from the Hopf link $f'_0$.)
We prove this in \S2 by showing that $f_0$ is not equivalent to $f_{1/k}$'s by
a small PL (possibly not locally flat) isotopy.
One cannot obtain such example by tying small knots on the image of $f_0$,
since they can be untied by a small PL pseudo-isotopy pushing them to points.
\endexample

\fig 3

\example{Example 1.4$'$}
Alternatively, recall the Hsiang--Shaneson--Wall--Casson--Kirby--Siebenmann
example of PL homeomorphisms of an $n$-torus, $n\ge 5$, arbitrarily close to
the identity (therefore small isotopic to the identity) but not PL isotopic to
it \cite{Ki2, proof of Theorem ~C}, \cite{KS, Appendix ~2 to Essay ~IV}.
(This example was the key ingredient in the elementary disproof of the
Hauptvermutung for manifolds \cite{Sieb1, \S2}, \cite{Ki2, \S0}, \cite{KS}.)
Perhaps the knots $f_{1/k}$ from the previous example can lead to a similar
construction, cf\. \cite{Da, \S12}.
See also \cite{Ed1, end of \S7} and \cite{CS}.
\endexample

Finally we remark that the straightforward way to disprove isotopic
realizability of a continuous map $f\:S^1\sqcup S^1\to\R^3$ is to measure the
way of linking of two simple closed curves by a positive integer $N$, tending
to infinity as they get closer to $f$.
In general, the problem of finding such invariants of a wild link seems to be
intricate and poorly studied (however, see \cite{KY, Part II} and \cite{MR} for
possible sources of such invariants, leaving alone $\alp(l)$ from Example 1.3).
One of the reasons for this difficulty is existence of the Bing sling
\cite{Bi2}, closely related to the knots $f_{1/k}$ from Example 1.4.

\example{Example 1.5} Given a link $l\:S^1\sqcup S^1\emb\R^3$, we define $N(l)$
to be the minimal number of intersections of distinct components under a
null-homotopy of the second component.
Let $h_0$ be the Hopf link and $h_i:S^1_1\sqcup S^1_2\emb\R^3$ be obtained from
$h_{i-1}$ by replacing a small regular neighborhood of $h_{i-1}(S^1_1)$ by the
solid torus containing the knot $f_{1/k(i)}$, where $k$ is some sufficiently
fast growing function, providing $\dist(h_{i-1},h_i)\le\frac{1}{2^i}$.
In other words, $h_{i-1}|_{S^1_1}$ is the axis of $h_i|_{S^1_1}$ in the sense
of \cite{Fox}, meanwhile $h_i|_{S^1_2}=h_{i-1}|_{S^1_2}$.
The limit of $h_i$'s is a wild link $h\:S^1\sqcup S^1\emb\R^3$ (compare to
\cite{Bi2, Fig\. 1}).
In the verification of Example 1.4 in \S2 we show that $N(f'_{1/i})=3$ for
each $i$.
Hence $N(h_i)=3^i$ and $N(h)=\infty$.
\endexample

\subhead B. Isotopic realization in higher codimensions\endsubhead

Whereas every map $f\:X^n\to Q^m$ is discretely realizable (even approximable
by PL embeddings) whenever $m\ge 2n+1$, the `stable range' for isotopic
and continuous realizability is, generally speaking, $m\ge 2n+2$ (this
restriction is sharp by the above examples).
Indeed, sufficiently close PL embeddings $X^n\emb Q^m$, $m\ge 2n+2$, are joined
by a small PL ambient isotopy (cf\. \cite{BK, 5.5}), and the statement follows
(cf\. \cite{Ke3, Lemma 1}).
We shall see that often the restriction $m\ge 2n+2$ can be weakened, especially
for maps satisfying some additional assumptions of `niceness'.

\proclaim{Theorem 1.6} Let $X^n$ be a compact polyhedron, $Q^m$ a PL manifold,
$m-n\ge 3$.

(a) Any PL realizable PL map $f\:X\to Q$ is PL continuously realizable.

(b) If $h\:X\to Y$ is a PL map into a polyhedron $Y$, $i\:Y\emb Q$ is an
embedding and $i\o h$ is realizable, then $i\o h$ is continuously realizable.
\endproclaim

In the part (a), under a stronger assumption $m>\frac{3(n+1)}2$, $Q=\R^m$, a
weaker conclusion of PL isotopic realizability was conjectured in
\cite{ReS2, 1.9d}.
A special case of (b) was proved in \cite{AhRS}: if $X$ is a closed smooth
manifold, $f\:X^n\to\R^{2n-1}$ a generic smooth map,
$i\:\R^{2n-1}\emb\R^{2n}$ the standard inclusion, $i\o f$ is realizable,
then $i\o f$ is isotopically realizabile.
The proof of Theorem 1.6 is based on the results of \cite{Ed1} (see \S3).
We reduce (b) to (a) and prove the latter using slicing techniques (see \S4).
For $(n,m)=(1,3)$ both statements of 1.6 fail by the above examples, and so
does isotopic realizability in (b), but using the proof of (a), it is easy
to verify that PL isotopic realizability in (a) holds in this case.

To approach the general case, where $f$ is an arbitrary continuous mapping,
we introduce $\eps$'s into the Haefliger--Harris theory of isovariant maps.
In the metastable range there is a certain correspondence between embeddings
and isovariant maps, and in \S7 we extend it for discrete and isotopic
realizability:

\proclaim{Criterion 1.7} Let $X^n$ be a compact polyhedron, $Q^m$ a PL
manifold, $f\:X\to Q$ a (PL) map, $m\ge\frac{3(n+1)}2$ in (a)'s and
$m>\frac{3(n+1)}2$ in (b)'s.

(a-) \cite{Harr} $f$ is (PL) homotopic to a (PL) embedding iff
$f^2\:X\x X\to Q\x Q$ is equivariantly homotopic to an isovariant map.

(a) $f$ is (PL) realizable iff $f^2$ is $\eps$-approximable by isovariant maps
for each $\eps>0$.

(a+) Moreover, for each $\eps>0$ there exists $\dta>0$ such that if $f^2$ is
$\dta$-close to an isovariant map, then $f$ is $\eps$-close to a PL embedding.

(b-) \cite{Harr}, \cite{Ed1} (PL) embeddings $g,h\:X\emb Q$ are (PL ambient)
isotopic iff $g^2,h^2$ are isovariantly homotopic. (In the TOP case, no
restrictions of local flatness are imposed on isotopy, in the spirit of
\cite{Miln}.)

(b) $f$ is (PL) isotopically realizable iff there is a homotopy
$\Phi_t\:X\x X\to Q\x Q$ such that $\Phi_1=f^2$ and $\Phi_t$ is isovariant for
$t<1$.

(b+) Moreover, for each $\eps>0$ there exists $\dta>0$ such that if
$g\:X\emb Q$ is a (PL) embedding and $g^2$ is $\dta$-homotopic to $f^2$
by a homotopy $\Phi_t$, isovariant for $t<1$, then $g$ is taken onto $f$ by a
(PL) $\eps$-pseudo-isotopy.
\endproclaim

A map $\Phi\:X\x X\to Q\x Q$ is {\it equivariant} if it commutes with
the involutions $(x,y)\inv (y,x)$ on $X\x X$ and $Q\x Q$, and {\it isovariant}
(cf\. \cite{Hae3}), if in addition $\Phi^{-1}(\Delta_Q)=\Delta_X$, where
$\Delta_X$ means the diagonal of the product $X\x X$.
The `only if' parts are evidently true without any dimensional restrictions.
The TOP case of (b-) follows from its PL case, proved in \cite{Harr}, and an
easy corollary (see 3.5a) of \cite{Ed1, 6.1+8.1}.
See \cite{Hae3} for smooth and \cite{Sk1}, \cite{Sk2} for various deleted
product versions of (a-) and (b-).

For $Q=\R^m$ the PL case of (a) was proved in \cite{ReS2}, and seemingly its
methods suffice to prove for $Q=\R^m$ the statement of (a+), cf\. \cite{ReS2,
pre-limit formulation of 1.2}.
On the other hand, the deleted product theory of \cite{ReS2} does not work
in an arbitrary $Q$, and, which seems to be more important, its natural
generalizations beyond the metastable range, the deleted $n$-th power
obstructions, turn out to be incomplete even in Euclidean space \cite{Sk2}.
That is why we reestablish the result of \cite{ReS2} in the more reliable
setting of isovariant maps.
To prove (b+), whose special case was conjectured in \cite{ReS2, 1.9c}, we
need, besides the straightforward boundary version of (a+), the controlled
version of the classical Concordance Implies Isotopy Theorem, which turns out
to be non-trivial and of independent interest (see part {\bf ~C}).

Since a constant map $X\to Q$ realizes discretely (or isotopically) iff $X$
embeds into $\R^m$, (a) and (b) generalize the case $Q=\R^m$ of (a-).
Consequently by \cite{SSS} (a), (b) are untrue for each $(n,m)$ such that
$3<m<\frac{3(n+1)}2$.
Counterexamples directly to (a), (b) for $(n,m)=(2,4)$ can be deduced from
\cite{AhRS}.

\proclaim{Corollary 1.8} Let $X^n$ be a compact polyhedron
and $Q^m$ a PL manifold.

(a) If $m=2n+1>3$, every continuous map $f\:X\to Q$ is continuously realizable.

(b) If $m>\frac{3(n+1)}2$, discrete realizability implies continuous for a map
$f\:X\to Q$ such that $\Delta(f)=\{(x,y)\in X\x X\mid f(x)=f(y)\}$ has an
equivariant, with respect to the factor exchanging involution, mapping cylinder
neighborhood in $X\x X$.
\endproclaim

The proof is given in \S8.
Under {\it equivariant mapping cylinder neighborhood} of an invariant
subspace $A$ of a space $B$ we mean a closed invariant neighborhood of $A$
in $B$ which is equivariantly homeomorphic to the mapping cylinder
$A\cup_{n\x 0\mapsto g(n)} N\x I$ of some equivariant map $g\:N\to A$.
In particular, the hypothesis of (b) is satisfied if $\Delta(f)$ is an
invariant subpolyhedron of $X\x X$.
Thus we obtain an alternative proof of 1.6b in the metastable range.
Analogously for 1.6a (the required PL version of 1.8b follows from the PL
part of 1.7 analogously to the proof of 1.8b); this yields a different proof
of \cite{ReS2, Conjecture 1.9d}.

\proclaim{Corollary \cite{AM}}
Let $X^n$ be a compact polyhedron, $Q^m$ a PL manifold, $m\ge\frac{3(n+1)}2$,
and $f\:X\to Q$ a discretely realizable map.
The composition of $f$ and the inclusion $Q=Q\x 0\emb Q\x\R$ is isotopically
realizable.
\endproclaim

\example{Example 1.9} Let us construct a map $S^1\x B^2\sqcup B^3\to\R^6$
(in codimension 3), realizable discretely but not isotopically.

Let $T=T_0\supset T_1\supset T_2\supset\dots$ be a sequence of solid tori
$T_i\cong S^1\x B^2$ such that each inclusion $T_i\i T_{i-1}$ induces
multiplication by $3$ in 1-dimensional homology.
The intersection $S=\bigcap T_i$ is the triadic solenoid (cf\. \cite{ES}).
We define a sequence of maps $f_i\:T\to\R^3\but 0$ as follows.
For each $i>0$ let $B^3_i$ be the $(2^{-i})$-neighborhood of the origin $0$
in $\R^3$, and $x_i$ be a point in $B^3_i\but 0$.
Let $f_0$ map $T$ onto $x_1$, and for $i>0$ put $f_i=f_{i-1}$ on $T\but T_i$ and
let $f_i|_{T_i}\:(T_i,\partial T_i)\to(B^2,\partial B^2)\to (B_i^3\but 0,x_i)$
be any map taking $T_{i+1}$ onto $x_{i+1}$ and inducing isomorphism in
relative 2-cohomology.

Then the limit map $f\:T\to\R^3$ meets $0$ in $f(S)$ and has the following
property: the absolute magnitude of the difference $d(\phi,f_0)\in
H^2(T,\partial T;\pi_2(\R^3\but 0))$ is arbitrarily great for every map
$\phi\:(T,\partial T)\to(\R^3\but 0,x_1)$, sufficiently close to $f$.
Indeed, from the fact that each homomorphism in the sequence
$$\dots\to H^2(T,T\but T_2)\to H^2(T,T\but T_1)\to H^2(T,\partial T)$$
is multiplication by $3$ of the group $\Z$ of integers, it follows, firstly,
that $d(f_i,f_0)=1+3+\dots+3^{i-1}=\frac{3^i-1}2$ for $i>0$ and, secondly, that
$d(\phi,\psi)\in 3^i\Z$ for each two maps $\phi,\psi\:(T,\partial T)\to
(\R^3\but 0,x_1)$ agreeing with $f$ on $T\but T_i$.
Given $i>0$, we take $\phi$ so close to $f$ that it can be homotoped, keeping
image in $\R^3\but 0$, to agree with $f$, hence with $f_i$, on $T\but T_i$.
Then $d(\phi,f_0)\in\frac{3^i-1}2+3^i\Z$ and, consequently,
$|d(\phi,f_0)|\ge\frac{3^i-1}2$.
It follows that there does not exist a homotopy $h_t\:T\to\R^3$ such that
$h_1=f$ and $\im h_t\i\R^3\but 0$ for $t<1$.

Now let us use $f$ and the standard inclusion $B^3\emb\R^3$ to obtain a new map
$F\:T\sqcup B^3\to\R^3\x 0\cup 0\x\R^3\emb\R^6$.
It is discretely realizable, for there are embeddings $F_i\:T\sqcup B^3\to\R^6$,
$i=1,2,\dots$, defined by $F_i|_T(p)=(f_i(p),g_i(p))$, where
$g_i\:T\emb B^3_i\i\R^3$ are some embeddings, and by $F_i|_{B^3}=F|_{B^3}$.
On the other hand, $F$ is not isotopically realizable, for otherwise it could
be assumed (see Remark 6.1) that the image of $B^3$ was fixed under the
pseudo-isotopy, and hence there would exist a homotopy $h_t$ as above.

Alternatively, isotopic realizability of $F$ would imply existence of a
homotopy $H_t\:T\x B^3\to\R^6$ such that $H_1(p,q)=F(p)-F(q)$ for each
$(p,q)\in T\x B^3$ and $\im H_t\i\R^6\but 0$ for $t<1$, which can be shown
to be impossible analogously to the above argument. ($H_t$ is
(pseudo-isotopy$|\ind{embedded $T$})\x($pseudo-isotopy$|\ind{embedded $B^3$})$,
composed with the projection $\R^6\x\R^6\to\R^6$ given by $(x,y)\mapsto x-y$.)
\qed
\endexample

\example{Example 1.9$'$} It is easy to see, that if in the above construction
the map $f_i|_{T_i}$ was replaced with one inducing multiplication by $k$ for
each $i$, where $k\not\equiv 1\mod 3$, or, instead, the triadic solenoid was
replaced with the dyadic one, then the resulting map, although isotopically
realizable, would still be not continuously realizable.
\endexample

Actually Example 1.9 can be improved to yield a series of maps $S^n\to\R^{2n}$,
$n\ge 3$, realizable discretely but not isotopically \cite{AM}.
On the other hand, this example, in view of Criterion 1.7, opens up the way to
a complete algebraic description of isotopically and continuously realizable
maps among discretely realizable maps in the metastable range.
Such a description was obtained recently, and for completeness we state it
briefly (for the case $X=S^n$, $Q=\R^m$; the general case is conceptually the
same but involves additional technicalities).

Given a continuous map $f\:S^n\to\R^m$, let us consider open sets
$$U=S^n\x S^n\but\Delta_{S^n}\qquad\supset\qquad U^f=S^n\x S^n\but\Delta(f)
\qquad\supset\qquad U^f_\eps=S^n\x S^n\but P_\eps,$$
where $P_\eps$ is some fixed closed polyhedral neighborhood of $\Delta(f)$,
containing $N_\eps=\{(x,y):||f(x)-f(y)||<\eps\}$ and contained in $N_{2\eps}$.
The composition $\tl f\:U^f\to S^{m-1}$ of the restriction $f^2|_{U^f}$ and the
obvious canonical homotopy equivalence
$$\R^m\x\R^m\but\Delta_{\R^m}\to\R^m\but 0\to S^{m-1}$$
is equivariant with respect to the factor exchanging involution $t$ on
$U^f\i S^n\x S^n$ and the antipodal involution $s$ on $S^{m-1}$.
On the quotient space $U/t$, let us consider the locally constant sheaf
$\Cal Z_m$ with each stalk isomorphic to $\Z$ and the action of $\pi_1(U/t)$ on
the stalks defined by
$$\alp\mapsto\cases s_*,\qquad&\dta_*(\alp)=1;\\0,&\text{otherwise.}\endcases$$
Here $\dta_*\:\pi_1(U/t)\to\Z_2$ denotes the connecting homomorphism from the
exact sequence of the bundle $U\to U/t$ and $1$ denotes the non-trivial element
of $\Z_2$, while $0$ denotes the trivial automorphism and $s_*$ the
automorphism $1\mapsto (-1)^m$ (induced by the involution $s$) of the group
$\Z=\pi_{m-1}(S^{m-1})$.
For each $X\i U$ we write $H^{m-1}\eq(X)$ for the cohomology group
$H^{m-1}(X/t;\Cal Z_m|_{X/t})$.
Let $®(f)\in H\eq^{m-1}(U^f)$ denote the first obstruction for equivariant
homotopy (cf\. \cite{CF, \S\S2,4}) of the maps $\tl f$ and $\tl i|_{U^f}$,
where $i\:S^n\emb\R^m$ is the standard (or any, in view of \cite{Ze})
inclusion.
Similarly, we denote by $o_\eps(f)\in H\eq^{m-1}(U^f_\eps)$ the first
obstruction for equivariant homotopy of the restrictions $\tl f|_{U^f_\eps}$,
$\tl i|_{U^f_\eps}$.
The latter obstruction can be equivalently defined in the spirit of van Kampen
(cf\. \cite{AM}, compare to \cite{ReS2, 1.4}).
Finally, for embeddings $g_1,g_2\:S^n\emb\R^m$, $\eps$-close to $f$, the
first obstruction $d(g_1,g_2)$ for equivariant $\eps$-homotopy of
$(\tl g_i|_{U^f})$'s is an element of $H^{m-1}\eq(U,U^f_\eps)$.
Let us write $G_k=H^{m-1}\eq(U,U^f_{2^{-k}})$ and let $j^l_k\:G_l\to G_k$,
$l>k$, be the forgetful homomorphism.

The following result, whose proof is based on Criterion 1.7,
shows in particular that from the algebraic viewpoint, maps yielding
negative solution to the metastable case of the IR Problem look quite similar
to phantom maps \cite{Gr}.

\proclaim{Theorem \cite{AM}}\footnote{{\it Added in v.2:}
In the case $m<2n$, the proof of this theorem in \cite{AM} contains a mistake
(on page 81, line 8), and it works to prove the theorem only under the
additional assumption that the given map $f$ is discretely $k$-realizable for
each $k=m+1,m+2,\dots,2n$ (see definition below).
Without this assumption, the parts (a) and (b) are incorrect already for
$m=2n-1$, while (d) fails for $m=2n-5\ge 9$ \cite{Me2}.
The part (c) (and hence its corollary below) is correct as stated, but its
proof in \cite{AM} is insufficient if without the additional assumption; the
correct proof rests on higher cohomology operations and appears in \cite{Me2}.
Let us state the required definition of discrete $k$-realizability.
Assume $2n>m>\frac{3(n+1)}2$ and $m\le k\le 2n$.
Let us call a PL map $f\:S^n\to\R^m$ a {\it $k$-embedding} if there is
a triangulation $T$ of $S^n$ such that $f$ is simplicial in some subdivision of
$T$ and embeds each simplex of $T$, and
$f(\sigma)\cap f(\tau)=f(\sigma\cap\tau)$ for any two simplices $\sigma^s$,
$\tau^t$ of $T$ such that $s+t\le k$.
Let us call a map $f\:S^n\to\R^m$ {\it discretely $k$-realizable}, if
$\forall\varepsilon>0$ $\exists\delta>0$ such that any $(k-1)$-embedding,
$\delta$-close to $f$, is PL $\varepsilon$-homotopic in the class of
$(k-2)$-embeddings to some $k$-embedding.
It is easy to see that $f$ is discretely $m$-realizable iff
$o_\eps(f)=0$ for all $\eps>0$.
}
Let $f\:S^n\to\R^m$ be a continuous map,
$m>\frac{3(n+1)}2$.

(a) $f$ is discretely realizable iff $o_\eps(f)=0$ for
$\eps=\frac{1}2,\frac{1}4,\frac{1}8,\dots$

(b) $f$ is isotopically realizable iff $o(f)=0$.

(c) Suppose that $f$ is discretely realizable.
$f$ is continuously realizable iff the inverse spectrum $\{G_k; j^l_k\}$
satisfies the Mittag-Leffler condition or, equivalently \cite{Gr},
$$\derlim\{G_k; j^l_k\}=0.$$

(d) Suppose that $f$ is discretely realizable.
$f$ is isotopically realizable iff $$0=O(f)\in\derlim\{G_k; j^l_k\}.$$
\endproclaim

The obstruction $O(f)$ can be defined (cf\. \cite{Me2}) as the class of the
sequence $d(g_1,g_2),d(g_2,g_3),\dots$, where $g_k\:S^n\emb\R^m$ is an
embedding, $(2^{-k})$-close to $f$.
See \cite{Ma} for definition and basic properties of the derived limit functor.
The fact that no obstructions arise in dimensions other than $m-1$ is due to
the Serre Theorem on finiteness of homotopy groups of spheres.
Using the fact that the forgetful homomorphism $H^*\eq(\cdot)\to H^*(\cdot)$
factors through the multiplication by $2$ in $H^*(\cdot)$, and the Alexander
duality, one immediately obtains the following

\proclaim{Corollary \cite{AM}} Let $f\:S^n\to\R^m$ be a discretely realizable
map, $m>\frac{3(n+1)}2$.
If the canonical epimorphism $H_{2n-m}(\Delta(f))\to\check H_{2n-m}(\Delta(f))$
between the reduced Steenrod (exact; cf\. \cite{Ma}) and the reduced \v{C}ech
(continuous; cf\. \cite{ES}) homology has trivial kernel, then $f$ is
continuously realizable.
\endproclaim

This puts the IR Problem in the metastable range in the context of the
discussion `continuity versus exactness' in Eilenberg and Steenrod
\cite{ES, p\. 265} (see \cite{Fe} for a modern version).
It follows e.g., that if a map $f\:S^n\to\R^{2n}$, $n\ge 4$, is realizable
discretely but not continuously, then the compactum $\Delta(f)$ cannot be
zero-dimensional or have countable Steenrod $0$-homology (cf\. \cite{Harl}).

\remark{Remark} In proving Criterion 1.7 we obtain a number of interesting
results in the PL category in the metastable range, among which are: sufficient
conditions for existence of an embedding in the $\eps$-homotopy class of a map
(7.1), of an embedding in the $\eps$-regular homotopy class of an immersion
(7.2), of an immersion in the $\eps$-homotopy class of a map (7.4),
of an $\eps$-ambient isotopy between close embeddings (7.9a).
These all are controlled versions of Harris' criteria \cite{Harr},
however we use some ideas, additional (proof of 7.2) and alternative
(proof of 7.4) to that of \cite{Harr}.
In fact, our proof of 7.4 is a new geometric proof of \cite{Harr, Th\. 2}
(roughly a half of 1.7a-) and a good candidate for generalization for
$k$-tuples of points.
\endremark

\proclaim{Conjecture 1.10} In each range $m\ge\frac{(k+1)(n+1)}k$
$(m>\frac{(k+1)(n+1)}k)$ the analogue of Criterion 1.7 holds for isovariant
maps $(X^n)^k\to (Q^m)^k$, provided $m-n\ge 3$.
\endproclaim

We call a map $\Phi\:(X^n)^k\to (Q^m)^k$ {\it isovariant} if it commutes
with the actions of the symmetric group $S_k$ on $(X^n)^k$ and $(Q^m)^k$, and
if $\Phi^{-1}(\Delta^S_Q)=\Delta^S_X$ for each $S\i\{1,2,\dots,k\}$, where
$\Delta^S_X=\{(x_1,\dots,x_k)\in (X^n)^k\mid i,j\in S\Rightarrow x_i=x_j\}$.
Such a result cannot be expected in codimension 2: the reader may wish to
verify that non-triviality of the link $f'_{1/k}$ from Example 1.5 (as well
as that of the $k$-th Milnor's link \cite{Miln} and of the link Whitehead$_k$,
cf\. \cite{KY}) is not detected by the isovariant homotopy class of the
$(k+2)$-th power mapping.
It is worth observing that, in contrast to the metastable range (where smooth
embeddability and quasi-embeddability are equivalent to PL embeddability), the
smooth and deleted versions of 1.10 are untrue \cite{Hae2}, \cite{Sk2}.

\subhead C. Other definitions of realizability and relations on close
embeddings\endsubhead

A (PL) map $F\:X\x I\to Q\x I$ will be called a {\it $(PL)$ pseudo-concordance}
if $F^{-1}(Q\x 1)=X\x 1$ and $F|_{X\x [0,1)}$ is an embedding.
We call a (PL) map $f\:X\to Q$ {\it $(PL)$ concordantly realizable}
if $f\x\id_1$ extends to a (PL) pseudo-concordance $F\:X\x I\to Q\x I$.
If, in addition, $F^{-1}(Q\x 0)=X\x 0$, then $f$ is
{\it $(PL)$ pseudo-concordant to} the unique (PL) embedding $g\:X\emb Q$ such
that $g\x\id_0=F|_{X\x 0}$.

\example{Example 1.11} The map $f$ from Example 1.3 is concordantly realizable.
Several slices $(\im F)\cap\R^3\x t$ of a pseudo-concordance $F$ are shown on
Fig\. 4.
\endexample

\bigskip
\fig 4
\bigskip\bigskip

In higher codimensions the situation is again quite different:

\proclaim{Theorem 1.12} Let $X^n$ be a compact polyhedron, $Q^m$ a PL manifold,
$m-n\ge 3$.
Then each (PL) concordantly realizable (PL) map $f\:X\to Q$ is
(PL) isotopically realizable.
Moreover, $\forall\eps>0$ $\exists\dta>0$ such that every (PL) embedding
$g\:X\emb Q$, (PL) $\dta$-pseudo-con\-cordant to $f$, can be taken onto $f$ by
a (PL) $\eps$-pseudo-isotopy.
\endproclaim

The proof is given in \S6; the TOP case is based on the following controlled
version 1.13a of the classical PL Concordance Implies Isotopy Theorem (CIIT):

\proclaim{Theorem 1.13}
For each $\eps>0$ and a positive integer $n$ there exists $\dta=\dta(n,\eps)>0$
such that the following holds.

(a) Let $X^n$ be a compact polyhedron and $Q^m$ a PL manifold, $m-n\ge 3$.
Then each two PL $\dta$-concordant embeddings $f,g\:X\emb Q$ are PL
$\eps$-ambient isotopic.

(b) Let $X^n$ be a compact smooth manifold, $Q^m$ a smooth manifold,
$m>\frac{3(n+1)}2$.
Then each two smoothly $\dta$-concordant embeddings $f,g\:X\emb Q$ are
smoothly $\eps$-ambient isotopic.
\endproclaim

It seems that Hudson's original proof of CIIT \cite{Hu} (as well as Lickorish's
proof of the case $Q=S^m$ \cite{Li, Theorem 6}) does not work to prove 1.13a
(compare to remarks in \cite{Mill1, Introduction}, \cite{ReS2, \S2}).
In \cite{Ro} Rourke sketched a new proof of CIIT, and in \cite{KeL, last
paragraph} it was `expected that, when the details of Rourke's proof are
published, they will apply' to prove 1.13a.
A special case of 1.13a was conjectured in \cite{ReS2, 1.9a}.
Perhaps 1.13a can be also proved by the methods of \cite{Ed1, proof of 7.1},
but hardly by that of \cite{Co, proof of Lemma 1}.

Be that as it may, in \S5 we present an explicit proof of 1.13a.
It is far from being a trivial extension of either known proof of CIIT (this is
clear at once from the statement of Lemma 5.8), and it is also a new proof of
CIIT (since Theorem 1.13a generalizes CIIT, by taking a metric on $Q$ with all
distances $<\dta$).
Theorem 1.13a, along with CIIT, is untrue in codimension 2 because of slice
knots and links.
It is worth observing that in our proof of 1.13a the main efforts are
applied to obtain $\eps$-ambient isotopy, rather than ambient $\eps$-isotopy.
In the proof of 1.13a we use Theorem 3.3a, which includes Miller's controlled
version \cite{Mill1, Theorem ~9} of Zeeman's Unknotting Balls \cite{Ze}.
In turn, the $\partial Q\neq\emptyset$ version of 1.13a, which is proved
analogously, immediately implies \cite{Mill1, Theorem ~9} (this was pointed out
in \cite{KeL, last paragraph}) and \cite{Co, Lemma ~1}.
In \cite{BS1, Idea of proof of Theorem ~3} a statement, similar to 1.13a was
used with `reference' to Hudson's CIIT (see \cite{Hu}); the misquotation
disappears in the revised proof \cite{BS2}.

Next we convert 1.13a to the smooth category and obtain Theorem 1.13b.
It answers, at least to some extent, a question of Kirby \cite{Ki1, discussion
preceding 2.1}: `suppose $m>\frac{3(n+1)}2$, is there a function $\eps$ of
$\dta$ such that any smoothly $\dta$-isotopic smooth embeddings are smoothly
$\eps$-ambient isotopic?'
(From the proof of 1.13b it follows that $\eps$ can be taken as
$c(n)*\dta$, where $c(n)$ is a constant depending on $n=\dim X$.)
It should be mentioned that our proof of 1.13b uses Kirby's partial answer
to his question (see Theorem 3.3b).
We conjecture that 1.13b holds in codimension $\ge 3$.

A corollary of 1.13b is the smooth version 3.2b+ of Edwards' Theorem 3.2a on
$\eps$-equivalence of PL embeddings, close to a TOP embedding (see also 3.5b+).
Using 1.13a we also obtain an alternative controlled version 3.7+ of CIIT.

\example{Example 1.14} In general, small (smooth or TOP/PL locally flat)
isotopy, in particular small concordance, does not imply small ambient isotopy.
(Of course, it implies a great smooth or TOP/PL ambient isotopy \cite{Hi},
\cite{EK}, \cite{RoS}.)
Indeed, take the standard circle $S^1\i\R^3$ and tie, near a point $x\in S^1$,
a small (e.g\. trefoil) knot on it to obtain an embedding $f_0\:S^1\emb\R^3$.
One can shift this small knot along $S^1$ by a (smooth or TOP/PL locally flat)
isotopy $f_t$, which, at each moment $t\in I$, has support in a small
neighborhood of the current position of the small knot on the circle.
(Such an isotopy cannot be obtained by means of rotation of the whole circle.)
But it is clear (see \S2 for details) that $f_0$ and $f_1$ can not be joined
by a small ambient isotopy.
\endexample

\fig 5

Finally, we relate the \v{S}\v{c}epin--\v{S}tan'ko definition of isotopic
realizability to the two Akhmet'ev's definitions \cite{Akh1}.
A map $f\:X\to Q$ of a compact smooth manifold into a smooth manifold is
{\it $A_1$-($A_2$-)isotopically realizable,} if there is a homotopy
$f_t\:X\to Q$ (resp\. $H_t\:Q\to Q$), called an
{\it $A_1$-($A_2$-)pseudo-isotopy,} such that $f_1=f$ (resp\. $H_0=\id$ and
$H_1\o g=f$ for some smooth embedding $g\:X\emb Q$), and which is a smooth
isotopy (resp. smooth ambient isotopy) for $t<1$.
Certainly, $A_i$-isotopic realizability is not equivalent to isotopic
realizability, see Example 1.2.
Evidently, $A_2$-isotopic realizability implies $A_1$-.
But the author does not see why the reverse implication holds, as claimed in
\cite{Akh1}.

\example{Example 1.15} The standard embedding $f\:S^1\emb\R^3$ is, of course,
isotopically realizable in either sense.
However, there is an $A_1$-pseudo-isotopy $f_t$ from $f$ to an emebedding,
which cannot be covered by a pseudo-isotopy (in particular, by an
$A_2$-pseudo-isotopy).

Indeed, rotate a small knot around the circle, as in Example 1.14, so that
its size tends to zero (hence $f_t\to f$) as $t\to 1$ and so that the speed
of its rotation, along with the number of turns, tends to infinity as $t\to 1$.
If $f_t$ was covered by a pseudo-isotopy $H_t\:Q\to Q$, we would obtain a
contradiction with 1.14. \qed
\endexample

In order to avoid too restrictive assumptions of smoothness, let us say that
a map $f\:X\to Q$ is {\it $M$-isotopically realizable}, if there exists a
homotopy $f_t\:X\to Q$, called an {\it $M$-pseudo-isotopy}, such that $f_1=f$
and for each $t<1$ the map $f_t$ is a topological embedding.
The letter `$M$' accounts for the fact that $f_t$ for $t\in [0,1)$ is an
isotopy in the sense of \cite{Miln}.

\proclaim{Question III} Does there exist an $M$-isotopically realizable map
which is not isotopically realizable?
\endproclaim

By 1.12 such a map cannot be found in the codimension $\ge 3$ range.
Furthermore, the DIFF case of the following theorem, proved in \S6, implies
(in view of 1.12) that in the metastable range all 4 definitions of
pseudo-isotopy, as well as all 4 definitions of isotopic realizability
($A_1$ and $A_2$ of Akhmet'ev, $M$ in the spirit of Milnor, and the classical
one of \v{S}\v{c}epin--\v{S}tan'ko) are equivalent.
For another application of Theorems 1.12 and 1.16 see Remark 6.1.

\proclaim{Theorem 1.16} Let $X^n$ be a compact polyhedron (compact smooth
manifold), $Q^m$ a PL (smooth) manifold, $f\:X\to Q$ a continuous map,
$m-n\ge 3$ (respectively $m\ge\frac{3(n+1)}2$, $n>1$ in (a), $m>\frac{3(n+1)}2$
in (b)).

(a) If $f$ is isotopically realizable, then there exists a pseudo-isotopy,
taking a PL (smooth) embedding $g\:X\emb Q$ onto $f$.

(b) If a PL (smooth) embedding $g\:X\emb Q$ is taken onto $f$ by a
pseudo-isotopy, then $g$ can be taken onto $f$ by a pseudo-isotopy
$H_t\:Q\to Q$ such that whenever $t\in [0,1)$, $H_t$ is a PL (smooth) isotopy.
\endproclaim

A continuous map $f\:X_1\sqcup\dots\sqcup X_k\to Q$ (where the components
$X_1,\dots,X_k$ are fixed, but not necessarily connected) is called
{\it disjoinable}, if it is approximable by link maps (cf\. \cite{ST},
\cite{DRS}, \cite{Sp}); a map $g\:X_1\sqcup\dots\sqcup X_k\to Q$ is called a
{\it (generalized) link map} if $g(X_i)\cap g(X_j)=\emptyset$ whenever
$i\neq j$ (cf\. \cite{Me1}).
We call $f$ {\it homotopically disjoinable} if there is a homotopy $f_t$ such
that $f_1=f$ and $f_t$ is a link map for $t<1$.

\example{Example 1.17} (i) The proof of Example 1.1$'$ allows to replace
`(isotopically) realizable' with `(homotopically) disjoinable' in its statement.

(ii) The construction of Example 1.9 yields a map $S^1\x B^2\sqcup pt\to\R^3$,
which is disjoinable but not homotopically disjoinable.

(iii) The map $f$ from Example 1.3 turns out to be homotopically disjoinable.
Indeed, we start from two disjoint arcs.
In the spirit of Example 1.11 we generate linking trefoils, keeping ends of
arcs fixed, by the price of self-intersections of components.
In the spirit of Example 1.3$'$ we compensate the increase of the linking
number by small loops tending to the singular point as the time approaches $1$.
\endexample

We conjecture that the analogues of 1.6--1.8 for homotopic disjoinability hold
and can be proved analogously.
Moreover, the remark (iii) in conjuction with the facts that for classical
links singular link concordance implies link homotopy and that
$\kappa$-invariant has trivial kernel up to link homotopy motivate a conjecture
that every map $S^1\sqcup\dots\sqcup S^1\to\R^3$ is homotopically disjoinable.

\remark{Remark} For completeness let us consider a concept, approximately dual
to homotopic disjoinability, that is, admitting distant self-intersections and
prohibiting close ones.
We call a map $f\:X\to Q$ {\it locally isotopically realizable} if there is a
homotopy $f_t\:X\to Q$ such that $f_1=f$ and $f_t$ is a topological immersion
for $t<1$.
From the $C^0$-dense $h$-principle for smooth immersions it follows
\cite{Akh1, proof of Lemma~2} that if a compact smooth manifold $X^n$ smoothly
immerses into a smooth manifold $Q^m$, $m-n\ge 1$, then each map $f\:X\to Q$
is locally isotopically realizable.
\endremark

\subhead Acknowledgements \endsubhead

I wish to thank P. ~M. ~Akhmet'ev, E. ~V. ~\v{S}\v{c}epin and A. ~B. ~Skopenkov
for suggesting the problem, attention, and help, which made the present work
possible.
I am especially indebted to Petya Akhmet'ev for numerous discussions which
eventually led to the present undersanding of the phenomena under study, and to
N. ~B. ~Brodskij, A. ~N. ~Dranishnikov and R. ~V. ~Mikhailov for pointing
out some gaps in my arguments.
I would like to acknowledge D. ~A. ~Botin, A. ~V. ~\v{C}ernavskij,
I. ~A. ~Chubarov, U. ~Koschorke, W. ~Lewis, D. ~V. ~Million\v{s}\v{c}ikov,
V. ~M. ~Nezhinskij, D. ~Repov\v{s}, R. ~R. ~Sadykov, K. ~R. ~Salikhov,
M. ~A. ~\v{S}tan'ko, E. ~G. ~Skl'arenko, Yu. ~A. ~Turygin and my undergraduate
advisor Yu\. ~P. ~Solovjov for useful discussions and important remarks.

\head 2. Verification of examples \endhead

\demo{Verification of 1.1$'$}
The argument below was inspired by an idea of A. Skopenkov
(compare \cite{RSS2} and \cite{ReS2}).
Suppose that $f_t\:P\x\{0,1\}\to\R^2$ is a homotopy such that $f_1=i\o p$ and
$f_t(P\x 0)\cap f_t(P\x 1)=\emptyset$ whenever $t<1$.
Let $F_t\:P\x P\to\R^2\x\R^2$ be the map defined by
$(p,q)\mapsto (f(p\x 0),f(q\x 1))$.
(In other words, $F_t=f_t^2|_{(P\x 0)\x (P\x 1)}$.)

Then $F_t^{-1}(\Delta_{\R^2})$ is empty for $t<1$ and equals
$\Delta_P=\{(p,p)\in P\x P\}$ for $t=1$.

Since $P$ is acyclic, $P\x P$ is acyclic and the map
$F_0\:P\x P\to\R^2\x\R^2\but\Delta_{\R^2}$ is null-homotopic in
$\R^2\x\R^2\but\Delta_{\R^2}$.
Therefore the map
$$F_1|_{P\x P\but\Delta_P}\:P\x P\but\Delta_P\to\R^2\x\R^2\but\Delta_{\R^2}$$
is also null-homotopic in $\R^2\x\R^2\but\Delta_{\R^2}$.
The latter map is equivariant with respect to the involutions $(p,q)\inv (q,p)$
on $P\x P\but\Delta_P$ and $\R^2\x\R^2\but\Delta_{\R^2}$, and the latter
space is equivariant homotopy equivalent to $S^1$ equipped with the antipodal
involution.
Thus we obtain an inessential equivariant map $P\x P\but\Delta_P\to S^1$.

By \cite{RSS1}, existence of such a map implies that $P\x P\but\Delta_P$ is not
connected.
Suppose that $(p_1,p_2)$ and $(q_1,q_2)$ lie in distinct connected components
of $P\x P\but\Delta_P$; without loss of generality $p_2\neq q_1$.
Then either $(p_1,p_2)$, $(q_1,p_2)$ or $(q_1,p_2)$, $(q_1,q_2)$ lie
in distinct components, say the first ones.
Consequently, $p_1$ and $q_1$ lie in distinct connected components of
$P\but p_2$.
But $P$ has no separating points, and we arrive at a contradiction. \qed
\enddemo

\demo{Verification of 1.3} We are to prove that for every wild link
$\bar g\:S^1_1\sqcup S^1_2\emb\R^3$
$$\limsup\Sb l\to\bar g,\\ l\in\Cal L\PL\endSb\ \alp(l)\ <\infty,$$
where $\Cal L\PL$ denotes the subspace of PL embeddings in the
space of all continuous maps $S^1\sqcup S^1\to\R^3$, equipped with the topology
of uniform convergence.

Assume on the contrary that there is a sequence of PL links $l_1,l_2,\dots$,
converging to $\bar g$ and such that $\alp(l_i)\to\infty$.
Let $N$ be a polyhedral neighborhood of $\bar g(S^1_1)$ in
$\R^3\but\bar g(S^1_2)$.
We can assume that $l_i(S^1_1)\i N$ for each $i$.
We fix a decomposition of $N$ into handles: $N=B^3\cup H_1\cup\dots\cup H_q$,
where $H_j\cong D^2\x I$ via a homeomorphism $h_j\:D^2\x I\to H_j$ and
$H_j\cap B^3=h_j(D^2\x\partial I)$ (without loss of generality there are no
$2$-handles).
By the definition of a prime knot, there is a collection of disjoint $3$-balls
$B_{i,1},B_{i,2},\dots,B_{i,\alp(l_i)}$ such that the boundary of each $B_{ij}$
meets $l_i(S^1_1)$ precisely in two points, separating one essential prime knot
$k_{i,e_j}\:(I,\partial I)\to(B_{ij},\partial B_{ij})$ from all other prime
knots in the decomposition of $l_i|_{S^1_1}$.
Since $k_{i,e_j}$ is essential, its image meets each $D^2$-fiber of some handle
$H_{m_{ij}}$.
By our assumption $\alp(l_i)\to\infty$ as $i\to\infty$, hence for some handle
$H_l$ the number $n_i$ of the knots $k_{i,e_j}$ such that $m_{ij}=l$ tends to
infinity as $i\to\infty$.
Let us fix the handle $H_l=h_l(D^2\x I)$ and denote the images of the latter
knots by $\kappa_{i,1},\kappa_{i,2},\dots,\kappa_{i,n_i}$.

Now for each $i,j$ the intersection $\kappa_{ij}\cap H_l$ is the union of some
PL arcs.
Each of these arcs meets either $h_l(D^2\x 0)$ or $h_l(D^2\x 1)$, and at least
one of these arcs, denoted by $a_{ij}$, meets $h_l(D^2\x\frac{1}2)$.
The diameter of $a_{ij}$ is therefore at least
$d=\dist(h_l(D^2\x\frac{1}2),h_l(D^2\x\{0,1\})$.
Since the arcs $a_{i,1},\dots,a_{i,n_i}$ are contained in disjoint curves
$\kappa_{i,1},\dots,\kappa_{i,n_i}\i l_i(S^1_1)$, the PL curve $l_i(S^1_1)$
contains at least $n_i$ disjoint subarcs, each of diameter at least $d$.
Since $n_i$ tends to infinity as $i\to\infty$, this is in contradiction with
the assumption of convergence of $l_i$'s to $\bar g$. \qed
\enddemo

\demo{Verification of 1.5}
Let $T$ be a small regular neighborhood of $f_0(S^1)$ and $l$ be a circle in
the complement of $T$, linked with $f_0(S^1)$ with linking number one.
If $H_t\:\R^3\to\R^3$ is a PL pseudo-isotopy taking a PL embedding $f_\eps$
onto $f_0$, then $h_t=H_t\o f_\eps$ is a PL (possibly not locally flat) isotopy
(with all points of failure of local flatness occuring for $t=1$), and moreover
if $H_t$ is sufficiently small, the image of $h_t$ lies in $T$.
The statement of Example 1.5 follows from Claims 2.1, 2.2 below.
\enddemo

\proclaim{Claim 2.1} If $g\:S^1\emb T$ is a PL embedding, the minimal number
$I(g)$ of transversal intersections of a singular disk, spanned by $l$, with
$g(S^1)$ is invariant under PL (generally not locally flat) isotopy in $T$.
\endproclaim

\demo{Proof}
Let $h_t\:S^1\emb T$ be a PL (possibly not locally flat) isotopy and $D$ be
a disk, spanned by $l$ and meeting $h_0(S^1)$ in $I(h_0)$ points.
It suffices to show that there is a disk $D'$, spanned by $l$ and meeting
$h_1(l)$ in $I(h_0)$ points.
Without loss of generality we can assume that $h_t$ is either locally flat,
or locally knotted at a unique point $a\in S^1$ in the moment $t=1/2$,
so that $h_t=h_0$ outside a small neighborhood $U$ of $a$ and $h_t(U)\i W$,
where $W$ is a regular neighborhood of $h_{1/2}(S^1)$ relative
$h_{1/2}(S^1\but U)$.

In the first case $h_t$ can be covered by an ambient isotopy $H_t$ \cite{RoS},
which carries the disk $D$ so that the number of intersections of
$D_t=H_t(D)$ with $h_t(S^1)=H_t(h_0(S^1))$ remains constant.
In the second case we modify the disk $D$ as follows.
First we shift any intersections with $h_0(S^1)$ along $h_0(S^1)$ out of
$h_0(U)$.
Then the shifted $D$, denoted by $\tl D$, does not meet $h_0(S^1)$ in $W$.
Next we push $\tl D$ out of $W$.
It is possible since the kernel of
$$incl_*\:\pi_1(\R^3\but(h_0(S^1)\cup W))\to\pi_1(\R^3\but h_0(S^1))$$
is trivial.
Indeed, introduction of commutativity relations into the subgroup of
$\pi_1(\R^3\but(h_0(S^1)))$, consisting of conjugates to the loops lying
in $W$, yields
$$comm_*\:\pi_1(\R^3\but(h_0(S^1)))\to\pi_1(\R^3\but(h_0(S^1)\cup W))$$
such that $incl_*\o comm_*=\id$.
Hence $incl_*$ has trivial kernel, consequently we can replace $\tl D$ by a
disk $D'$ avoiding $W$.
Since $h_t$ has its support in $W$, the disk $D'$ meets $h_0(S^1)$ and
$h_1(S^1)$ in the same points. \qed
\enddemo

Now let us recall the knots $f_{1/k}$ from Example 1.5.

\proclaim{Claim 2.2} $I(f_{1/k})=3$ for each $k=1,2,\dots$ (while $I(f_0)$ is
clearly $1$).
\endproclaim

\demo{Proof}
It is clear that $I(f_{1/k})\le 3$ for each $k=1,2\dots$.
There is a $k$-fold cover $p\:T\to T$ such that $p^{-1}(f_1(S^1))=f_k(S^1)$,
hence $I(f_{1/k})\ge I(f_1)$.
It remains to show that $I(f_1)$ is not less than 3.
Actually $f_1$ is the trefoil knot, and $l$ represents $a^{-1}b^2$ in its group
$G=\langle a,b\mid aba=bab\rangle$.
A disk, spanned by $l$, cannot meet $f_1(S^1)$ in 2 points, for this would
imply an even linking number of $l$ and $f^1(S^1)$.

Suppose that there is a disk, spanned by $l$ and meeting $f_1(S^1)$
transversely in one point.
Then $l$ is homotopic to a loop representing an element of $G$ of type
$g^{-1}b^\eps g$, where $g\in G$, $\eps=1\text{ or }-1$.
Since $[[l]]=[a^{-1}b^2]=[b]$ in $G/[G,G]=H_1(\R^3\but f(S^1))$, necessarily
$\eps=1$, hence for some $g\in G$ the equality $a^{-1}b^2=b^g$ holds in $G$
(here $b^g$ denotes $g^{-1}bg$).
We will show this to be impossible by considering a representation of $G$.

It is easy to see that the formulae $a\mapsto (123)$, $b\mapsto (432)$
yield a representation $\phi\:G\to A_4\i S_4$ in the symmetric group
(the well-known representation $\psi$ in $S_3$ is insufficient, since
$\psi(a^{-1}b^2)=\psi(a^{-1})=\psi(b^{ba})$).
We have $\phi(a^{-1}b^2)=(412)$ and $\phi(b^g)=(432)^{\phi(g)}$.
But $(432)$ and $(412)=(432)^{(13)}$ are not conjugate in $A_4$,
which is a contradiction. \qed
\enddemo

\demo{Verification of 1.14}
Suppose that the small isotopy $f_t$ can be covered by a small ambient isotopy
$H_t\:\R^3\to\R^3$, $H_0=\id$, $H_1\o f_0=f_1$ (we omit the epsilonics).
Denote by $\pi$ the fundamental group $\pi_1(\R^3\but f_0(S^1))$.
Let $a\in\pi$ be the class of a small circle around $S^1$ far from $x$.
Let $b$ be any element of $\pi$ which is not a power of $a$, and represent $b$
by a small loop $l$ (which necessarily lies near $x$).
Then $H_1(l)$ lies possibly little farther from $x$, but still near it.
Now $f_1$ has its small knot far from $x$, hence far from $l$.
This means that $l$ should represent a power of $a$ in
$\pi=\pi_1(\R^3\but f_1(S^1))$, which is a contradiction. \qed
\enddemo

\head 3. Some facts on close PL, DIFF and TOP embeddings \endhead

In this section we recall some approximation theorems to be heavily used in
the rest of the paper.
Exceptions are 3.2b+, 3.5b+, 3.7+, which are not used in the sequel; on the
contrary, their proofs require Theorem 1.13, proved in \S5.

\proclaim{Theorem 3.1} (a) \cite{Ch4}, \cite{Mill2}, \cite{Br1}, \cite{Ed1, 8.1},
\cite{Br2}
Let $(X^n,Y^{n-1})$ be a polyhedral pair, $(Q^m,\partial Q)$ a PL manifold,
$m-n\ge 3$.
Then any TOP embedding $f\:(X,Y)\emb (Q,\partial Q)$ is $\eps$-approximable,
for each $\eps\:X\to(0,\infty)$, by a PL embedding $g\:(X,Y)\emb
(Q,\partial Q)$.
Moreover if $Z$ is a subpolyhedron of $X$ and $f|_Z$ is PL, then it can
be assumed that $g|_Z=f|_Z$.

(b) \cite{Hae1}, \cite{Ki1, 2.2}
Let $X^n$ be a compact smooth manifold, $Q^m$ a smooth manifold and
$m\ge\frac{3(n+1)}2$.
Any TOP embedding $f\:X\emb Q$ is $\eps$-approximable, for each $\eps>0$,
by a smooth embedding $g\:X\emb Q$.
Moreover if $Z$ is a closed subset of $X$ and $f$ is smooth on the
$\dta$-neighborhood of $Z$, then it can be assumed that $g|_Z=f|_Z$.
\endproclaim

\proclaim{Theorem 3.2} (a) \cite{BD}, \cite{BS2}, \cite{Mill2}, \cite{Ed1, 6.1}
Let $(X^n, Y^{n-1})$ be a polyhedral pair, $(Q^m,\partial Q)$ a PL manifold,
$m-n\ge 3$, $Z$ a subpolyhedron of $X$, and $f\:(X,Y)\emb (Q,\partial Q)$ a
TOP embedding.
For each $\eps\:Q\to(0,\infty)$ there exists $\dta\:X\to(0,\infty)$ such that
any PL embeddings $g,h\:(X,Y)\emb (Q,\partial Q)$, $\dta$-close to $f$, are PL
$\eps$-ambient isotopic.
Moreover if $g|_Z=h|_Z$, then the isotopy can be chosen fixing $g(Z)=h(Z)$.

(b)
Let $X^n$ be a compact smooth manifold, $Q^m$ a smooth manifold,
$m>\frac{3(n+1)}2$, $Z$ a closed subset of $X$, and $f\:X\emb Q$ a TOP
embedding.
For each $\eps>0$ there exists $\dta>0$ such that any smooth embeddings
$g,h\:X\emb Q$, $\dta$-close to $f$, are smoothly $\eps$-isotopic.
Moreover if $g=h$ on the $\dta$-neighborhood of $Z$, then the isotopy can be
chosen fixing $g(Z)=h(Z)$.

(b+) In the (b) part, `$\eps$-isotopic' can be replaced with `$\eps$-ambient
isotopic'.
\endproclaim

The (b) and (b+) parts are proved later in this section.
We point out the following special case of 3.2.

\proclaim{Theorem 3.3}
(a) \cite{Ch1}, \cite{Mill1}, \cite{Co}
Let $X^n$ be a compact polyhedron, $Q^m$ a PL manifold, $m-n\ge 3$,
and $f\:X\emb Q$ a PL embedding.
For each $\eps>0$ there exists $\dta>0$ such that any PL embedding
$f'\:X\emb Q$, $\dta$-close to $f$, is PL $\eps$-ambient isotopic to $f$.

(b) \cite{Ki1, 2.1}
Let $X^n$ be a compact smooth manifold, $Q^m$ a smooth manifold,
$m>\frac{3(n+1)}2$, and $f\:X\emb Q$ a smooth embedding.
For each $\eps>0$ there exists $\dta>0$ such that any smooth embedding
$f'\:X\emb Q$, $\dta$-close to $f$, is smoothly $\eps$-ambient isotopic to $f$.
\endproclaim

\proclaim{Theorem 3.4} (a)
Let $X^n$ be a compact polyhedron, $Q^m$ a PL manifold, $m-n\ge 3$.
Any TOP isotopy $f_t$ between PL embeddings $f_0,f_1\:X\emb Q$
is $\eps$-approximable, for each $\eps>0$, by a PL isotopy $g_t$ between
$f_0$ and $f_1$.
Moreover if $f_t$ fixes a subpolyhedron $Z$ of $X$, then $g_t$ can be chosen
fixing $Z$.

(b) \cite{Hae1}, \cite{Ki1, 2.3}
Let $X^n$ be a compact smooth manifold, $Q^m$ a smooth manifold and
$m>\frac{3(n+1)}2$.
Any TOP isotopy $f_t$ between smooth embeddings $f_0,f_1\:X\emb Q$ is
$\eps$-approximable, for each $\eps>0$, by a smooth isotopy $g_t$ between
$f_0$ and $f_1$.
Moreover if $f_t$ fixes the $\dta$-neighborhood of a closed subset $Z$ of $X$,
then $g_t$ can be chosen fixing $Z$.
\endproclaim

\demo{Proof of 3.4a} (Compare to \cite{Mill3, proof of Theorem 3}, \cite{Lu};
see Remark 3.8 for an alternative proof.)
For each $t\in I$ let $U(t)$ denote an open neighborhood of $t$ in $I$ such
that for each $s\in U(t)$ the embedding $f_s$ is $\bta(t)$-close to $f_t$,
where $2\bta(t)=\dta\ind{3.2a}$ is given by 3.2a for
$\eps\ind{3.2a}=\frac{\eps}2$ and $f\ind{3.2a}=f_t$;
we can assume $\bta(t)<\frac{\eps}4$.

Since $I$ is compact, it can be covered by a finite number $k$ of open
intervals $U_1, U_2,\dots, U_k$, where $U_i=U(s_i)$ for some $s_i\in I$,
$s_1=0$, $s_k=1$ and $U_i\cap U_{i+1}\neq\emptyset$ for each $i=1,\dots,k-1$;
let $t_i$ be a point in $U_i\cap U_{i+1}$.
By 3.1a for each $i=1,\dots,k-1$ there is a PL embedding $g_i\:X\emb Q$,
agreeing with $f_{t_i}$ on $Z$ and such that
$\dist(g_i,f_{t_i})<\min(\bta(s_i),\bta(s_{i+1}))$.
We put $g_0=f_0$ and $g_k=f_1$.
Then $g_i$ is $2\bta(s_i)$-close to $f_{s_i}$ and $2\bta(s_{i+1})$-close to
$f_{s_{i+1}}$ for each $i=0,\dots,k$.
By 3.2a for each $i=0,\dots,k-1$ the embeddings $g_i$ and $g_{i+1}$ are PL
$\frac{\eps}2$-isotopic fixing $Z$.
The stacked composition of these isotopies is the required isotopy,
$\eps$-close to $f_t$. \qed
\enddemo

\proclaim{Theorem 3.5}
(a) Suppose that $X^n$ is a compact polyhedron, $Q^m$ a PL manifold,
$m-n\ge 3$, $Z$ a subpolyhedron of $X$, and $f\:X\emb Q$ a TOP embedding.
For each $\eps>0$ there exists $\dta>0$ such that for any PL embedding
$g\:X\emb Q$, $\dta$-close to $f$, there is an $\eps$-isotopy
$f_t\:X\to Q$ such that $f_0=g$, $f_1=f$ and such that $f_t$ for $t<1$ is a
PL isotopy.
Moreover, if $g|_Z=f|_Z$, then $f_t$ fixes $g(Z)=f(Z)$.

Furthermore, $f_t$ is covered by an $\eps$-homotopy $H_t\:Q\to Q$ such that:

$H|_{Q\x [0,1)}$ is a PL homeomorphism;

$H_0=\id_Q$ and $H_1\o g=f$;

if $g|_Z=f|_Z$, then $H_t$ fixes $g(Z)=f(Z)$.

(b) Suppose that $X^n$ is a compact smooth manifold, $Q^m$ a smooth manifold,
$m>\frac{3(n+1)}2$, $Z$ a closed subset of $X$, and $f\:X\emb Q$ a TOP
embedding.
For each $\eps>0$ there exists $\dta>0$ such that for any smooth embedding
$g\:X\emb Q$, $\dta$-close to $f$, there is an $\eps$-isotopy $f_t\:X\emb Q$
such that $f_0=g$, $f_1=f$ and such that $f_t$ for $t<1$ is a smooth isotopy.
Moreover if $f=g$ on the $\dta$-neighborhood of $Z$, then $f_t$ can be chosen
fixing $Z$.

(b+) In the (b) part it can be assumed that $f_t$ is covered by an
$\eps$-homotopy $H_t\:Q\to Q$ such that $H|_{Q\x [0,1)}$ is a diffeomorphism,
$H_0=\id_Q$, and $H_1\o g=f$.
\endproclaim

Theorem 3.5a is an immediate corollary of 3.1a and 3.2a.

\demo{Proof of 3.2b}
Triangulate $X$ and ambient isotop $g$ onto a PL embedding $H\o g$.
Let $Z'$ be a subpolyhedron of $X$ such that $g|_{Z'}=h|_{Z'}$ and
$Z\i Z'$.

By 3.1a and 3.5a $H\o h$ is TOP isotopic, by an arbitrarily small isotopy,
fixing $Z'$, to a PL embedding $h'$.
Hence $h$ and $H^{-1}\o h'$ can be assumed TOP $\frac{\eps}3$-isotopic
fixing $Z'$.

By 3.2a $H\o g$ and $h'$ can be assumed topologically (even PL)
isotopic, by a sufficiently small isotopy, fixing $Z'$.
Hence $g$ and $H^{-1}\o h'$ can be assumed TOP $\frac{\eps}3$-isotopic
fixing $Z'$.

Finally, by 3.4b the obtained TOP $\frac{2\eps}3$-isotopy between $g$ and $h$
can be approximated by a smooth $\eps$-isotopy fixing $Z$. \qed
\enddemo

Theorem 3.2b+ follows immediately from 3.2b and 1.13b.

Theorem 3.5b follows immediately from 3.1b and 3.2b.

Theorem 3.5b+ follows immediately from 3.1b and 3.2b+.

\proclaim{Theorem 3.6} Let $X^n$ be a compact polyhedron, $Q^m$ a PL
manifold, $m-n\ge 3$.

(a) Suppose that $f\:X\to Q$ is an embedding.
For each $\eps>0$ there exists $\dta>0$ such that if an embedding $g\:X\emb Q$
is $\dta$-close to $f$, then for any $\gma>0$ there is an $\eps$-ambient
isotopy, taking $g$ onto an embedding, $\gma$-close to $f$.

(a$\,'\!$) In addition, there is an $\eps$-ambient isotopy, taking $f$ onto an
embedding, $\gma$-close to $g$.

(b) Suppose that $f\:X\to Q$ is a map and $g\:X\to Q$ an embedding.
For each $\eps>0$ there exists $\dta>0$ such that if an embedding $g'\:X\emb Q$,
$\dta$-close to $g$, is taken onto $f$ by a pseudo-isotopy $H'_t$, then $g$ is
taken onto $f$ by a pseudo-isotopy $H_t$, $\eps$-close to $H'_t$.
\endproclaim

\demo{Proof of (a)}
Let $2\dta=\dta\ind{3.5a}$ be given by 3.5a for
$f\ind{3.5a}=f$ and $\eps\ind{3.5a}=\frac{\eps}2$.
In addition let $\dta'=\dta\ind{3.5a}$ be given by 3.5a for
$f\ind{3.5a}=g$ and $\eps\ind{3.5a}=\frac{\eps}2$.
We can assume that $\dta'<\dta$.
By 3.1a $g$ is $\dta'$-close to a PL embedding $h\:X\emb Q$.
By 3.5a $h$ can be taken by an $\frac{\eps}2$-pseudo-isotopy $G_t$ onto $g$
and by an $\frac{\eps}2$-pseudo-isotopy $F_t$ onto $f$.

Let $U\i Q\x I$ be the closed neighborhood of $G_1\o h(X)\x [0,1)$ in $Q\x I$
such that $U\cap Q\x 1=G_1\o h(X)\x 1$.
Then $G|_U$ is injective, and since $U$ is compact, the map
$G^{-1}|_U\:U\to G^{-1}(U)$ is uniformly continuous.
Hence for each $\bta>0$ there is a number $t_0<1$ such that the embedding
$h'=G^{-1}_{t_0}\o g$ is $\bta$-close to $h=G^{-1}_1\o g$.
The map $G'_t=G^{-1}_{t_0}\o G_{t_0(1-t)}$, $t\in I$, yields an
$\frac{\eps}2$-ambient isotopy taking $g$ onto $h'$.
Finally, since $F$ is uniformly continuous, the number $\bta$ can be chosen so
that for each $t\in I$ the embeddings $F_t\o h'$, $F_t\o h$ are
$\frac{\gma}2$-close, while $h'$ and $F_t\o h'$ are clearly
$\frac{\eps}2$-ambient isotopic.
If $t<1$ is such that $F_t\o h$ and $F_1\o h=f$ are $\frac{\gma}2$-close, then
$f$ and $F_t\o h'$ are $\gma$-close, while $F_t\o h'$ and $g$ are
$\eps$-ambient isotopic. \qed
\enddemo

\demo{Proof of (a$\,'\!$)} Proceed as in the proof of (a) until $F_t$, $G_t$
are constructed, and after that exchange their roles. \qed
\enddemo

\demo{Proof of (b)}
We can assume that $Q$ is compact, hence $H'_t$ is uniformly continuous.
For any fixed $t_0<1$ the map $(H'_t)^{-1}$, $t\in [0,t_0]$, is uniformly
continuous, and so is the map $H'_{st}=H'_t\o (H'_s)^{-1}$, $s\in [0,t_0]$,
$t\in I$.
For $k=0,1,\dots$ let $\lda_k>0$ be such number that
$\dist(H'_{st}(p),H'_{st}(q))<\lda\dist(p,q)$ whenever $p,q\in Q$,
$s\in [0,1-2^k]$, $t\in I$.
Let $\eps_k$, $k=0,1,\dots$ be a sequence of reals such that
$$\sum\limits_{k=0}^{\infty}\eps_k\lda_k<\eps.\tag{$*$}$$
Let $\dta=\dta\ind{3.6a}$ be given by the (a$'$) part for
$f\ind{3.6a}=g$ and $\eps\ind{3.6a}=\eps_0$.
Let $\dta_k=\dta\ind{3.6a}$ be given by the (a) part for
$f\ind{3.6a}=g'_k=H_{1-2^{-k}}\o g'$ and
$\eps\ind{3.6a}=\eps_k$, $k=1,2,\dots$.

For $k=0,1,\dots$ put $\gma_k=\frac{\dta_{k+1}}{\lda_{k+1}}$.
Then for any embedding $g_k$, $\gma_k$-close to $g'_k=H_{1-2^{-k}}\o g'$, and
for each $t\in [1-2^{-k},1-2^{-k-1}]$ the embedding $H'_{1-2^{-k}\!,\,t}\o g_k$
is $\dta_{k+1}$-close to $H'_t\o g'$.
In particular, the embedding $\bar g_{k+1}=H'_{1-2^{-k}\!,\,1-2^{-k-1}}\o g_k$
is $\dta_{k+1}$-close to $g'_{k+1}$, $k=0,1,\dots$.
The statement `$\bar g_k$ is $\dta_k$-close to $g'_k$' holds also for $k=0$ if
we put $\bar g_0=g$, $\dta_0=\dta$.

Now by (a) and (a$'$) for $k=0,1,\dots$ one can take an embedding $\bar g_k$,
which is $\dta_k$-close to $g'_k$, onto an embedding $g_k$, which is
$\gma_k$-close to $g'_k$, by an $\eps_k$-ambient isotopy $G^k_t\:Q\to Q$.
Then the stacked composition of isotopies
$$G^0_t;\ H'_t, t\in [0,\tfrac{1}2];\ G^1_t;\ H'_t, t\in
[\tfrac{1}2,\tfrac{3}4];\ G^2_t;\ \dots$$ yields a pseudo-isotopy $H_t$
(the possibility of continuous extension as $t\to 1$ is guaranteed by ($*$),
cf\. \cite{Ke3, Lemma 1}).
By ($*$) $H_t$ is $\eps$-close to $H'_t$.
Since $\dta_k\to 0$ as $k\to\infty$ and $H_t\o g=H'_{1-2^{-k},t}\o g_k$ is
$\dta_{k+1}$-close to $H'_t\o g'$, where $k=[-\log_2 (1-t)]$, we obtain that
$H_1\o g=H'_1\o g'=f$. \qed
\enddemo

\proclaim{Theorem 3.7} \cite{Ed1, 4.1} For each $\eps>0$ and a positive integer
$n$ there exists $\dta>0$ such that the following holds.

Let $X^n$ be a compact polyhedron and $Y^{n-1}$ its subpolyhedron,
$(Q^m,\partial Q)$ a PL manifold, $m-n\geq 3$, and $f\:(X,Y)\x [-1,1]\to
(Q,\partial Q)\x\R$ a PL embedding such that $\Pi\o f$ is $\dta$-close to $\pi$.
(Here $\Pi\:Q\x\R\to\R$, $\pi\:X\x [-1,1]\to [-1,1]\i\R$ denote the
projections).
Then there is a PL $\eps$-ambient isotopy $H_t$ with support in
$Q\x [-\eps,\eps]$, taking $f$ onto a PL embedding $g$ such that
$g^{-1}(Q\x J)=X\x (J\cap [-1,1])$ for each $J=(-\infty,0], 0, [0,+\infty)$.

Furthermore, for each $\gma>0$ given in advance it can be assumed that $P\o H_t$
moves points less than $\gma$, where $P\:Q\x\R\to Q$ denotes the projection.
Moreover, if $f^{-1}(\partial Q\x J)=Y\x J$ for each $J$ as above, $H_t$ can be
chosen to fix $\partial Q\x\R$.
\endproclaim

Theorem 3.7, called Slicing Lemma in \cite{Ed1}, \cite{Mill3}, was one of the
key steps in the proof of the (a) parts of Theorems 3.1 and 3.2.
(In the statement \cite{Ed1, 4.1(3)} one should read `$(h_1\o g)^{-1}$'
instead of `$g^{-1}$'.)
In view of an analogy between Lemmas 4.1 and 7.6 below, one can regard 3.7 as a
geometric version of the Freudental Suspension Theorem.
The proof of 3.7 in \cite{Ed1} is somewhat similar to the proof of the
Penrose-Whitedead-Zeeman-Irwin Embedding Theorem, meanwhile Miller proves a
generalization of 3.7 in \cite{Mill3} using his controlled version (see
\cite{Mill1}) of sunny collapsing (see \S5).

The following curious statement, not required in the rest of this paper, can
be regarded as an alternative controlled version of the Concordance Implies
Isotopy Theorem.
Call a concordance $F\:X\x I\emb Q\x I$ {\it $\eps$-level-disturbing} if for
each $t\in I$ there is a neighborhood $U(t)$ of $t$ in $I$ such that
$\Pi\o F(X\x t)\i U(t)$, $\Pi\:Q\x I\to I$ being the projection, and
$F|_{X\x U(t)}$ moves points less than $\eps$.
Notice that a 0-level-disturbing concordance is an isotopy, and the property
of being $\eps$-level-disturbing is independent on the choice of metric in $I$.

\proclaim{Theorem 3.7+} For each $\eps>0$ and a positive integer $n$ there
exists $\dta>0$ such that the following holds.
Consider a compact polyhedron $X^n$, a PL manifold $Q^m$, $m-n\ge 3$, and PL
embeddings $f,g\:X\emb Q$.
Then any PL $\dta$-level-disturbing concordance between $f,g$ is $\eps$-close
to a PL isotopy between $f,g$.
\endproclaim

\demo{Proof}
Let $F\:X\x I\emb Q\x I$ be the given concordance.
Without loss of generality $F(X\x I)\i Q\x I$.
Let $\gma=\dta\ind{1.13a}$ be given by 1.13a for
$\eps\ind{1.13a}=\frac{1}2\eps$ and suppose $\gma<\frac{1}2\eps$.
If $F$ moves points less than $\frac{1}2\gma$, then by 1.13a $F$ is
$(\frac{1}2\gma+\frac{1}2\eps)$-close to an isotopy between $f,g$.

Otherwise we can split $I$ into pieces $J_i$ such that $F|_{X\x J_i}$ moves
points less than $\frac{1}2\gma$ but more than $\frac{1}4\gma$ for each $i$.
Then for each $i$ there is a point $x_i\in X$ such that
$\diam P\o F(x_i\x J_i)>\frac{1}4\gma$, $P\:Q\x I\to Q$ being the projection.
Hence for each $i$ and any positive number $\dta<\frac{1}4\gma$ one can find
numbers $d_{i1}<\dots<d_{im}$, where $m=[\frac{\gma}{4\dta}]$, such that
$J_i=(d_{i1},d_{im})$ and
$\dist (P\o F(x_i\x d_{ij}), P\o F(x_i\x d_{i,j+1}))\ge\dta$.

Now suppose that $F$ is $\dta$-level-disturbing, then for each $t\in I$ there
exists $U(t)$ such that $\Pi\o F(X\x t)\i U(t)$, and $F|_{X\x U(t)}$ moves
points less than $\dta$.
Then for each $t\in I$, $U(t)$ contains at most one point $d_{ij}$.
Choose a metric on $I$ such that $\dist (d_{ij},d_{i,j+1})=\frac{\dta}2$, then
$\diam U(t)<\dta$, while $\diam J_i=\frac{1}2\dta m>\frac{1}{10}\gma$, provided
$\dta\le\frac{4}5\gma$.
Finally let $\dta=\dta\ind{3.7}$, which is obtained from Lemma 3.7 for
$\eps\ind{3.7}=\frac{1}{20}\gma$.

Then by 3.7 $F$ is $\frac{1}{20}\gma$-ambient isotopic (hence
$\frac{1}{20}\gma$-close) to a concordance $G$ between $f,g$ such that
$G(X\x J_i)\i Q\x J_i$ for each $i$.
Thus $G$ splits into concordances $G|_{X\x J_i}$, each moving points less than
$\frac{1}{20}\gma+\frac{1}2\gma+\frac{1}{20}\gma$.
By 1.13a each of them is $(\frac{12}{20}\gma+\frac{1}2\eps)$-close to an isotopy.
Together these isotopies yield an isotopy $\Phi\:X\x I\emb Q\x I$,
$(\frac{12}{20}\gma+\frac{1}2\eps)$-close to $G$, hence
$(\frac{13}{20}\gma+\frac{1}2\eps)$-close (thus $\eps$-close) to $F$. \qed
\enddemo

\remark{Remark 3.8}
Given a homotopy $H\:X\x I\to Q\x I$, for each $\eps>0$ it is easy to
find $\dta>0$ (depending on $H$) such that any concordance,
$\dta$-close to $H$, is $\eps$-level-disturbing.
Taking into account 3.1a and 3.7+, we thus obtain an alternative proof of 3.4a.
\endremark

\head 4. Proof of Theorem 1.6 \endhead

The following lemma is a corollary of Theorem 3.7:

\proclaim{Lemma 4.1} Let $X^n$ be a finite simplicial complex, $Q^m$ a
combinatorial manifold, $m-n\geq 3$, $f\:X\to Q$ a simplicial map and
$C$ a union of some top-dimensional dual cells of $Q$.
Then for each $\eps>0$ there exists $\dta>0$ such that the following holds.
Suppose that $g\:X\emb Q$ is a PL embedding, $\dta$-close to $f$ and such
that $g^{-1}(B)=f^{-1}(B)$ for each dual cell $B$ of $C$.
Then $g$ is PL $\eps$-ambient isotopic, keeping $C$ fixed, to a PL embedding
$h\:X\emb Q$ such that $h^{-1}(B)=f^{-1}(B)$ for each dual cell $B$ of $Q$.
\endproclaim

Suppose that $\dim f(X)=k$, and $\Cl{Q\but C}$ consists of $l$ top-dimensional
dual cells.
Denote by 4.1$(k,l)$ the statement of 4.1 for $k$ and $l$.
Then 4.1$(i,0)$ and 4.1$(0,j)$ are trivial for any $i$, $j$.
Assuming that 4.1$(i,j)$ is proved for $i<k$ and arbitrary $j$, and for $i=k$
and $j<l$, let us prove 4.1$(k,l)$.

\demo{Proof of 4.1$(k,l)$}
Choose any vertex $v$ of $Q$ outside $C$.
Let $D=\st(v,Q')$ be its dual cell, and write
$E=\Cl{\partial D\but\partial C}$.
Notice that the pair $(E,\partial E)$ is bi-collared in
$(\Cl{Q\but C},\partial C)$.
By Theorem 3.7, for any $\gma>0$ the number $\dta$ can be chosen so
that $g$ is PL $\gma$-ambient isotopic, keeping $C$ fixed, to a PL embedding
$\phi\:X\emb Q$ such that $\phi^{-1}(E)=f^{-1}(E)$.
It follows that, in addition, $\phi^{-1}(D)=f^{-1}(D)$.

Pseudo-radial projection \cite{RoS} yields a PL homeomorphism
$\partial D\to\partial\st(v,Q)$ (in general the latter complex does not
coincide with $\lk(v,Q)$), which takes the intersection of $D$ with a simplex
of $Q$ onto a simplex of $\partial\st(v,Q)$ and a dual cell of $Q$, lying in
$\partial D$, onto a dual cell of $\partial\st(v,Q)$.
We apply 4.1$(k-1,l')$ in $\partial D$, equipped with triangulation, inherited
from $\partial\st(v,Q)$, where $l'$ is the number of dual cells of $Q$ in
$C\cap\partial D$.
Using collaring, we obtain that for each $\bta>0$ we can choose $\gma+\dta$
so that $\phi$ (which is $(\gma+\dta)$-close to $f$), is PL $\bta$-ambient
isotopic, keeping $C$ fixed, to a PL embedding $\psi\:X\emb Q$ such that
$\psi^{-1}(B)=f^{-1}(B)$ for each dual cell $B$ of $C\cup D$.

By 4.1$(k,l-1)$, for any $\alp>0$ the number $\bta+\gma+\dta$ can be chosen
so that $\psi$ (which is $(\bta+\gma+\dta)$-close to $f$) is PL $\alp$-ambient
isotopic, keeping $C\cup D$ fixed, to a PL embedding $h\:X\emb Q$
such that $h^{-1}(B)=f^{-1}(B)$ for each dual cell $B$ of $Q$.
Thus $g$ is $\eps$-ambient isotopic to $h$, keeping $C$ fixed, provided
$\alp+\bta+\gma<\eps$. \qed
\enddemo

\proclaim{Lemma 4.2}
Let $X^n$ be a finite simplicial complex, $Q^m$ a combinatorial manifold,
$m-n\geq 3$, and $f\:X\to Q$ a simplicial map.
If $h\:X\emb Q$ is a PL embedding such that $h^{-1}(B)=f^{-1}(B)$ for each
dual cell $B$ of $Q$, then $h$ is taken onto $f$ by a PL pseudo-isotopy
$H_t\:Q\to Q$ such that $H_t(B)=B$ for each dual cell $B$ of $Q$.
\endproclaim

\demo{Proof}
Put $H_0=\id_Q$ and $H_t=\id$ outside $N=\N(f(X),Q)$.
Let $A_1,\dots,A_m$ be the dual cells of $N$, except for those in $\partial N$,
arranged in an order of increasing dimension.
Assuming that $H$ is defined on $A_j\x I$ for all $j<i$ (hence on
$(\partial A_i)\x I$), extend it to $A_i\x I$ as follows.

Denote the cone point of $A_i$ by $a_i$.
Let $R$ be a relative regular neighborhood in $A_i$ of
$\partial A_i\cup h(f^{-1}(A_i))$ modulo $h(f^{-1}(a_i))$, and put
$P=\Cl{A_i\but R}$.
Then $h^{-1}(P)=f^{-1}(a_i)$, and we define $H_t(p)=t*a_i+(1-t)*p$ for
each $p\in P$ (we use here the cone structure $a_i*\partial A_i$ on $A_i$).

The quotient space $A_i/P$ is PL homeomorphic keeping $\partial A_i$ fixed to
$A_i=a_i*\partial A_i$, and $f^{-1}(A_i)/f^{-1}(a_i)$ is PL homeomorphic
keeping $f^{-1}(\partial A_i)$ fixed to the cone $a_i*f^{-1}(\partial A_i)$.
Denote these homeomorphisms by $\phi$ and $\psi$, respectively, and let
$$i\:A_i\but P\emb A_i/P,\qquad j\:f^{-1}(A_i)\but f^{-1}(a_i)\emb
f^{-1}(A_i)/f^{-1}(a_i)$$ be the natural inclusions.
Let $h'\:a_i*f^{-1}(\partial A_i)\emb a_i*\partial A_i$ be the embedding
defined by the identity on $a_i$ and by $\phi\o i\o h\o j^{-1}\o\psi^{-1}$
elsewhere.
By the Lickorish Cone Unknotting Theorem \cite{Li} there is a PL homeomorphism
$\lambda\:a_i*\partial A_i\to a_i*\partial A_i$ keeping $\partial A_i$ fixed
and such that $\lambda\o h'$ is the conical map
$\id_{a_i}*h'|_{f^{-1}(\partial A_i)}$.
Define an isotopy $\Lambda\:A_i\x I\to A_i\x I$ by $$\Lambda=(\id_{a_i\x 1}*
(\lambda^{-1}\x\id_0\cup\id_{\partial A_i\x I}))\o(\lambda\x\id_I),$$ then
$\Lambda\o (h'\x\id_I)$ is the conical map $\id_{a_i\x 1}*
(h'\x\id_0\cup h'|_{f^{-1}(\partial A_i)}\x\id_I)$.

Extend $H'=H|_{\partial A_i\x I\cup A_i\x 0}$ conewise to obtain the map
$$\id_{a_i\x 1}*H'\:(a_i\x 1)*(\partial A_i\x I\cup A_i\x 0)\to
(a_i\x 1)*(\partial A_i\x I\cup A_i\x 0).$$
Finally, define $H$ on $(A_i\but P)\x I$ by
$$H=(\id_{a_i\x 1}*H')\o\Lambda\o(\phi\x\id_I)\o (i\x \id_I).$$
Clearly, $H$ is well-defined, is PL, and $H_1$ takes $h$ onto $f$.

Assuming that $H|_{A_j\x I}$ is level-preserving for each $j<i$, we see
from the construction above that so is $H|_{A_i\x I}$.
Assuming that $H$ is a homeomorphism on $A_j\x [0,1)$ for each $j<i$ and
recalling that it is a homeomorphism on $Q\x 0$, we see that $H$ is
a homeomorphism on $A_i\x [0,1)$. \qed
\enddemo

\demo{Proof of 1.6a}
Triangulate $X$ and $Q$ so that $f$ is simplicial and each dual cell of $Q$ is
of diameter $<\frac{\eps}2$.
Let $\dta$ be less than $\dta\ind{4.1}$, which is given by Lemma 4.1
for $\eps\ind{4.1}=\frac{\eps}2$, and apply Lemma 4.1 to
obtain a PL $\frac{\eps}2$-isotopy, taking $g$ onto a PL embedding
$h\:X\emb Q$ such that $h^{-1}(B)=f^{-1}(B)$ for each dual cell $B$ of $Q$.
Finally, apply Lemma 4.2 to obtain a PL pseudo-isotopy $H_t\:Q\to Q$ taking
$h$ onto $f$.
Then $H$ moves no point as much as the maximal diameter of a dual cell of $Q$,
which, in turn, is less than $\frac{\eps}2$. \qed
\enddemo

\demo{Proof of 1.6b}
By 3.6b and 3.1a we can assume that the $\dta$-close to $i\o f$ embedding
$g\:X\emb Q$ is PL.
Without loss of generality $f$ is surjective, hence we can assume
$\dim Y\le n\le m-3$.
Then by 3.1a and 3.5a there is a PL embedding $j\:Y\emb Q$ and a
pseudo-isotopy $H_t$, taking $j$ onto $i$.
For any $\gma>0$ we can assume that $H_t$ moves points less than $\gma$ and
$j\o f$ is $\gma$-close to $i\o f$.

The PL embedding $g$ is $(\gma+\dta)$-close to the PL map $j\o f$, and one
could attempt to apply Theorem 1.6a here.
But this is impossible, for one cannot make $\gma$ as small as required keeping
$j\o f$ unchanged.
The solution is to use uniform continuity (as in the proof of 3.6a).

Let $U\i Q\x I$ be a closed neighborhood of $H_1\o j(Y)\x [0,1)$ in $Q\x I$
such that $U\cap Q\x 1=H_1\o j(Y)\x 1$.
Then $H|_U$ is injective, and since $U$ is compact, the map
$H^{-1}|_U\:U\to H^{-1}(U)$ is uniformly continuous.
Clearly, for each $t_0<1$ the number $\dta>0$ can be chosen so that the image
of the embedding $g\x t_0\:X\to Q\x t_0$ lies in $U$.
Now $g\x t_0$ is $(1-t_0+\dta)$-close to $(i\o f)\x 1$, therefore for each
$\bta>0$ the numbers $t_0$, $\dta$ can be chosen so that
$g'=H^{-1}_{t_0}\o g$ is $\bta$-close to $j\o f=H^{-1}_1\o (i\o f)$.
The map $G_t=H^{-1}_{t_0}\o H_{t_0(1-t)}$, $t\in I$, yields a $\gma$-ambient
isotopy taking $g$ onto $g'$.

By 1.6a for any $\alp>0$ the number $\bta$ can be chosen so that $g'$ is taken
onto $j\o f$ by an $\alp$-pseudo-isotopy $F_t\:Q\to Q$.
Then the `diagonal' $(\alp+\gma)$-pseudo-isotopy $\Phi_t=H_t\o F_t\:Q\to Q$
takes $g'$ onto $i\o f$.
Since $g$ is $\gma$-ambient isotopic to $g'$, there is an $\eps$-pseudo-isotopy
taking $g$ onto $i\o f$, provided $2\gma+\alp<\eps$. \qed
\enddemo

\head 5. Proof of Theorem 1.13. \endhead

\definition{Definition}
A subcomplex $Y$ of a simplicial complex $X$ is said to be {\it locally of
codimension $\geq k$} in $X$, if every $n$-simplex of $Y$ faces some
$(n+k)$-simplex of $X$ \cite{Li}.
We call $S(f)=\Cl{\{x\in X\mid f^{-1}f(x)\neq x\}}$ the {\it singular set} of
a map $f\:X\to Y$.
\enddefinition

\proclaim{Lemma 5.1} \cite{HL}, \cite{Bi3}
Let $X^n$ be a compact polyhedron and $Q^m$ a PL manifold, $m-n\ge 3$,
and $p\:X\x I\to X$, $P\:Q\x I\to Q$ the projections.
For any PL embedding $F\:X\x I\emb Q\x I$ and any $\eps>0$ there is a
PL level-preserving $\eps$-homeomorphism $H\:Q\x I\to Q\x I$ such that
$S(P\o H\o F)$ is locally of codimension $\geq 2$ in $X\x I$.

Moreover, one can choose $H$ so that $p|_{S(P\o H\o F)}$ is non-degenerate.
Furthermore, the preimage of any point under $P\o H\o F$ contains at most
$\phi(n)=[\frac{n+1}{3}]+1$ points.
\endproclaim

\definition{Definition}
Let us think of the second factor of $Q\x I$ as of height (that is, a point
$(q_1,t_1)$ lies below a point $(q_2,t_2)$ if $q_1=q_2$ and $t_1<t_2$).
If $X\i Q\x I$, let $\sh X$ denote a {\it shadow} of $X$, the set
of points of $Q\x I$ lying below some point of $X$.

We say that a collapse $X\col Y$ in $Q\x I$ is a {\it simple sunny collapse},
if no point of $X\but Y$ lies in $\sh X$.
A sequence of simple sunny collapses is called a {\it sunny collapse}
\cite{Ze}.
Let us say that a sunny collapse is {\it $m$-complex}, if it consists of at
most $m$ simple sunny collapses.
Repeating the same for $\Cl{X\but Y}$ instead of $X\but Y$, we define a
(simple/$m$-complex) {\it stable} sunny collapse \cite{Me1}.
\enddefinition

\example{Example 5.2} (compare to \cite{Ze, Remark on p\. 510})
Let us illustrate the relation between sunny collapsing and unknotting.
Evidently, $I$ collapses onto $0$.
Let $F\:I\to I^3$ be a PL embedding such that $F(i)\i\Int(I^2\x i)$, $i=0,1$.
It turns out that if a collapse $F(I)\col F(0)$ is sunny, $F$ is unknotted.
Indeed, define a PL isotopy $H_t\:I\emb I^3$ by $s\mapsto F(s)$ for $s\leq 1-t$
and by mapping $(1-t,1]$ linearly onto points lying above $F(1-t)$.
Then $H_0=F$ and $H_1$ is linear.
Clearly, $H$ is locally flat, hence by \cite{RoS} it extends to a PL ambient
isotopy, which `unknots' $F$.

Surprisingly, if $F\:I\emb I^3$ maps $0$ into $\Int(I^2\x 1)$ and $1$ into
$\Int(I^2\x 0)$, then $F$ can be knotted even if there is a sunny collapse
$F(I)\col F(0)$.
However, in all other cases of PL embeddings $F\:(I,\partial I)\emb
(I^3,\partial I^3)$ existence of a sunny collapse $F(I)\col F(0)$ implies that
$F$ is unknotted.
Indeed, in the case $F(1)\i\Cl{\partial I^3\but (I^2\x0)}$ we use that
$\N(0,I)$ is not overshadowed by $I$ to shift $F(0)$ upwards into
$(\partial I^2)\x 1$.
Then we apply the above construction of $H_t$ for $t\leq 1-\eps$, where
$\eps>0$ is the minimal distance between vertices in $F(I)$.
Now $H_{1-\eps}$ consists of two linear pieces, hence is unknotted.
To manage with the case $F(0)\i\Cl{\partial I^3\but (I^2\x1)}$, notice that a
collapse $F(I)\col F(0)$ is sunny, iff sunny is the analogous collapse
$U\o F\o u(I)\col U\o F\o u(0)$, where $u\:I\to I$ and $U\:I^3\to I^3$ are
defined by $t\mapsto 1-t$ and $(r,s,t)\mapsto (r,s,1-t)$, respectively.
\endexample

\proclaim{Lemma 5.3} \cite{HL, Lemma 2}
Let $X$ be a simplicial complex, $Q$ a combinatorial manifold, $p\:X\x I\to X$
and $P\:Q\x I\to Q$ simplicial projections.
Let $G\:X\x I\emb Q\x I$ be a simplicial embedding satisfying the conclusion
of 5.1 and such that $G(X\x 0)\i Q\x 0$.
Then there is a sunny collapse $G(X\x I)\col G(X\x 0)$ such that

(i) $\tr Z\x I\i\N(Z\x I, X\x I)$ for any simplex $Z$ of $X$.
\endproclaim

Speaking informally, the main idea of the proof of Lemma 5.3 was to use
codimension 2 (that is, connectedness of $G(X\x I\but S)$) to have a
simultaneous collapsing access to all the $m$-simplices of $G(S)$
succesively for $m=n-1,\dots,0$, which enabled to collapse them in the order
they overshadow each other.
See \cite{Ze, proof of Lemma ~9} for a detailed proof of a similar statement.

\proclaim{Addendum to 5.3} The sunny collapse $G(X\x I)\col G(X\x 0)$
can be chosen $\psi(n)=\frac{\ (n+7)^2}{6}$-complex.
\endproclaim

\demo{Proof} Arrange the simplices of $K=G(S)$ in the following order.
Assuming that the order is defined in the case $\dim K\leq m$, define it
when $\dim K=m$, as follows.
Let first go all the top-dimensional non-overshadowed simplices, then all the
top-dimensional once overshadowed and so on, up to the top-dimensional
simplices, overshadowed by $\phi(m)-1$ ones.
After that put all at most $(m-1)$-dimensional simplices of $K$, arranged in
the order given by the inductive assumption.
The proof of Lemma 5.3 actually allows to collapse the simplices of $K$ in
any order of decreasing dimension, given in advance, particularly in the above.
Clearly, the obtained collapse is sunny, and since
$\psi(n)\geq\phi(n-1)+\dots+\phi(0)$, it is $\psi(n)$-complex. \qed
\enddemo

It turns out that any sunny collapse can be improved to a stable sunny one.
We prove this by following the given collapse with a slight but precisely
calculated lag.

\proclaim{Lemma 5.4}
Let $Q$ be a combinatorial manifold, $P\:Q\x I\to Q$ a simplicial projection
and $K_0\supset\dots\supset K_N$ be a sequence of subcomplexes of $Q$ such that
$\sh K_i\cap K_0\i K_{i+1}$ for each $i\leq N$, where $K_{N+1}=\emptyset$, and
suppose that $P|_{K_0}$ is non-degenerate.
Then there is a sequence $K_0=U_0\supset\dots\supset U_M=K_N$ of subpolyhedra
of $Q$ such that $\sh U_j\cap U_0\i\Int U_{j+1}$ for each $j\leq M$, where
$U_{M+1}=\emptyset$ and

(i) $K_0\!\col\!\dots\!\col\! K_i\!\col\!\dots\!\col\! K_N$ (simplicially)
implies $U_0\!\col\!\dots\!\col\! U_j\!\col\!\dots\!\col\! U_M$;

(ii) the trace of any simplex $Z$ of $K_0$ under $U_0\!\col\! U_M$
lies in that under $K_0\!\col\! K_N$.
\endproclaim

Actually, in the application of 5.4 the hypothesis of (i) will be fulfilled; we
allow it not to be fulfilled only to carry out induction in the proof of 5.4.
The prototypes of Lemma 5.4 can be found in \cite{Hu, proof of Proposition 5.1}
and \cite{Me1, Lemma 4.1}.
To prove 5.4 we need a couple of preliminary observations.

\proclaim{Claim 5.5} There is a second derived subdivision $\alp K_0$ of $K_0$
such that for any subcomplex $Y$ of $K_0$, the inclusion $\sh Y\cap K_0\i Y$
implies $$\sh\N(Y,\alp K_0)\cap K_0\i\Int\N(Y,\alp K_0).$$
\endproclaim

\demo{Proof} Let $K_0'$ be the barycentrically derived subdivision of $K_0$
and construct a derived subdivision $\alp K_0$ of $K_0'$ as follows.
For each simplex $A$ of $K_0$ define a map $f_A\:A\to\R^1$ by
$\partial A_j\mapsto -1$, $\hat A\mapsto\phi\o\Pi(a)$ and extending linearly,
where $\hat A$ denotes the barycenter of $A$, $\Pi\:Q\x I\to I$ denotes the
projection, and $\phi$ maps $[0,1]$ linearly onto $[\frac{1}{100},1]$.
Let $\Cal F=(A_0\supsetneq\dots\supsetneq A_m)$ run over the flags of simplices
in $K_0$ and let $B_i=B_i(\Cal F)=\hat A_i*\dots*\hat A_m$.
Then $B_0$ runs over the simplices of $K_0'$.
Define a derivation point of $B_0$ by
$d(B_0)=\hat A_0*\hat B_1\cap f^{-1}_{A_0}(0)$, unless $m=0$.
The subdivision $\alp K_0$ is defined.

It is easy to see that $d(B_m)=B_m$ lies in a subcomplex $Y$ of $K_0$ if
and only if $d(B_0)*\dots*d(B_m)$ (or, equivalently, $d(B_0)$) lies in
$\N(Y,\alp K_0)$. Notice that $\hat B_0\in\hat A_0*\hat B_1$ and
$$f_{A_0}(\hat B_0)=\frac{f_{A_0}(\hat A_0)-m}{m+1}<0$$ (unless $m=0$), hence
$\hat B_0\in\Int(d(B_0)*\hat B_1)$.
An induction on $m$ implies that $\hat B_0\in\Int(d(B_0)*\dots*d(B_m))$,
and moreover, that if $x\in\Int(d(B_0)*\hat B_1)$ then
$x\in\Int(d(B_0)*\dots*d(B_m))$.

Now suppose that $d=d(B_0)$ overshadows a point $d^*$ of $K_0$.
Then there are a flag $\Cal F^*=(A_0^*\supsetneq\dots\supsetneq A_m^*)$ and
the simplices $B_i^*=B_i(\Cal F^*)$, overshadowed respectively by a flag
$\Cal F=(A_0\supsetneq\dots\supsetneq A_m)$ and the simplices $B_i=B_i(\Cal F)$
(in the sense that each each point, say, of $A_0^*$ is overshadowed by, or
coincides with a point of $A_0$) and such that $d^*\in\hat A_0^**\hat B_1^*$.
If $d\in\N(Y,\alp K_0)$ then $B_m\i Y$, and since $B_m$ overshadows, or
equals, to $B_m^*$, we obtain $B_m^*\i Y$.
Consequently $d(B_0^*)*\dots*d(B_0^*)\i\N(Y,\alp K_0)$.
Finally, $\hat A_0^*\neq\hat A_0$, hence
$f_{A_0^*}(\hat A_0^*)<f_{A_0}(\hat A_0)$, therefore
$f_{A_0^*}(d^*)<f_{A_0}(d)$ and $f_{A_0^*}(d^*)<0$.
Thus by the above $d^*\in\Int\N(Y,\alp K_0)$. \qed
\enddemo

\proclaim{Claim 5.6} Let $K$ be a simplicial complex and $A$ its simplex.

(a) If $V\supset W$ in $\lk(A,K'')$ then
$\N(\partial A,K'')\cup A*V\col\N(\partial A,K'')\cup A*W$.
Moreover, $\tr Z\i Z$ for any simplex $Z$ of $K$.

(b) If $K\col L$ simplicially, then $K\col\N(L,K'')$.
Moreover, the trace of any simplex $Z$ of $K$ under the second collapse lies in
that under the first.
\endproclaim

\demo{Proof of (a)} Suppose that a simplex $A$ (strictly) faces a simplex $B$.
Since a ball collapses onto its face, $\N(A,B'')\col\N(A,\partial B'')\cup
\N(\partial A,B'')$.
Applying this to $B$ runnung over the simplices of $K$ which are faced by $A$
and meet $V\but W$, in order of decreasing dimension, we obtain the required
collapse. \qed
\enddemo

\demo{Proof of (b)} For each elementary collapse $K_i\col K_{i+1}$ it suffices
to prove that $N(K_i,K'')\col N(K_{i+1},K'')$.
Suppose that $K_i\col K_{i+1}$ goes from $A_i$ along $B_i$.
Apply the full collapse of (a) first to $A=A_i$ and then to $A=B_i$ to obtain
a collapse $\N(K_i,K'')\col\N(K_{i+1},K'')\cup K_i$.
Finally, since a ball collapses onto its face, $V\cup\N(W,K'')\col\N(W,K'')$.
By the moreover part of (a) and since $\tr A_i\i B_i$ under the last collapse,
the trace of any simplex $Z$ of $K_0$ under the obtained collapse
$K\col\N(L,K'')$ lies in that under $K\col L$. \qed
\enddemo

\demo{Proof of 5.4} Assume that 5.4 is proved for $\dim K_0<n$ and prove
it for $\dim K_0=n$.
We will construct a descending sequence of subpolyhedra $U_*$ (with several
indices) in $K_0$, arranged lexicographically, so that the lexicographic
unwrapping of indices yields the required sequence
$U_0\supset\dots\supset U_M$.

Let $\alp K_0$ be the subdivision given by 5.5. Define
$U_i=K_i\cup\N(K_{i+1},\alp K_0)$, $i\leq N$ and insert between them
$U_{i,0}=\N(K_{i+1},\alp K_0)$, $i<N$.
Then $U_0=K_0$, $U_N=K_N$ and $U_i\but U_{i,0}\i K_i\but K_{i+1}$.
By 5.5 and since $\sh K_i\cap K_0\i K_{i+1}\i\Int U_{i,0}$, we obtain
$\sh U_i\cap U_0\i\Int U_{i,0}$.
By 5.6b, $K_i\col K_{i+1}$ implies $K_i\col\N(K_{i+1},\alp K_i)$, or,
equivalently, $U_i\col U_{i,0}$, and the trace of any simplex $Z$ of $K$
under the last collapse lies in that under the first.

It remains to insert subpolyhedra in between $U_{i-1,0}$ and $U_i$. Let
$A_1,\dots,A_T$ be the simplices of $K_i\but K_{i+1}$, arranged in an order of
decreasing dimension and put $B_j=\cup\{A_k\mid k>j\}$.
Define $U_{ij}=K_i\cup\N(K_{i+1}\cup B_j,\alp K_0)$, $j\leq T$, then $U_{i,0}$
is same as above and $U_{i,T}=U_{i+1}$.
By 5.5 $\sh U_{ij}\cap U_0\i\Int U_{ij}$ for each $j\leq T$. Unfortunately
$$\sh U_{ij}\cap U_0\not\i\Int U_{i,j+1}$$ in general, so we should insert yet
more subpolyhedra.
Put $L_l=\lk(A_j,K_l)$ for each $l\leq i$.
Since $\sh K_l\cap K_0\i K_{l+1}$, we have that $\sh L_l\cap L_0\i L_{l+1}$.
Now $\dim L_0<n$ and we can apply the inductive hypothesis to obtain a sequence
of subpolyhedra $L_0=V_0\supset\dots\supset V_R=L_i$ such that
$\sh V_k\cap V_0\i\Int_{V_0}V_{k+1}$, $k\leq R$, where $V_{R+1}=\emptyset$.
Here `$\Int_{V_0}$' denotes topological interior in $V_0$.
Since $P|_{K_0}$ is simplicial and non-degenerate,
$$\sh(A_j*V_k)\cap\Int(A_j*V_0)\i\Int(A_j*V_{k+1}).$$
We put $W_k=A_j*V_k\cap\N(A_j,\alp K_0)$ for $k\leq R$, then
$W_0\cup\N(\partial A_j,\alp K_0)=\N(A_j,\alp K_0)$ and $W_R\i K_i$.
Furthermore, $\sh W_k\cap W_0$ lies in
$$\multline \sh(A_j*V_k)\cap\sh\N(A_j,\alp K_0)\cap\N(A_j,\alp K_0)\i
\sh(A_j*V_k)\cap\Int\N(A_j,\alp K_0)\i\\
\i\sh(A_j*V_k)\cap(\Int(A_j*V_0)\cap\Int\N(A_j,\alp K_0)\cup
                         \Int\N(\partial A_j,\alp K_0))\i\\
\i\Int(A_j*V_{k+1})\cap\Int\N(A_j,\alp K_0)\cup\Int\N(\partial A_j,\alp K_0)=
\Int(W_{k+1}\cup\N(\partial A_j,\alp K_0)). \endmultline$$
Finally, define $U_{ijk}=U_{i,j+1}\cup W_k$, $k\leq R$.
Then by the above $U_{ij,0}=U_{ij}$ and $U_{ij,R}=U_{i,j+1}$, while
$\sh U_{ijk}\cap U_0\i\Int U_{ij,k+1}$.
By 5.6a $U_{ijk}$ collapses onto $U_{ij,k+1}$ for all $k<R$ and $\tr Z\i Z$
under this collapse for any simplex $Z$ of $K_0$. \qed
\enddemo

\proclaim{Addendum to 5.4} $M$ can be chosen equal to $\xi(N,n)=N^{n+1}n!$.
\endproclaim

\demo{Proof} Prove this by induction on $n$.
Clearly, we can choose $M=N$ if $n=0$.
Since $M$ originally depends on an arbitrarily great number $T$, we should
redefine the subpolyhera $U_*$ so that it does not.
Notice that since $\sh A_j\cap K_0\i K_{i+1}$ for each simplex $A_j$ of
$K_i\but K_{i+1}$, by 5.5
$\sh\N(A_j,\alp K_0)\cap K_0\i\Int\N(K_{i+1}\cup A_j,\alp K_0)$.
Hence
$$\sh W_k\cap K_0\i\Int(W_{k+1}\cup\N(K_{i+1}\cup\partial A_j,\alp K_0)),$$
where $W_k=W_k(A_j)$ is defined for each simplex $A_j$ of $K_i\but K_{i+1}$ as
in the proof of 5.4, $k\leq R$ (by the inductive hypothesis we can choose
$R=\xi(N,n-1)$ to be the same for all $A_j$).
We redefine the subpolyhedra $U_*$ by
$$\alignedat 2
U_{ij}&=K_i\cup\N(K_{i+1}\cup(K_i\but K_{i+1})^{(n-j)},\alp K_0),
\qquad &j\leq n,\,\\
U_{ijk}&=U_{i,j+1}\cup\{W_k(A_l)\mid \dim A_l=n-j\},
\qquad &k\leq R.\endalignedat$$
Then by the above $\sh U_{ijk}\cap U_0\i U_{ij,k+1}$ and the statement is
fulfilled for the new sequence of subpolyhedra $U_*$, while there are only
$M=NnR$ of them.
Hence we can choose $\xi(N,n)=Nn\xi(N,n-1)=N^{n+1}n!$. \qed
\enddemo

\demo{Proof of 1.13a}
Let $F\:X\x I\to Q\x I$ be the given concordance between $F_0=f$ and
$F_1=g$.
Write $G=H\o F$ and $S=S(P\o G)$.
Let $\zeta(n)=\xi(\psi(n),n)$ and $\dta=\frac{\eps}{5\zeta(n)+1}$.
Subdivide $Q\x I$ and $X\x I$ so that $G$, $p$ and $P$ are simplicial and
$\gma=\mesh(Q\x I)$ is less than $\min(\dta,\frac{1}3\dta\ind{3.3a})$,
where $\dta\ind{3.3a}$ is obtained from Theorem 3.3a for
$f\ind{3.3a}=f$ and $\eps\ind{3.3a}=\dta$.

Apply Lemma 5.4 to the sunny collapse of Lemma 5.3 to obtain a stable sunny
collapse $G(X\x I)=U_0\col\dots\col U_M=G(X\x 0)$.
By Addenda, it consists of at most $\zeta(n)$ simple stable sunny collapses.
Also, $\tr G(Z\x 1)\i N_\gma(G(Z\x I))$ for any simplex $Z$ of $X$.
Since $F$ is a $\dta$-concordance and $H$ is arbitrarily, say,
$\dta$-close to the identity, $G$ is a $2\dta$-concordance.
Hence $\tr G(Z\x 1)\i\N_{2\dta+\gma}((G_0\x\id_I)(Z\x I))$.

By the simple stable condition, the projection $P\:Q\x I\to Q$, restricted to
$\Cl{U_i\but U_{i+1}}\cup\im_i G(X\x 1)$, where $\im_i$ denotes the image under
the first $i$ collapses, is a homeomorphism for each $i<M$.
Notice that
$$\Cl{U_i\but U_{i+1}}\cup\im_i G(X\x 1)\i U_i\but(\Int U_{i+1}\but X\x 1)$$
(the inclusion follows by an induction on $i$).
Since $\Cl{U_i\but U_{i+1}}\cup\im_i G(X\x 1)$ collapses onto
$\im_{i+1} G(X\x 1)$, there is a sequence of at most $\zeta(n)$ collapses
$$P(\Cl{U_i\but U_{i+1}}\cup\im_i G(X\x 1))\col P(\im_{i+1} G(X\x 1)),
\tag{$*$}$$
each of diameter at most $2\dta+\gma<3\dta$.
Subdivide $Q$ so that these collapses are simplicial. Our next goal is to
obtain a sequence of isotopies, using the following
\enddemo

\proclaim{Lemma 5.7} \cite{Mill1, Proposition 1}, \cite{Mill3, proof of
Theorem 11}
Let $Q$ be a combinatorial manifold, $V$ and $W$ be its subcomplexes such that
$V$ collapses simplicially to $W$.
Then there is a PL ambient isotopy $H_t$ of $Q$, $H_0=\id_Q$ such that for any
subcomplex $Z$ of $V$ and arbitrary $t\in I$

(i) $H_1\N(V,Q'')=\N(W,Q'')$;

(ii) $H_t$ is the identity outside
$\N(\N(\text{\rm vertices in }(V\but W)',Q''),Q''')$;

(iii) $H_1\N(Z,Q'')\i\N(\im Z,Q'')$;

(iv) $H_t\N(\N(Z,Q''),Q''')\i\N(\N(\tr Z,Q''),Q''')$.
\endproclaim

\proclaim{Addendum to 5.7} \cite{Mill1, Corollary 2} If the diameter of
the collapse is less than $\alp$ and $\mesh Q<\gma$, then $H$ moves points
less than $\alp+2\gma$. \qed
\endproclaim

\demo{Proof of 1.13a (continued)} By 5.3, 5.4 and 5.7(iii) there is
a sequence of ambient isotopies $h_t^i$ of $Q$ such that for any simplex $Z$ of
$X$ $$h_1^i\N(P(\im_i G(Z\x 1)),Q'')\i\N(P(\im_{i+1}G(Z\x 1)),Q'').$$
Let $h_t$ be the stacked composition of $h_t^0,\dots,h_t^M$.
Then by the addendum to 5.7 $h_t$ is a composition of at most $\zeta(n)$ of
$5\dta$-ambient isotopies.
An induction on $i$ implies
$$h_1\N(P\o G(Z\x 1),Q'')\i\N(P(\im_M G(Z\x 1)),Q'').$$
By 5.1 the homeomorphism $H$ can be chosen arbitrarily, say,
$\frac{1}4\gma$-close to the identity.
By 5.3(i) and 5.4(ii) $\im_M G(Z\x 1)\i\N(G(Z\x 0),G(X\x 0))$.
Thus, $$\multline h_1\o P\o F(Z\x 1)\i h_1\N(P\o G(Z\x 1),Q'')\i\\
\i\N(G(\N(Z\x 0,X\x 0)),Q'')\i\N_{2\gma}(F(Z\x 0)),\endmultline$$
that is, $h_1\o P\o F(a\x 1)\i\N_{3\gma}(F(a\x 0))$ for any $a\in X$.
Recalling that $f=F_0$ and $g=P\o F_1$ are the given embeddings, we obtain that
$h_1\o g$ is $3\gma$-close to $f$.

Finally, we use 3.3a to obtain a $\dta$-ambient isotopy $\phi_t$, taking
$h_1\o g$ onto $f$.
Since $(5\zeta(n)+1)\dta=\eps$, the stacked composition of
$h_t$ and $\phi_t$ is an $\eps$-ambient isotopy, taking $g$ onto $f$. \qed
\enddemo

In the proof of 1.13b we will need the following observation.

\proclaim{Lemma 5.8}
For each positive integer $n$ there is a number $\rho(n)$ such that the
following holds.
Let $X^n$ be a compact polyhedron, $Q^m$ a PL manifold, $m-n\ge 3$, and
$f,g\:X\emb Q$ two PL $\dta$-concordant embeddings.

Given $\bta>0$, there are PL ambient isotopies $H_t^1,\dots,H_t^{\rho(n)}$
such that for each $i=1,\dots,\rho(n)$, $H_0^i=\id_Q$ and the isotopy $H_t^i$
has support in the disjoint union of sets of diameter $<7\dta$, and the
composition $H_1^{\rho(n)}\o\dots\o H_1^1\o g$ is $\bta$-close to $f$.
\endproclaim

\demo{Proof} This is clear from the proof of 1.13a, provided the following
modification is made.
(We use the notation from the proof of 1.13a.)

We can divide each $i$-th collapse ($*$), $i=1,\dots,\zeta(n)$, which we denote
for simplicity by $K_i\col L_i$, into $n+1$ collapses
$K_i\col K_{1,i}\cup L_i\col\dots\col K_{n,i}\cup L_i\col L_i$,
where $K_{ij}=\tr_{{\sssize K_i\searrow L_i}} P(\im_i G(X\x 1)^{(n-j)})$.
For each $i=1,\dots,\zeta(n)$, $j=0,\dots,n$, the set $K_{ij}\but K_{i,j+1}$
is the disjoint union of the sets $T_{i,Z}\but K_{i,j+1}$, where
$T_{i,Z}$ denotes $\tr_{{\sssize K_i\searrow L_i}} P(\im_i G(Z\x 1))$ and $Z$
runs over the $(n-j)$-simplices of $X$.

By 5.3(i) and 5.4(ii) each $T_{i,Z}$ is of diameter at most $\dta+4\gma$.
Now the sets $\N(\N(\text{\rm vertices in }(T_{i,Z}\but K_{i,j+1})',Q''),Q''')$
are each of diameter at most $\dta+6\gma<7\dta$, and are disjoint for distinct
$(n-j)$-simplices $Z$ of $X$.
Hence applying Lemma 5.7 for each $i=1,\dots,\zeta(n)$, $j=0,\dots,n$ to the
collapse $K_{ij}\cup L_i\col K_{i,j+1}\cup L_i$, we obtain a PL ambient isotopy
$H_t^{(n+1)(i-1)+j+1}$ with support in the disjoint union of the sets each of
diameter at most $7\dta$.
As in the proof of 1.13a it follows that
$H_1^{(n+1)\zeta(n)}\o\dots\o H_1^1\o g$ is $3\gma$-close to $f$.
Finally, the statement follows if we put $\rho(n)=(n+1)\zeta(n)$ and
$\gma=\mesh Q<\frac{\bta}3$. \qed
\enddemo

\demo{Proof of 1.13b}
Put $\dta=\frac{\eps}{29\rho(n)+2}$.
Let $\bta=\min(\dta,\frac{1}3\dta\ind{3.3b})$, where $\dta\ind{3.3b}$ is
obtained from 3.3b for $f\ind{3.3b}=f$ and $\eps\ind{3.3b}=\dta$.
By 3.1a and 3.5a $f$ is topologically $\bta$-isotopic to a PL embedding $f\PL$.
Let $\alp=\min(\dta,\frac{1}2\dta\ind{3.3b})$, where $\dta\ind{3.3b}$ is
obtained from 3.3b for $f\ind{3.3b}=g$ and $\eps\ind{3.3b}=\dta$.
By 3.1a and 3.5a $g$ is topologically $\alp$-isotopic to a PL embedding $g\PL$.
Now $g\PL$ is TOP $(\alp+\bta+\dta)$-concordant to $f\PL$.
Hence by (the relative case of) 3.1a the embeddings $f\PL$ and $g\PL$ are PL
$4\dta$-concordant.

Apply Lemma 5.8 to obtain a sequence of PL ambient isotopies
$H_t^1,\dots,H_t^{\rho(n)}$ such that for each $i=1,\dots,\rho(n)$ the isotopy
$H_t^i$ has support in the disjoint union of sets of diameter $<28\dta$,
and the composition $H_1^{\rho(n)}\o\dots\o H_1^1\o f\PL$ is
$\bta$-close to $g\PL$.
Each $H_t^i$ has a compact support, hence is uniformly continuous.
Let $\gma_{\rho(n)}=\bta$, and assuming that $\gma_i$ is defined, define
$\gma_{i-1}$ to be a number such that under $H_1^i$ any $\gma_{i-1}$-close
points are thrown into $\frac{1}2\gma_i$-close points.

Put $\gma=\min(\gma_0,\alp)$ and $\gma$-approximate $g\PL$ by a smooth
embedding $g_0\:X\emb Q$.
By 3.3b $g_0$ is smoothly $\dta$-ambient isotopic to $g$.
Now $g_0$ is topologically isotopic to the embedding $H_1^1\o g_0$ by means of
the isotopy $H_t^1\o g_0$, supported by the disjoint union $U$ of sets of
diameter $<28\dta$.
The embedding $H_1^1\o g_0$ is smooth outside $U$.
Hence by (the relative case of) 3.1b and by (the relative case of) 3.5b
this embedding is TOP $\frac{1}2\gma_1$-isotopic, fixing the exterior of an
arbitrarily small neighborhood $U'$ of $U$, to a smooth embedding
$g_1\:X\emb Q$.
The embedding $g_1$ is $\gma_1$-close to $H_1^1\o g\PL$, therefore by (the
relative case of) 3.4b $g_1$ is smoothly isotopic to $g_0$ by an isotopy $g_t$,
$t\in I$, fixing the exterior of an arbitrarily small neighborhood $U''$ of
$U'$.
Hence $g_t$ extends to a smooth ambient isotopy $G_t\:Q\to Q$ with support in
arbitrarily small neighborhood $U'''$ of $U''$, such that $G_0=\id_Q$ and
$G_1\o g_0=g_1$.
We can assume that $U'''$ is the disjoint union of sets of diameter $<29\dta$.
Consequently $g_0$ is smoothly $29\dta$-ambient isotopic to $g_1$.

Repeating the same construction for $i=2,3,\dots,\rho(n)$, we obtain a sequence
of smooth embeddings $g_2,\dots,g_{\rho(n)}\:X\emb Q$ such that $g_i$ and
$g_{i+1}$ are smoothly $29\dta$-ambient isotopic for each $i=0,\dots,\rho(n)-1$
and such that $g_{\rho(n)+1}$ is $\bta$-close to
$H_1^{\rho(n)}\o\dots\o H_1^1\o g\PL$.
Therefore $g_{\rho(n)+1}$ is $3\bta$-close to $f$.
Hence by 3.3b $g_{\rho(n)+1}$ and $f$ are smoothly $\dta$-ambient isotopic.
Thus $g$ is smoothly $\dta$-ambient isotopic to $g_0$, which, in turn, is
smoothly $29\dta\rho(n)$-ambient isotopic to $g_{\rho(n)}$, which is smoothly
$\dta$-ambient isotopic to $f$. \qed
\enddemo

\head 6. Proofs of Theorems 1.12 and 1.16 \endhead

\demo{Proof of 1.12}
We proceed with the first and the `moreover' parts simultaneously (in the first
part, let $\eps>0$ be any number).
Let $F\:X\x I\to Q\x I$ be the given (PL) pseudo-concordance.
The proof splits into two cases.
\enddemo

\demo{PL case} (compare to \cite{RoS, proof of Lemma 4.23 on level-preserving
collars}).
Fix some triangulations of $X\x I$, $Q\x I$ such that $\mesh Q<\frac{\eps}2$
and $F$ and the projections $p\:X\x I\to X$, $P\:Q\x I\to Q$ are simplicial.
Let $(X\x I)'$, $(Q\x I)'$ denote derived subdivisions of $X\x I$, $Q\x I$
which project simplicially onto the barycentrically derived subdivisions
$X'$, $Q'$ of $X$, $Q$.
For each simplex $A$ of $X\x I$ (resp\. $Q\x I$), we denote by $d_A$ its
derivation point in $(X\x I)'$ (resp\. in $(Q\x I)'$).
Let $\gma>0$ be so small that no vertex of $(X\x I)'$ lies in $X\x (1-\gma,1)$
and no vertex of $(Q\x I)'$ lies in $Q\x (1-\gma,1)$.
Then for each simplex $A$ of $X\x I$ (resp\. $Q\x I$) meeting $X\x 1$ (resp\.
$Q\x 1$) in a simplex $B$, the join $d_A*d_B$ meets $X\x\{1-\gma\}$ (resp\.
$Q\x\{1-\gma\}$) precisely in one point, which we denote by $d^+_A$.

We define a new PL pseudo-concordance $F_+\:X\x I\emb Q\x I$ as follows.
Put $F_+|_A=F|_A$ for any simplex $A$ not meeting $X\x\{1-\gma\}$.
Let $A_1,\dots,A_M$ be the simplices of $X\x I$ meeting $X\x\{1-\gma\}$,
arranged in some order of increasing dimension.
Assuming that $F_+|_{\partial A_i}$ is defined, define $F_+|_{A_i}$ by
$d^+_{A_i}\mapsto d^+_{F(A_i)}$ and extending linearly.
Then $F_+^{-1}(Q\x\{1-\gma\})=X\x\{1-\gma\}$ (moreover, $F_+|_{X\x [1-\gma,1]}$
is level-preserving, but we do not use this fact).
Also, $p\o F_+^{-1}(B\x\{1-\gma\})=p\o F_+^{-1}(B\x 1)=f^{-1}(B)$ for each
dual cell $B$ of $Q$.
Let $h$ denote the unique embedding $X\emb Q$ such that
$h\x\{1-\gma\}=F_+|_{X\x\{1-\gma\}}$.
By Lemma 4.2 $h$ is PL $\frac{\eps}2$-pseudo-isotopic to $f$, and this
completes the proof of the first part.

To prove the `moreover' part, notice that $F_+|_{X\x [0,1-\gma]}$ yields a PL
$\dta$-concordance between $g$, $h$.
By 1.13a $\dta$ can be chosen so that $g$ is PL $\frac{\eps}2$-ambient
isotopic to $h$. \qed
\enddemo

\demo{TOP case} In the `moreover' part we can assume (by 3.6b, 3.5a and 3.1a)
that the embedding $g\:X\emb Q$ is PL.
By the non-compact relative case of Theorem 3.1a the embedding $F|_{X\x [0,1)}$
can be assumed PL.

Put $t_0=t_1=0$.
Assuming $t_i$, $i>1$, to be already defined, put
$$t_{i+1}=\cases\sup\Pi\o F(X\x [0,\frac{1}2(1-t_i)])&\text{ if $i$ is even,}\\
\sup\pi\o F^{-1}(Q\x [0,\frac{1}2(1-t_i)])&\text{ if $i$ is odd,}\endcases$$
where $\pi\:X\x I\to I$, $\Pi\:Q\x I\to I$ denote the projections.
Clearly, $t_i<t_{i+1}<1$ for all $i=0,1,\dots$.
The main property of $t_i$'s is: for even $i$
$$F(X\x [t_{i-1},t_{i+1}])\i Q\x [t_{i-2},t_{i+2}].$$
Let $4\gma=\dta\ind{3.7}$ be given by Theorem 3.7 for
$\eps\ind{3.7}=\frac{1}2$.
Define a PL homeomorphism $\lda\:[0,1)\to [0,+\infty)$ by mapping
$t_i\mapsto\gma (i-1)$ for all odd $i$, an extending linearly.
Analogously, define a PL homeomorphism $\mu\:[0,1)\to [0,+\infty)$ by mapping
$t_i\mapsto\gma i$ for all even $i$, an extending linearly.
Let $$G=(\lda^{-1}\x \id_X)\o F\o (\mu\x \id_Q)\:
X\x [0,+\infty)\to Q\x [0,+\infty).$$
We obtain that $G(X\x [\gma(i-2),\gma i])\i Q\x [\gma(i-2),\gma(i+2)]$
for each even $i>0$.
Then by Theorem 3.7 there exists an ambient isotopy taking $G$ onto an
embedding
$G_+\:X\x [0,+\infty)\to Q\x [0,+\infty)$ such that $G_+^{-1}(Q\x i)=X\x i$ for
all $i=1,2,\dots$ (in the `moreover' part for $i=0$ in addition).

Let us assume $i$ to run over the positive integers in the proof of the
first part, and over the nonnegative integers in the `moreover' part.
Let $\alp_i=\dta\ind{1.13a}$ be given by Theorem 1.13a for
$\eps\ind{1.13a}=\frac{\eps}{2^{i+1}}$.
(In the `moreover' part we put $\dta=\frac{\alp_0}2$ in addition.)
Since $F$ is continuous, we can choose a sequence of integers $n_i$ such that
$P\o G|_{X\x [n_i,+\infty)}$ is $\frac{\alp_i}2$-close to
$f\x\id_{[n_i,+\infty)}$, where $P\:Q\x\R^1\to Q$ denotes the projection.
(In the `moreover' part the hypothesis allows to take $n_0=0$.)
By the furthermore part of 3.7 we can assume without loss of generality that
$P\o G_+|_{X\x [n_i,+\infty)}$ is $\frac{\alp_i}2$-close to
$P\o G|_{X\x [n_i,+\infty)}$ for each $i$.
Therefore $G_+|_{X\x [n_i,n_{i+1}]}$ is an $\alp_i$-concordance.
Hence by 1.13a for each $i$ the embeddings $f_i=P\o G_+|_{X\x n_i}$ and
$f_{i+1}$ are PL $\frac{\eps}{2^{i+1}}$-ambient isotopic.
It follows (cf\. \cite{Ke3, Lemma 1}) that $f_1$ ($f_0$ in the `moreover' part)
can be taken onto $f$ by an $\eps$-pseudo-isotopy. \qed
\enddemo

\demo{Proof of 1.16a}
We prove PL and DIFF cases simultaneously.
Suppose that an embedding $g'\:X\emb Q$ is taken onto $f$ by a pseudo-isotopy
$H'_t\:Q\to Q$.
By 3.1 and 3.5a there is a PL (smooth) embedding $g\:X\emb Q$, which is taken
onto $g'$ by a pseudo-isotopy $G_t\:Q\to Q$.
Then the `diagonal' pseudo-isotopy $H_t=H'_t\o G_t$ takes $g$ onto $f$. \qed
\enddemo

\demo{Proof of 1.16b}
The PL case was actually proved in the above proof of 1.12, TOP case.
Or it can be proved analogously to the below proof of the DIFF case:
\enddemo

\demo{DIFF case}
Let $H'_t\:Q\to Q$ be the given pseudo-isotopy and put $h_t=H'_t\o g$, where
$g\:X\emb Q$ is the given smooth embedding.
Let $\alp_i$, $i=1,2,\dots$, be a monotonely decreasing sequence of reals
(defined below) and let $t_i$, $i=1,2,\dots$ be such that $h_{t_i}$ is
$\alp_i$-close to $f$.
By 3.1b and 3.5b, $h_{t_i}$ is TOP $\alp_i$-isotopic to a smooth embedding
$g_i$, $i=1,2,\dots$.
Therefore $g_i$ and $g_{i+1}$ are TOP $3\alp_i$-isotopic for each
$i=1,2,\dots$, and $g$ is TOP isotopic to $g_1$.
By the relative case of 3.1b and by 1.13b $\alp_i$ can be chosen so that
$g_i$ and $g_{i+1}$ are smoothly $2^{-i-1}$-ambient isotopic,
while $g$ and $g_1$ are smoothly ambient isotopic.
It follows that $g$ can be taken onto $f$ by a pseudo-isotopy $H_t$ which is
smooth whenever $t\in [0,1)$. \qed
\enddemo

\remark{Remark 6.1}
We point out one useful observation following from Theorem 1.16 and the
relative version of Theorem 1.12 (which is proved analogously to 1.12,
using the relative version of 1.13a).
Suppose that $m-n\ge 3$, $X^n$ is a compact polyhedron, $Q^m$ a PL manifold,
and $f\:X\to Q$ a map which PL embeds a subpolyhedron $Z$ of $X$.
It turns out that if $f$ is isotopically realizable, then there is an embedding
$g\:X\emb Q$, agreeing with $f$ on $Z$, and a pseudo-isotopy $H_t\:Q\to Q$,
taking $g$ onto $f$ and keeping $g(Z)$ fixed.

Indeed, by 1.16 there is a PL embedding $g'$ and a pseudo-isotopy $H'_t$, PL
whenever $t\in [0,1)$ and taking $g'$ onto $f$.
From the non-compact case of Theorem 3.2a it follows that $H'_t\o g'|_Z$,
regarded as a level-preserving embedding $Z\x I\to Q\x I$, can be topologically
isotoped (not level-preserving in general) onto the embedding $f|_Z\x\id_I$.
Thus we obtain a pseudo-concordance from an embedding to $f$, keeping $Z$
fixed, and it remains to apply the relative version of 1.12, TOP case. \qed

Notice that if $Z$ is a manifold or $m-n\ge 4$, instead of Theorem 3.2a one
could apply its parametric version \cite{Mill3}, \cite{La}, \cite{Ste}, thus
making the application of the relative version of 1.12 no longer necessary.
\endremark

\head 7. Proof of Criterion 1.7 \endhead

The following is the controlled version of \cite{Harr, Theorem 1(R) on
page ~21}.

\proclaim{Theorem 7.1}
For each $\eps>0$ and a positive integer $n$ there exists $\dta>0$ such that
the following holds.
Let $X^n$ be a compact polyhedron, $(Q^m,\partial Q)$ a PL manifold,
$m\geq\frac{3(n+1)}2$, and $f\:X\to Q$ a PL map which embeds a subpolyhedron
$Z$ of $X$.

If $f^2$ is equivariantly $\dta$-homotopic to an isovariant map $H_1$ by a
homotopy $H\:X^2\x I\to Q^2$ which is isovariant on $Z^*_f\x I$, then $f$
is PL $\eps$-homotopic keeping $Z$ fixed to an embedding $g\:X\emb Q$ such
that $g^2$ is isovariantly $\eps$-homotopic to $H_1$.
\endproclaim

If $f\:X\to Q$ is a map into a manifold and $Z\i X$, we denote by $Z^*_f$ the
set $Z^2\cup (Z\cap f^{-1}\partial Q)\x X\cup X\x (Z\cap f^{-1}\partial Q)$.

Compare 7.1 to \cite{ReS2, pre-limit version of 1.2}.

Theorem 7.1 follows from 7.2 and 7.4 below, where $U\ind{7.4}=X^2$
and $U\ind{7.2}=V\ind{7.4}$.

\proclaim{Theorem 7.2}
For each $\eps>0$ and a positive integer $n$ there exists $\dta>0$ such that
the following holds.
Let $X^n$ be a compact polyhedron, $(Q^m,\partial Q)$ a PL manifold,
$m\geq\frac{3(n+1)}2$ and $f\:X\imm Q$ a PL immersion which embeds a
subpolyhedron $Z$.

Suppose that $f^2$ is equivariantly $\dta$-homotopic to an isovariant
map $H_1$ by a homotopy $H\:X^2\x I\to Q^2$ which is isovariant on
$(U\cup Z^*_f)\x I$, where $U$ is some neighborhood of $\Delta_X$ in $X^2$.

Then $f$ is PL $\eps$-regularly homotopic, keeping $Z$ fixed, to an embedding
$g\:X\emb Q$ such that $g^2$ is isovariantly $\eps$-homotopic to $H_1$.
\endproclaim

\demo{Proof}
This is a modification of \cite{Harr, proof of Theorem 3} in the spirit of
\cite{ReS2}.
Familiarity with \cite{Harr, proof of Proposition 9} is assumed.

Fix a triangulation of $X$ such that $f$ is simplicial, $Z$ is
a subcomplex, and the diameter of each simplex is less than $\dta$.
Arrange the simplices of $X$ so that first go that in
$Z\cap f^{-1}\partial Q$ in some order of increasing dimension, then the rest
in $Z$ in like order, and then the rest in $X$ in like order.
Equip $X^2$ with a cell structure induced by the triangulation of $X$.

Now 7.2 follows from the case $(p,q)=(n,n)$ of the below statement.
\enddemo

\proclaim{Claim 7.3} Let $p,q$ be some integers, $-1\le p\le q\le n$.
There is a positive integer $c_{pq}$ such that the following holds.

The immersion $f$ is PL $c_{pq}\dta$-regularly homotopic, keeping $L_p$ fixed,
to an immersion $f_{pq}\:X\imm Q$ such that $f_{pq}(A)\cap f_{pq}(B)=\emptyset$
for any cell $A\x B$ of $L_{pq}$, and such that $f_{pq}^2$ is equivariantly
$c_{pq}\dta$-homotopic to $H_1$ by a homotopy $H_{pq}\:X^2\x I\to Q^2$
which is isovariant on $(V\cup Z^*_f\cup L_{pq})\x I$.
\endproclaim

Here $L_p=X^{(p)}\cup Z$ and $L_{pq}=L_p\x L_q\cup L_q\x L_p\cup L_{p-1}\x X
\cup X\x L_{p-1}$.

\demo{Proof} The case $(p,q)=(-1,-1)$ follows from the hypothesis of 7.2,
assuming $c_{-1,-1}=1$.
The transition from $(p,n)$ to $(p+1,p+1)$ can be regarded as that from
$(p+1,p)$ to $(p+1,p+1)$.
We thus assume 7.3$(p,q)$ and prove 7.3$(p,q+1)$.

If two simplices $A,B$, where $A\x B$ is a cell of $L_{p,q+1}$, are mapped by
$f_{pq}$ sufficiently far from each other (namely, so that the minimal
distance between points in $f_{pq}(A)$, $f_{pq}(B)$ is greater than
$2c_{pq}\dta$), then $f_{pq}$ embeds $A\cup B$ and $H_{pq}$ is isovariant on
$(A\x B\cup B\x A)\x I$.

On the other hand, let $J_{pq}$ denote the set of pairs $(A_i,A_j)$ of
simplices of $X$ such that $\dim A_i=p$, $\dim A_j=q+1$, $A_j$ is not a simplex
of $Z$, if $p=q$ then $A_i$ preceds $A_j$ in the ordering, and such that
$\diam f_{pq}(A_i\cup A_j)<4c_{pq}\dta$.

By \cite{ReS2, 3.2} (which obviously generalizes for embedding into a PL
manifold) we can assume that there is a collection of PL balls $B_{ij}\i Q$,
where $(A_i,A_j)$ run over $J_{pq}$, such that $\diam B_{ij}<16c_{pq}\dta$,
$B_{ij}\cap B_{kl}=\emptyset$ whenever $(i,j)\neq (k,l)$,
$f_{pq}(A_i)\cap f_{pq}(A_j)$, if nonempty, lies in $\Int B_{ij}$, and $B_{ij}$
meets $A_i$, $A_j$ in PL balls $B_i^p\i\Int A_i$, $B_j^{q+1}\i\Int A_j$.

Now the argument of \cite{Harr, proof of Proposition 9} can be applied to each
pair $(A_i,A_j)=(\sma^p,\sma^{q+1})$ independently, so that there is a sequence
of ambient isotopies $H_{ij}$, $H'_{ij}$ of $Q$, supported by arbitrarily small
neighborhoods $N_{ij}$ of balls $B_{ij}$ (which can be chosen still disjoint
from each other and of diameter $<17c_{pq}\dta$).

The immersion $f_{p,q+1}$ is defined by $f_{p,q+1}=H'_{ij}\o H_{ij}\o f_{pq}$
on $N_{ij}\cap\st(A_j,X')$ and by $f_{p,q+1}=f_{p,q}$ elsewhere, in particular,
on $N_{ij}\cap\st(A_i,X')$.
It follows from the construction of $H$, $H'$ in \cite{Harr} that
$f_{p,q+1}(A)\cap f_{p,q+1}(B)$ is empty for any simplices $A$ of
$L_p=X^{(p)}\cup Z$ and $B$ of $L_{q+1}$ and for any simplices $A$ of $L_{p-1}$
and $B$ of $X$.
Also $f_{p,q+1}$ is $17c_{pq}\dta$-regularly homotopic to $f_{pq}$, keeping
$Z$ fixed.

The equivariant homotopy $H_{p,q+1}$ is defined analogously to as $F'$ in
\cite{Harr, proof of Proposition 9}, and it follows from the construction of
$F'$ in \cite{Harr} that $H_{p,q+1}$ is isovariant on
$(V\cup Z^*_f\cup L_{pq})\x I$.
From the obtained control of $f_{p,q+1}$ and the construction of $F'$ in
\cite{Harr} it follows that $H_{p,q+1}$ moves points less than
$17\sqrt{2}c_{pq}\dta$.

Hence if we take $c_{p,q+1}>(17\sqrt{2}+1)c_{pq}$, all conditions are
satisfied. \qed
\enddemo

\proclaim{Addendum to 7.2}
Let $H'\:X^2\x I\to Q^2$ denote the isovariant homotopy between $g^2$ and
$H_1$, let $r_t\:X\to Q$ denote the regular homotopy between $f$ and $g$,
and let $G\:X^2\x I\to Q^2$ be defined by $G_t=(r_{1-2t})^2$ for
$t\in [0,\frac{1}2]$ and $G_t=H_{2t-1}$ for $t\in [\frac{1}2,1]$.

Then $G$ and $H'$ are equivariantly $\eps$-homotopic by a homotopy isovariant
on $(X^2\x\partial I\cup (V\cup Z^*_f)\x I)\x I$, where
$V$ is some smaller neighborhood of $\Delta_X$ in $X^2$.
\endproclaim

This addendum is needed to carry out induction in the proof of 7.4 below.
It follows from the proofs of 7.2 and \cite{Harr, Proposition 9}.
(The similar addenda to 7.1 and 7.4 can be also shown to hold, but are not
required in this paper).

\proclaim{Theorem 7.4}
Let $X^n$ be a compact polyhedron, $(Q^m,\partial Q)$ be a PL manifold,
$m\geq\frac{3(n+1)}2$, and $f\:X\to Q$ be a PL map immersing a subpolyhedron
$Z$ of $X$.

Suppose that $f^2$ is equivariantly homotopic to map $H_1$ which is isovariant
on a neighborhood $U$ of $\Delta_X$ by a homotopy $H\:X^2\x I\to Q^2$ which
is isovariant on $(Z^*_f\cap U)\x I$.

Then for any $\eps>0$, $f$ is PL $\eps$-homotopic, keeping $Z$ fixed, to an
immersion $g\:X\imm Q$ such that $g^2$ is equivariantly homotopic to $H_1$
by a homotopy which is $\eps$-close to $H$ and isovariant on
$(V\cup (Z^*_f\cap U))\x I$ for some smaller neighborhood $V$ of $\Delta_X$.
\endproclaim

We minimize the program in \cite{Harr} by using a different method.
The idea is to construct a PL immersion by inductive gathering of certain PL
regular homotopies from immersions to embeddings. (This idea traces back to
discussions of A. ~Skopenkov and the author, and was first explicitely realized
by Skopenkov; see \cite{Sk2}.)
More precisely, the PL immersion will be constructed via an inductive
application of (7.2 + Addendum).

Actually we do not use the control in this application, that is, we could apply
simply (\cite{Harr, Theorem 3} + non-controlled Addendum).
Hence the below proof of 7.4, together with Harris' original proof of
\cite{Harr, Theorem 3}, yields a new short proof of \cite{Harr, Theorems 1, 2}.
It seems that an application of the control in (7.2 + Addendum) would be
necessary in a proof of a version of 7.4 with $C^1$-control.
Perhaps 7.4 can be also proved by a modification of \cite{Harr, proof of
Proposition 11} in the spirit of \cite{ReS2}.

We prove explicitely only the case $Z=\emptyset$.
Lack of $Z^*_f$ allows to use the following convention: we call a map
$f\:X^2\to Q^2$ (a homotopy $h\:X^2\x I\to Q^2$) {\it locally isovariant} if
there is a neighborhood $W$ of $\Delta_X$ in $X^2$ such that $f|_W$ (resp\.
$h|_{W\x I}$) is isovariant.

\demo{Proof}
Without loss of generality we can assume that $f$ is non-degenerate.
Triangulate $X$ and $Q$ so that $f$ is simplicial and diameter of each dual
cell $C$ of $Q$ is less than $\frac{\eps}{\sqrt 2}$.

Let $C_1,\dots,C_M$ be the dual cells of $Q$ arranged in an order of
increasing dimension, and let $C_{i1},\dots,C_{iM_i}$ be the dual cones of $X$,
whose disjoint union is $f^{-1}(C_i)$.
Each cell $C_i$ is dual to a simplex of $Q$ which we denote by $A_i$, whose
barycenter we denote by $b_i$, and whose link in $Q'$ we denote by $B_i$,
so that $b_i*B_i=C_i$.
We also write $D_i=b_i^2*B_i^2\i C_i^2$.
We define $A_{ij}$, $b_{ij}$, $B_{ij}$ and $D_{ij}$ analogously.

We write $\Cal S$ for $\{C_i\x C_j,b_i^2,D_i\mid i,j=1,\dots, M\}$.
We arrange $C_{ij}$ lexicographically, and denote $L_{pq}=C_{11}\cup\dots\cup
C_{pq}$.

We use the following notation. By $\eps(b*B)$ we denote, for any
$\eps\in (0,1]$ and any cone $b*B=B\x [0,1]/_{B\x 1}$, its subpolyhedron
$B\x [1-\eps,1]/_{B\x 1}$.
By $H_{\eps(b*B)}$ we denote the natural homeomorphism $(b*B\but\eps(b*B),B)\to
(B\x I,B\x 0)$.

Theorem 7.4 follows from 7.5 and the case $p=M$, $q=M_M$ of 7.8.
\enddemo

\proclaim{Claim 7.5} For each $\dta>0$ there is a locally isovariant homotopy
$G\:X^2\x I\to Q^2$, $\dta$-close to $H$ and such that $G_1=H_1$.
\endproclaim

\demo{Proof}
Clearly, $H|_{\Delta_X\x I}$ can be extended (analogously to the Borsuk Lemma)
to some locally isovariant homotopy $H'\:X^2\x I\to Q^2$ with $H'_1=H_1$.
There is a neighborhood $W$ of $\Delta_X$ in $X^2$ such that
$H'|_{W\x I}$ is $\dta$-close to $H|_{W\x I}$.

Moreover, if $\dta$ is sufficiently small, $H'|_{(\Fr W)\x I}$ is equivariantly
homotopic, by the `linear' homotopy, to $H|_{(\Fr W)\x I}$.
Hence by the Borsuk Lemma we can assume without loss of generality (but with
possible decrease of the neighborhood $W'$ of $\Delta_X$ in $X^2$ such
that $H'|_{W'\x I}$ is isovariant) that $H'|_{(\Fr W)\x I}=H|_{(\Fr W)\x I}$.

Define $G$ by $G=H'$ on $W\x I$ and by $G=H$ elsewhere, then $G$ is locally
isovariant and $\dta$-close to $H$. \qed
\enddemo

\proclaim{Lemma 7.6} Let $X^n$ be a finite simplicial complex,
$(Q^m,\partial Q)$ a combinatorial manifold, $m-n\ge 2$, $f\:X\to Q$ a
simplicial map and $C$ a union of some top-dimensional dual cells of $Q$.
Then for each $\eps>0$ there is $\dta>0$ such that the following holds.

Suppose that $G\:X^2\to Q^2$ is a locally isovariant map, $\dta$-close to $f^2$
and such that $G^{-1}(C_1\x C_2)=(f^2)^{-1}(C_1\x C_2)$, where $C_i$ run over
(dual cells of $C$)$\cup\{Q\}$.
Then $G$ is locally isovariantly $\eps$-homotopic, keeping $C^2$ fixed
and preserving $C\x Q\cup Q\x C$, to a map $F\:X^2\to Q^2$ such that
$F^{-1}(C_1\x C_2)=(f^2)^{-1}(C_1\x C_2)$ for each dual cells $C_1$, $C_2$
of $Q$.
\endproclaim

The proof of Lemma 7.6 modulo the below Lemma 7.7a is analogous to the proof of
Lemma 4.1 modulo Theorem 3.7 and we leave it for the reader.

\proclaim{Lemma 7.7} Suppose that $m-n\ge 1$.

(a) Let $X^n$ be a compact polyhedron, $Y$ its collared subpolyhedron,
$(Q^m,\partial Q)$ a PL manifold and $J_0=\{(-\infty,0],0,[0,+\infty)\}^2$.
Any locally isovariant map
$$f\:(X^2,X\x Y\cup Y\x X)\x [-1,1]^2\to (Q^2,\partial(Q^2))\x\R^2$$ such that
$f^{-1}((\partial Q)^2\x J)=Y^2\x (J\cap [-1,1]^2)$ for each $J\in J_0$ is
locally isovariantly homotopic, keeping $(\partial Q)^2\x\R^2$ fixed and
preserving $\partial(Q^2)\x\R^2$, to a map $g$ such that
$g^{-1}(Q^2\x J)=X^2\x (J\cap [-1,1]^2)$ for each $J\in J_0$.

(b) Let $c*X^n$ be cone over a compact polyhedron, $c*Q^m$ cone over a PL
sphere or a PL ball, and
$$f\:((c*X)^2,X\x (c*X)\cup (c*X)\x X)\to ((c*Q)^2,Q\x (c*Q)\cup (c*Q)\x Q)$$
a locally isovariant map such that $f^{-1}(Q^2)=X^2$.
Then $f$ is locally isovariantly homotopic, keeping $(Q\x (c*Q)\cup (c*Q)\x Q)$
fixed, to a mapping $g$ such that $g^{-1}(c^2*Q^2,c^2)=(c^2*X^2,c^2)$.

(c) Let $X^n$ be a compact polyhedron, $Y$ its subpolyhedron, $(Q^m,\partial Q)$
a PL manifold, $J_d=\{\{(x,y)\in\R^2\mid x\ge y\},\Delta_\R,
\{(x,y)\in\R^2\mid x\le y\}\}$, $K$ an equivariant subpoly\-hedron of $\R^2$.
Any locally isovariant map $f\:X^2\x [-1,1]^2\to Q^2\x\R^2$ such that
$$\gather
f^{-1}(\partial(Q^2)\x J_\Delta)=(X\x Y\cup Y\x X)\x (J_\Delta\cap [-1,1]^2)\\
\text{and }f^{-1}(Q^2\x(K\cap J_\Delta))=X^2\x(K\cap J_\Delta\cap [-1,1]^2)
\endgather$$
for each $J_\Delta\in J_d$, is locally isovariantly homotopic, keeping
$\partial(Q^2)\x\R^2\cup Q^2\x K$ fixed, to a map $g$ such that
$g^{-1}(Q^2\x J_\Delta)=X^2\x (J_\Delta\cap [-1,1]^2)$ for each
$J_\Delta\in J_d$.
\endproclaim

Lemma 7.7 can be regarded as a homotopy-theoretic version of Theorem 3.7.
We reduce (a) and (b) to (c), which is proved analogously to the geometric
proof of the Freudental Suspension Theorem.

\demo{Proof of 7.7a}
Let us write $R=\partial(Q^2)$ and $Z=Y\x X\cup X\x Y$.
Without loss of generality $f^{-1}(R\x J)=Z\x (J\cap [-1,1]^2)$ for each
$J\in J_0$.

By 7.7c with $K=(-\infty,0]\x[0,+\infty)\cup [0,+\infty)\x(-\infty,0]$ there
is a locally isovariant homotopy $f_t\:Z\x [-1,1]^2\to R\x\R^2$, preserving
$R\x J$ for each $J\in J_0$, from $f|_{Z\x [-1,1]^2}$ to a map $f_1$
such that $f_1^{-1}(R\x J_\Delta)=Z\x (J_\Delta\cap [-1,1])$ for each
$J_\Delta\in J_d$.

Hence (using collars on $Z$ and $R$) without loss of generality we can assume
that $f^{-1}(R\x J_\Delta)=Z\x (J_\Delta\cap [-1,1])$ for each
$J_\Delta\in J_d$.

Then we can apply 7.7c ($K=\emptyset$) to obtain a locally isovariant homotopy,
keeping $R\x\R^2$ fixed, from $f$ to a map $h$ such that
$h^{-1}(Q^2\x J_\Delta)=X^2\x(J_\Delta\cap [-1,1])$ for each $J_\Delta\in J_d$.

There is no obstruction in homotoping $h$ onto required $g$
keeping $R\x\R^2$ fixed. \qed
\enddemo

\demo{Proof of 7.7b} We consider the case $Q=S^m$, since the case $Q=B^m$ is
its corollary.
By general position we can assume that $f^{-1}(c^2)=c^2$,
by \cite{Harr, Appendix} we can assume that $f$ is PL, and using pseudo-radial
projection \cite{RoS} we can assume that $f^{-1}(\eps(c*S^m)^2)=\eps(c*X)^2$
for sufficiently small $\eps>0$.
Let $$\gather h\:\partial\eps(c*X)^2\to X\x(c*X)\cup(c*X)\x X,\\
H\:\partial\eps(c*S^m)^2\to S^m\x (c*S^m)\cup (c*S^m)\x S^m\endgather$$
be the natural homeomorphisms.
The map $$F=H\o f|_{\partial\eps(c*X)^2}\o h^{-1}\:X\x (c*X)\cup (c*X)\x X\to
S^m\x (c*S^m)\cup (c*S^m)\x S^m$$ is isovariant.
Notice that $X\x (c*X)\cup (c*X)\x X$ is homeomorphic to
$X*X$, i.e\. to $X\x X\x I/(X\x 0\cup X\x 1)$, the product $X\x X$ being thrown
onto the middle section $X\x X\x\frac{1}2$ of the join.
Now the pair $$(S^m*S^m\but(\Delta_{S^m}\x\tfrac{1}2),
(S^m\x S^m\but\Delta_{S^m})\x\tfrac{1}2)$$ with induced involution
$(x,y,t)\inv (y,x,1-t)$ equivariantly deformation retracts onto
$\nabla_{S^m}\x (I,\frac{1}2)$ (where $\nabla_{S^m}$ is the {\it antidiagonal}
$\{(x,-x)\in S^m\x S^m\}\cong S^m$) with involution $(z,t)\inv (-z,1-t)$.
If $r$ denotes the retraction, $r\o F|_{X*X\but\Delta_X\x\frac{1}2}$ is
equivariantly homotopic to a map $R$ such that $R^{-1}(\nabla_{S^m}\x\frac{1}2)
=X\x X\x\frac{1}2$.
Since $r$ is homotopy invertible, it follows that there is an isovariant
homotopy $F_t$ from $F$ to a map $F_+$ such that $F_+^{-1}((S^m)^2)=X^2$.

Let us consider the map $$f_+=H^{-1}\o F_+\o h\:\partial\frac{\eps}2(c*X)^2\to
\partial\frac{\eps}2(c*S^m)^2.$$
Let us extend $f_+$ linearly over $\frac{\eps}2(c*X)^2$, by $f_+=f$ on
$(c*X)^2\but\eps(c*X)^2$ and using $F_t$ on
$\eps(c*X)^2\but\frac{\eps}2(c*X)^2$ (more precisely, by
$f_+=H^{-1}_{\eps(c*S^m)^2}\o F_{2t}\o H_{\eps(c*X)^2}$).
Then $f_+$ is locally isovariant, $$f_+^{-1}(\frac{\eps}2(c*S^m)^2,
\frac{\eps}2(c^2*(S^m)^2))=(\frac{\eps}2(c*X)^2,\frac{\eps}2(c^2*X^2)),$$ and
$F_t$ yields a locally isovariant homotopy from $f$ to $f_+$.
To obtain the required map $g$, it remains to apply 7.7c with
$K=1\x[-1,1]\cup [-1,1]\x 1\cup [-1,-1+\eps]^2$ to
$f_+|_{(c*X)^2\but\frac{\eps}2(c*X)^2}$. \qed
\enddemo

\demo{Proof of 7.7c}
Let us write $R$ for $\partial (Q^2)\x\R^2\cup Q^2\x K$
and $Z$ for $$(X\x Y\cup Y\x X)\x\R^2\cup X^2\x K.$$
By \cite{Harr, Appendix} we can assume that $f$ is PL, meanwhile by
\cite{KuL, Theorem ~1} there exists a PL tangent bundle to $Q$, i.e\. a
collection of PL open ball pairs $(B^{2m}_i,D^m_i)\i (Q^2,\Delta_Q)$ such that
for each $i$ there exists a PL homeomorphism $H_i\:(B_i,D_i)\to(\R^{2m},\R^m)$,
making the following diagram commutative:
$$\CD \R^{2m}@>H_i>>B_i@>\text{\rm inclusion}>>Q^2\\
@V\text{\rm projection}VV @. @VV(q,r)\mapsto(q,q)V\\
\R^m@>H_i>>D_i@>\text{\rm inclusion}>>\Delta_Q \endCD$$
Fix equivariant triangulations on $X^2\x [-1,1]^2$ and $Q^2\x\R^2$ in which
$f$ and the projection $\dta\:Q^2\x\R^2\to Q^2\x\Delta_\R$ are simplicial.
We define $$\alignat 2 &\hat Q=\N(\Delta_{Q\x\R},Q^2\x\R^2),\qquad
&&\hat X=\N(\Delta_{X\x [-1,1]},X^2\x[-1,1]^2), \\
&\hat Q_0=\hat Q\cap Q^2\x\Delta_\R,
&&\hat X_0=\hat X\cap X^2\x\Delta_{[-1,1]}.\endalignat$$
Notice that $\hat X$ is a connected component of $f^{-1}(\hat Q)$.
We are to modify $f$ so that $\hat X_0$ is a connected component of
$f^{-1}(\hat Q_0)$.

We have that $f^{-1}(\Fr\hat Q\cap R)=\Fr\hat X\cap Z$.
It suffices to homotop $f|_{\Fr\hat X}$ onto a map
$$F_+\:(\Fr\hat X,\Fr\hat X_0)\to(\Fr\hat Q,\Fr\hat Q_0)$$
by a sufficiently small homotopy $F_t$ keeping $\hat Q\cap R$ fixed.
For if $F_t$ is so small that each simplex is left by it in the same
$B_i\x\R^2$, then it can be extended `linearly' (by an induction on joins
$Y*Z$, where $Y$ and $Z$ run, in orders of increasing dimension, over simplices
of $\Delta_Q$ and $\lk(Y,Q)\cap\Fr\hat X$, respectively) to an isovariant
homotopy $f_t$ taking $f|_{\hat X}$ onto a map
$$f_+\:(\hat X,\hat X_0)\to(\hat Q,\hat Q_0).$$
To arrange such a homotopy $F_t$, let us consider the analogues of the north
and the south poles:
$$\alignat 2 &\hat Q^n=(\Fr\hat Q)\cap\Delta_Q\x ((-\infty)\x(+\infty),0^2],
\quad &&\hat X^n=(\Fr\hat X)\cap\Delta_X\x [(-1)\x 1,0^2],\\
&\hat Q^s=(\Fr\hat Q)\cap\Delta_Q\x [0^2,(+\infty)\x(-\infty)),
&&\hat X^s=(\Fr\hat X)\cap\Delta_X\x [0^2,1\x(-1)].\endalignat$$
Since $\dta$ is simplicial, these are subcomplexes of $\Fr\hat Q,\Fr\hat X$,
while $\Fr\hat X\but(\hat X^n\cup\hat X^s)$ is equivariantly homeomorphic to
the cylinder $C=\Fr\hat X_0\x (0,1)$ equipped with the involution
$((x,y,s),t)\inv((y,x,s),-t)$.
We can assume (adjusting $K$, if necessary, without loss of generality)
that the homeomorphism takes $Z\cap(\Fr\hat X\but(\hat X^n\cup\hat X^s))$ onto
$Z_0\x (0,1)$ for some subpolyhedron $Z_0$ of $\Fr\hat X_0$.

By general position \cite{Bi3} the inequality $(2n-m+1)+(n+1)<2n+2$ implies
that $f^{-1}(\hat Q^n)$ does not meet $\hat X^n\cup\hat X^s$, while
$(2n-m+1)+(2n-m+1)<2n+2$ implies that $f^{-1}(\hat Q^n\cup\hat Q^s)$, regarded
as a subset of $\Fr\hat X_0\x (0,1)$, is not self-overshadowing.
Therefore there is an equivariant isotopy $h_t\:\Fr\hat X\to\Fr\hat X$, keeping
$\hat X^n\cup\hat X^s\cup Z_0\x (0,1)$ fixed and preserving the generators of
the cylinder $C$, and taking $f^{-1}(\hat Q^n)$ (resp\. $f^{-1}(\hat Q^s)$)
into a small neighborhood of $\hat X^n$ (resp\. $\hat X^s$).
The homotopy $H_t=f\o h_t^{-1}$ carries $f|_{\Fr\hat X}$ (equivariantly and
keeping $R\cap\Fr\hat Q$ fixed) onto a map $f_1\:\Fr\hat X\to\Fr\hat Q$
such that $f_1^{-1}(\hat Q^n)$ (resp\. $f_1^{-1}(\hat Q^s)$) is close to
$\hat X^n$ (resp\. $\hat X^s$), hence does not meet $\Fr\hat X_0$.

By stretching a neighborhood of $\hat Q^n$ over the `hemisphere'
$(\Fr\hat Q)\cap Q^2\x\{(x,y)\mid x\le y\}$ and similarly for the symmetric
neighborhood of $\hat Q^s$ (so that the complement to these neighborhoods is
squeezed onto the `equator' $\Fr\hat Q_0$) we obtain a homotopy (equivariant,
keeping $R\cap\Fr\hat Q$ fixed) onto a map $f_2$ such that
$f_2(\Fr\hat X_0)\i\Fr\hat Q_0$.
Final straightening in the spirit of the Alexander trick yields an equivariant
homotopy, keeping $R\cap\Fr\hat Q$ fixed, onto the required map
$f_+\:\Fr\hat X\to\Fr\hat Q$ such that $f_+^{-1}(\Fr\hat Q_0)=\Fr\hat X_0$.
\qed
\enddemo

\proclaim{Claim 7.8}
Let $F$ denote the map, obtained from Lemma 7.6 for $G\ind{7.6}=G_0$
and $C\ind{7.6}=\emptyset$, and refined by application of Lemma 7.7b
to cones $c*X\ind{7.6}$ running over all dual cones $C_{ij}$.
Then for each $p=0,\dots,M$, $q=0,\dots,M_p$ the following holds.

The map $f$ is PL homotopic to a map $f_{pq}\:X\to Q$, which immerses $L_{pq}$,
by a homotopy $h_{pq}\:X\x I\to Q$ such that $h_{pq}^{-1}(C_i)=f^{-1}(C_i)\x I$
for each $i\le M$.

Moreover, $f_{pq}^2$ is equivariantly homotopic to $F$ by a homotopy
$F_{pq}\:X^2\x I\to Q^2$ such that $F_{pq}^{-1}(S)=(f^2)^{-1}(S)\x I$ for each
$S\in\Cal S$, and which is locally isovariant on $L_{pq}^2\x I$.
\endproclaim

\demo{Proof} By lexicographic induction on $p,q$.
The base $p,q=0$ follows by taking $f_{00}=f$ and constructing $F_{00}$
inductively by the Alexander trick.
(Given polyhedra $Y\supset Y_1$, $Z\supset Z_1$ and maps
$\Psi_0,\Psi_1\:c*Y\to c*Z$ such that $\Psi_i^{-1}(c)=c$ and
$\Psi_i^{-1}(c*Z_1)=c*Y_1$, then any homotopy $\psi_t\:Y\to Z$ between
$\Psi_0|_Y$, $\Psi_1|_Y$ such that $\psi^{-1}(Z_1)=Y_1\x I$ can be extended
to a homotopy $\Psi_t\:c*Y\to c*Z$ between $\Psi_0$, $\Psi_1$ such that
$\Psi^{-1}(Z)=Y\x I$, $\Psi^{-1}(c*Z_1)=(c*Y_1)\x I$, and
$\Psi_t^{-1}(c)=c\x I$.)

To prove the inductive step, notice that $f_{p,q-1}|_{B_{pq}}\:B_{pq}\to B_p$
is an immersion, and that the homotopy $F_{p,q-1}|_{B_{pq}^2\x I}\:
B_{pq}^2\x I\to B_p^2$ is locally isovariant.

Also we have that $F_{p,q-1}|_{B_{pq}^2\x 1}=F|_{B_{pq}^2}$
is locally isovariantly homotopic to an isovariant map.
Indeed, define a homotopy $\Psi\:B_{pq}^2\x I\to B_p^2$ by
$$\Psi=\pi_{B_p^2}\o H_{\eps D_p}\o F\o H_{\eps D_{pq}}^{-1},$$
where $\eps$ is so small that $\eps D_{pq}$ lies in a neighborhood $W$ of
$\Delta_X$ in $X\x X$ such that $F|_W$ is isovariant.
Then by 7.2 there is a PL regular homotopy $r\:B_{pq}\x I\to B_p\x I$ from
$f_{p,q-1}|_{B_{pq}}$ to an embedding $r_1$.
We define $f_{pq}$ by $H_{\eps C_p}^{-1}\o r\o H_{\eps C_{pq}}$ on
$C_{pq}\but\eps C_{pq}$, conewise on $\eps C_{pq}$, and by $f_{p,q-1}$
outside $C_{pq}$.
Since $r$ is a PL regular homotopy and $r_1$ is a PL embedding,
$f_{pq}|_{C_{pq}}$ is a PL immersion.

Also $f_{pq}$ agrees with $f_{p,q-1}$ on $B_{pq}$, hence $f_{pq}|_{C_{pq}}$ is
linearly homotopic, keeping $B_{pq}$ fixed, to $f_{p,q-1}|_{C_{pq}}$.
We denote this homotopy by $l\:C_{pq}\x I\to C_p$, and we denote its extension
over $X$ by identity by $L\:X\x I\to Q$, so that $L_0=f_{p,q-1}$ and
$L_1=f_{pq}$.
Then $L$ together with $h_{p,q-1}$ yield the required PL homotopy
$h_{pq}\:X\x I\to Q$ between $f$ and $f_{pq}$.

Since $r$ is level-preserving, $l^{-1}(tC_p)=tC_{pq}\x I$ for each $t\in I$.
It follows that $(l^2)^{-1}(D_p)=(D_{pq})$ and hence
$(L^2)^{-1}(S)=(f^2)^{-1}(S)\x I$ for each $S\in\Cal S$.
Now $L^2$ and $F_{p,q-1}$ together give an equivariant homotopy
$F^+_{p,q-1}\:X^2\x I\to Q^2$ between $(f_{pq})^2$ and $F$ such that
$(F^+_{p,q-1})^{-1}(S)=(f^2)^{-1}(S)\x I$ for each $S\in\Cal S$,
and which is locally isovariant on $L_{p,q-1}^2\x I$.
We want to make it locally isovariant on $L_{pq}^2\x I$.

We define an equivariant homotopy $\Phi\:B_{pq}^2\x [-1,2]\to B_p^2$ by
$$\Phi_t=\cases r_{-t},\qquad &t\in [-1,0];\\
(F_{p,q-1}|_{B_{pq}^2\x I})_t,\qquad &t\in [0,1];\\
\Psi_{t-1},\qquad &t\in [1,2].\endcases$$
By 7.2 and Addendum, $\Phi$ is locally isovariantly homotopic to an isovariant
map $\Xi_1$ by a homotopy $\Xi\:B_{pq}^2\x [-1,2]\x I\to B_p^2$, which is in
addition isovariant on $B_{pq}^2\x\{-1,2\}\x I$.
Denote $B_p^2\x I\cup D_p\x\partial I$ by $E_p$.
Define an embedding $\phi_p\:B_p^2\x [-1,2]\emb E_p$ by the identity on
$B_p^2\x [0,1]$, by $H_{\eps D_p\x 0}^{-1}\o (t\mapsto -t)$ on
$B_p^2\x [-1,0]$, and by $H_{\eps D_p\x 1}^{-1}\o (t\mapsto t-1)$ on
$B_p^2\x [1,2]$.
Define analogously $E_{pq}$ and $\phi_{pq}$.
Then twice extending the homotopy $\phi_p\o\Xi_t\o\phi_{pq}^{-1}$ conewise, we
obtain a locally isovariant homotopy $\Xi^+\:E_{pq}\x I\to E_p\x I$ from
$F^+_{p,q-1}|_{E_{pq}}\:E_{pq}\to E_p$ to an isovariant map $\Xi^+_1$.

We define the required $F_{pq}$ by $F^+_{p,q-1}$ outside $C_{pq}^2$, by
$H_{\eps(D_p\x I)}^{-1}\o \Xi^+\o H_{\eps(D_{pq}\x I)}$ on
$(D_{pq}\x I)\but\eps(D_{pq}\x I)$, conewise on $\eps(D_{pq}\x I)$, and by
the relative Alexander trick on the rest of $C_{pq}^2\x I$.
Then $F_{pq}$ is locally isovariant on $L_{p,q-1}^2\x I$ and on
$C_{pq}^2\x I$.
Also $F_{pq}^{-1}(S)=(f^2)^{-1}(S)\x I$ for all $S\in\Cal S$,
and $\Delta_{L_{pq}}\i L_{p,q-1}^2\cup C_{pq}^2$, therefore $F_{pq}$
is locally isovariant on $L_{pq}^2\x I$. \qed
\enddemo

\demo{Proof of 1.7a+} The PL case follows immediately from 7.1, because
any equivariant map, $\dta$-close to a given equivariant map $f^2$, is
equivariantly $\dta$-homotopic to $f^2$, provided $\dta>0$ is sufficiently
small.

The TOP case follows from the PL case and simplicial approximation. \qed
\enddemo

It is convenient to call a homotopy $h_t:X^2\to Q^2$ {\it
pseudo-isovariant} if $h_t$ is isovariant for $t<1$.
If $f_0,f_1\:X\to Q$ are maps, let us say that an equivariant ($\dta$-)homotopy
$\phi_t\:X^2\to Q^2$ between $f_0^2$ and $f_1^2$ {\it ($\dta$-)holonomic}
if the homotopy $\phi_t$ is homotopic, in the class of equivariant
($\dta$-)homotopies between $f_0^2$ and $f_1^2$, to a homotopy $(f_t)^2$, where $f_t\:X\to Q$ is
some homotopy between $f_0$ and $f_1$.

\proclaim{Corollary 7.9}
For each $\eps>0$ and a positive integer $n$ there exists $\dta>0$ such that
the following holds.
Let $X^n$ be a compact polyhedron, $Q^m$ a PL manifold, and $m>\frac{3(n+1)}2$.

(a) If $f,g\:X\emb Q$ are PL embeddings such that $f^2$ and $g^2$
are isovariantly $\dta$-holonomically homotopic, then $f$ and $g$ are PL
$\eps$-ambient isotopic.

(b) If $f\:X\emb Q$ is a PL map and $g\:X\to Q$ a PL embedding
such that $f^2$ and $g^2$ are pseudo-isovariantly $\dta$-holonomically
homotopic, then $f$ and $g$ are PL $\eps$-pseudo-isotopic.
\endproclaim

Compare 7.9a to \cite{ReS2, Conjecture 1.9b}.
Notice that the statement of 7.9a is concerned with embeddings only, although
its proof and applications contain the idea of map realization.
Actually 7.1 and 7.9a are not only controlled versions of Harris' results,
but also generalize them (by considering an appropriate metric on $Q$).
Notice that from 7.1 and 7.9a at once follow 3.1a and 3.2a restricted
to the metastable range, respectively (compare to \cite{Harr, Corollary 4}).

Observe that if $Q=\R^m$ or if $\dta$ is allowed to depend on the metric on
$Q$, the holonomy condition is fulfilled automatically and can be dropped
from the hypotheses of 7.9a, 7.9b.

\demo{Proof of (a)} This is analogous to \cite{Harr, proof of Corollary 1}
(but not to \cite{Sk1, proof of theorems 1.0e and 1.1c}).

Let $\phi_t\:X^2\to Q^2$ be the given isovariant $\dta$-homotopy
between $f^2$ and $g^2$.
Let $\Omega_{s,t}\:X^2\to Q^2$ be the given homotopy, in the class
of equivariant $\dta$-homotopies between $f^2$ and $g^2$, from $\phi_t$
to $(f_t)^2$, where $f_t\:X\to Q$ is a homotopy between $f$ and $g$.
Let $F\:X\x I\to Q\x I$ be the PL map defined by $F(x,t)=(f_t(x),t)$.
Define an equivariant homotopy $\Psi_t\:(X\x I)^2\to (Q\x I)^2$ from
$F^2=\Psi_0$ to a map $\Psi_1$ by $$\Psi_t(x,u,y,v)=
(F(x,u(1-\tfrac{t}2)+v\tfrac{t}2),F(y,v(1-\tfrac{t}2)+u\tfrac{t}2)).$$
Then $\Psi_1(x,u,y,v)=(F(x,\frac{u+v}2),F(y,\frac{u+v}2))$,
which equals $((f_{(u+v)/2})^2(x,y),u,v)$.
Define an equivariant homotopy $\Phi_t\:(X\x I)^2\to (Q\x I)^2$ from
$\Psi_1=\Phi_0$ to an isovariant map $\Phi_1$ by $$\Phi_t((x,u,t,v))=
(\Omega_{t,(u+v)/2}(x,y),u,v).$$

Together $\Psi_t$ and $\Phi_t$ yield an equivariant $3\dta$-homotopy
$H_t$, which is isovariant on $(X\x\partial I)\x (X\x I)\cup (X\x I)\x
(X\x\partial I)$, from $F^2$ to an isovariant map $H_1=\Phi_1$.
Then by 7.1 for any $\gma>0$, $\dta$ can be chosen so that $F$ is
$\gma$-homotopic, fixing $X\x\partial I$, to a PL embedding, i.e\. to a PL
$(\gma+\dta)$-concordance between $f$ and $g$.
By 1.13a $\gma+\dta$ can be chosen so that $f$ and $g$ are PL $\eps$-ambient
isotopic. \qed
\enddemo

\demo{Proof of (b)}
Let $F$ be as above, and use the argument above to obtain an equivariant
homotopy $H_t\:(X\x I)^2\to (Q\x I)^2$, isovariant on
$$(X\x 0)\x (X\x I)\cup (X\x I)\x (X\x 0),$$ between $F^2=H_0$ and a map $H_1$,
isovariant on $$(X\x [0,1))\x (X\x I)\cup (X\x I)\x (X\x [0,1)).$$

Introduce an equivalence relation $\sigma$ on $X\x I$ by
$$(x,t)\sim (y,s): t=s=1\text{\rm \ and } f(x)=f(y).$$
Then the quotient space $M=X\x I/_\sigma$ is a compact polyhedron (the mapping
cylinder of $f$) and the identifying map $\Sigma\:X\x I\to M$ is PL.
By the construction of $H_t$, if $(x,u)\sim (x',u')$ and $(y,v)\sim (y',v')$,
then $H_t$ maps $(x,u,y,v)$ and $(x',u',y',v')$ to the same point of $Q^2$.
Hence there is the (unique) homotopy $H_t/_\sigma\:M^2\to (Q\x I)^2$ such that
$H_t/_\sigma\o\Sigma^2=H_t$.
Furthermore, $H_t$ maps $(x,u,y,v)$ to $\Delta_Q$, where either of $u$, $v$
equals to 1, only if $(x,u)\sim (y,v)$.
Therefore $H_t/_\sigma$ is isovariant on $$(X\x\partial I/_\sigma)\x
(X\x I/_\sigma)\cup (X\x I/_\sigma)\x (X\x\partial I/_\sigma)$$
and $H_1/_\sigma$ is isovariant everywhere.
Hence by 7.1 for any $\gma>0$, $\dta$ can be chosen so that
$F/_\sigma=H_0/_\sigma\:M\to Q\x I$ is $\gma$-homotopic, fixing
$X\x\partial I/_\sigma$, to a PL embedding, i.e\. to a PL
$(\gma+\dta)$-pseudo-concordance between $f$ and $g$.
By 1.12 (the PL case), $\gma+\dta$ can be chosen so that $f$ is PL
$\eps$-pseudo-isotopic to $g$. \qed
\enddemo

\demo{Proof of 1.7b+}
Let $\chi_t:X^2\to Q^2$ denote the equivariant $\dta$-homotopy such that
$\chi_1=f^2$, $\chi_0=g^2$ and $\chi_t$ is isovariant for $t<1$.
The proof splits into two cases.
\enddemo

\demo{PL case}
Since $\dta$ is allowed to depend on the metric on $Q$, the PL case follows
at once from 7.9b. \qed
\enddemo

\remark{Remark 7.10}
Notice that to prove only the PL case of 1.7b, rather than that of 1.7b+, one
can simplify the argument above by using the non-controlled version of 7.9b,
in whose proof one can use the first part of 1.12 instead of its `moreover'
part, and the Harris' Theorem itself instead of its controlled version 7.1.
This idea provides a somewhat simpler proof (not using the control in Harris'
Theorem) of the PL case of 1.7a for $m>\frac{3(n+1)}2$ (instead of
$m\ge\frac{3(n+1)}2$).
\endremark

\demo{TOP case}
By 3.6b, 3.5a and 3.1a we can assume that the embedding $g$ is PL.

For each positive integer $k$ let $\alp_k>0$ be some number (defined below).
Choose $t_k\in (0,1)$ so that $\chi_t$ moves points less than
$\frac{1}2\alp_k$ for $t\in [t_k,1]$.
Then $\chi_t$ is equivariantly $\alp_k$-homotopic to $f_k^2$ for any PL map
$f_k$, which is $\frac{1}2\alp_k$-homotopic to $f$, provided $\alp_k$ is
sufficiently small.

Then by 7.1, for any $\bta_k>0$, $\alp_k$ can be chosen so that there is a
PL embedding $\phi_k\:X\emb Q$, $\frac{1}3\bta_k$-close to $f$ and such that
$\phi_k^2$ is isovariantly $\frac{1}3\bta_k$-homotopic to $\chi_{t_k}$.
Also, $\chi_{t_k}$ is isovariantly $\frac{1}3\bta_k$-homotopic
to $\chi_{t_{k+1}}$ (provided $\frac{1}3\bta_k\ge\frac{1}2\alp_k$), which,
in turn, is isovariantly $\frac{1}3\bta_k$-homotopic to $\phi_{k+1}^2$
(provided $\bta_{k+1}<\bta_k$).
Thus, $\phi_k^2$ and $\phi_{k+1}^2$ are isovariantly $\bta_k$-homotopic.

Hence by 7.9a for any $\gma_k>0$, the number $\bta_k=\bta_k(\gma_k, n,
\text{ metric on } Q)$ can be chosen so that $\phi_k$ and $\phi_{k+1}$
are PL $\gma_k$-ambient isotopic.
Similarly, for any $\gma_0>0$, the number $\dta+\frac{1}3\bta_1$ can be
chosen so that $\phi_1$ and $g$ are PL $\gma_0$-isotopic.
Thus if we take $\gma_k=\frac{\eps}{2^{k+1}}$, where $k=0,1,2,\dots$, then
composing these isotopies for all $k$, we obtain an $\eps$-pseudo-isotopy
taking $f$ onto $g$. \qed
\enddemo

\head 8. Proof of Corollary 1.8 \endhead

In view of Criterion 1.7, it suffices to prove the following

\proclaim{Lemma 8.1}
Let $X^n$ be a compact polyhedron, $Q^m$ a PL manifold, $\phi:X^2\to Q^2$
an equivariant map.
If either

(a) $S=\phi^{-1}(\Delta_Q)$ has an equivariant mapping cylinder neighborhood
$N_1$ in $X^2$, or

(b) $m=2n+1$,

\noindent
then for each $\eps>0$ there exists $\dta>0$ such that any isovariant map
$\psi$, $\dta$-close to $\phi$, is $\eps$-homotopic to $\phi$ by a homotopy,
isovariant for $t<1$.
\endproclaim

Recall that for some neighborhood $U$ of $\Delta_Q$ in $Q\x Q$ the projection
$\tau$ onto the first factor is a PL (in general not vector) bundle over $Q$,
called tangent bundle in the PL category \cite{KuL}.
For each $t\in [0,1]$ let us denote by $U_t$ the total space of some subbundle
of this bundle with each fiber $U_{pt}$ of diameter $<t$, cf\. \cite{KuL}.

\demo{Proof of (a)}
By the definition, $N_1$ is equivariantly homeomorphic, by a homeomorphism $G$,
to the mapping cylinder $N\x I/_\sim$ of an equivariant map $g\:N\to S$, where
$(n_1,t)\sim(n_2,s)$ denotes `$g(n_1)=g(n_2)$ and $s=t=0$'.
Let $N_t$ be the $G$-preimage of $N\x [0,t]/_\sim$ for any $t\in (0,1]$ and
$N_{pt}$ be the $G$-preimage of $p\x [0,t]/_\sim$ for any $p\in N$,
$t\in (0,1]$.
Let us write simply $(p,t)$ for $G^{-1}((p,t)/_\sim)\in N_1$, if $t\in (0,1]$.
If $\dta>0$ is sufficiently small, there exists a `linear' equivariant
$\dta$-homotopy $h_t$ between $\phi$ and $\psi$.
Furthermore, by taking $\dta$ small enough we can achieve that
$h_t^{-1}(\Delta_Q)\i\Int N_1$ for each $t\in I$.

Let us construct an isovariant $2\dta$-homotopy $\psi_t$ between $\psi_0=\psi$
and $\psi_1$ such that $\psi_1$ is $2\dta$-close to $\phi$ and equals to $\phi$
outside $N_1$.
Define two functions on triangles: for each $t\in [0,1]$, let $\alp_t$ map
$[\frac{1}3,\frac{2}3t+1(1-t)]$ linearly onto $[\frac{1}3,1]$ and $\bta_t$ map
$[\frac{2}3t+1(1-t),1]$ linearly onto $[0,t]$.
Put $\psi_t=\psi_0$ on $N_{1/3}$, $\psi_t(p,s)=\psi_0(p,\alp_t(s))$ for
$s\in [\frac{1}3,\frac{2}3t+1(1-t)]$, $\psi_t(p,s)=h_{\bta_t(s)}(p,1)$ for
$s\in [\frac{2}3t+1(1-t),1]$, and $\psi_t=h_t$ outside $N_1$.
Then $\psi_t$ is clearly isovariant.
Moreover, $h_t$ is $2\dta$-close to $\phi$, provided $N_1$ is so small that
$\phi(p,s)$ and $\phi(p,t)$ are $\dta$-close for any $p\in N$, $s,t\in (0,1]$.

Choose $N_1$ so small that $\phi(N_1)$ and $\psi_1(N_1)$ lie in $U_1$.
Let $\alp,\bta\:[1,+\infty)\to (0,1]$, be such homeomorphisms that
$\phi(N_{\alp(t)})\i U_{\bta(t)}$ and $\diam\phi(N_{p\alp(t)})$,
$\diam U_{p\bta(t)}$ are less than $r_t=\frac{\eps}{11*2^{t+1}}$ for each $t$.
Let us write $(p,t),N_t,N_{pt},U_t,U_{pt}$ instead of $(p,\alp(t)),
N_{\alp(t)},N_{p\alp(t)},U_{\bta(t)},U_{p\bta(t)}$.
If $\dta<\frac{\eps}{44}$, the statement follows from
\enddemo

\proclaim{Claim 8.2}
Let $k$ be a positive integer and $\psi_k\:X^2\to Q^2$ an isovariant map,
$r_k$-close to $\phi$ and coinciding with $\phi$ outside $N_k$.
Then $\psi_k$ is isovariantly $(11r_k)$-homotopic to a map $\psi_{k+1}$,
$r_{k+1}$-close to $\phi$ and coinciding with $\phi$ outside $N_{k+1}$.
\endproclaim

\demo{Proof}
Define a homotopy $\Psi_t\:X^2\to Q^2$, $t\in I$, by $\Psi_0=\psi_k$,
the identity outside $N_k$, the identity on $N_{k+2}$ and by
$$\Psi_t(p,s)=\left\{
\alignedat 2
&\psi_k(p,s-(s-k)*2t),
 &\ t\le \tfrac{1}2,
 &\ k\le s\le k+\tfrac{3}2      \\
&\psi_k(p,s-(k+2-s)*3*2t),
 &\ t\le \tfrac{1}2,
 &\ k+\tfrac{3}2\le s\le k+2    \\
&\phi(p,k+(s-k)*(2t-1)),
 &\ t\ge \tfrac{1}2,
 &\ k\le s\le k+1               \\
&\phi(p,k+(k+\tfrac{3}2-s)*2*(2t-1)),
 &\ t\ge \tfrac{1}2,
 &\ k+1\le s\le k+\tfrac{3}2    \\
&\psi_k(p,s-(k+2-s)*3*1),
 &\ t\ge \tfrac{1}2,
 &\ k+\tfrac{3}2\le s\le k+2    \\
\endalignedat \right. $$
In other words, $\Psi_t$ stretches $\psi_k(N_{k+\frac{3}2})$ over $\psi_k(N_k)$
and takes $\psi_k(\Cl{N_k\but N_{k+\frac{3}2}})$ onto
$\phi(\Cl{N_k\but N_{k+1}})$.
The result is that $\Psi_1$ coincides with $\phi$ outside $N_{k+1}$.
Clearly, $\Psi_t$ is isovariant and moves points less than $3r_k$, meanwhile
$\Psi_1$ is $3r_k$-close to $\phi$.

We have that $\phi(N_{k+1})\i U_{k+1}$.
Homotop $\Psi_1$ isovariantly, fixing the exterior of $N_{k+1}$, `linearly'
towards $\phi$ onto a map $\Psi_2$ such that $\Psi_2(N_{k+1})\i U_{k+1}$.
This homotopy moves points less than $3r_k$, and $\Psi_2$ is still
$3r_k$-close to $\phi$.

Now both $\phi(N_{k+1})$ and $\Psi_2(N_{k+1})$ lie in $U_{k+1}$.
Since tangent bundle is locally trivial, we can homotop $\Psi_2$ isovariantly,
fixing the exterior of $N_{k+1}$, leaving the image of $N_{k+1}$ inside
$U_{k+1}$, and moving points less than $5r_k$, onto a map $\psi_{k+1}$ such
that $\tau\o\psi_{k+1}|_{N_{k+1}}=\tau\o\phi|_{N_{k+1}}$.
Then the images of any point under $\phi$ and $\psi_{k+1}$ lie in the same
$U_{p,k+1}$, which implies that $\psi_{k+1}$ is $r_{k+1}$-close to $\phi$. \qed
\enddemo

\demo{Proof of (b)}
The following argument was partially inspired by \cite{AgRS, proof of
Proposition 3.5}.
Given a compactum $C$ in the body of a simplicial complex $K$, we denote by
$\N(C,K)$ the union of all simplices of $K$, meeting $C$.
Let $K_1$ be a triangulation of $X^2$, and for each positive integer $k$ let
$K_{k+1}$ be a subdivision of $K_k$ such that for any simplex $A$ of $K_k$
$\diam\phi(A)<n_k=\frac{\eps}{14*2^{k+1}}$ and $\phi(A)\i U_{n_k}$.

Let $S=\phi^{-1}(\Delta_Q)$ and $N_k=\N(S,K_k)$, then $\phi(N_k)\i U_{n_k}$.
If $\dta>0$ is sufficiently small, there exists a `linear' equivariant
$\dta$-homotopy $h_t$ between $\phi$ and $\psi$ such that
$h_t^{-1}(\Delta_Q)\i\Int N_1$ for each $t\in I$.
By the equivariant Borsuk Lemma, $\psi$ is isovariantly $2\dta$-homotopic to a
map $\psi_1\:X^2\to Q^2$ which coincides with $\phi$ outside $N_1$ and is
$2\dta$-close to $\phi$, provided $\mesh K_1<\dta$.
If $\dta<\frac{\eps}{56}$, the statement follows from
\enddemo

\proclaim{Claim 8.3}
Let $k$ be a positive integer and $\psi_k\:X^2\to Q^2$ an isovariant map,
$n_k$-close to $\phi$ and coinciding with $\phi$ outside $N_k$.
Then $\psi_k$ is isovariantly $(14n_k)$-homotopic to a map $\psi_{k+1}$,
$n_{k+1}$-close to $\phi$ and coinciding with $\phi$ outside $N_{k+1}$.
\endproclaim

\demo{Proof}
Observe that in each simplex of $N_k$ there is at least one simplex of
$N_{k+1}$.
Therefore $P_{k+1}=\Cl{X^2\but N_{k+1}}$ equivariantly strong deformation
retracts onto $P_k\cup Y$, where $Y^{2n-1}$ is an equivariant subpolyhedron of
$N_k$ of dimension $2n-1$.
Let $r_t\:P_{k+1}\to P_{k+1}$ denote the deformation retraction, so that
$r_1(P_{k+1})=P_k\cup Y$; we can assume that $r_t$ moves points $<n_k$.

Now $\psi_k$ coincides with $\phi$ on $P_k$ and is $n_k$-close to $\phi$
elsewhere.
Hence there is a `linear' equivariant $n_k$-homotopy $g_t\:X^2\to Q^2$ from
$\psi_k$ to $\phi$, keeping $P_k$ fixed.
By equivariant general position we can assume that $g_t|_{Y^{2n-1}}$ does not
meet $\Delta_Q$.

Define an equivariant $(3n_k)$-homotopy $\Psi_t\:P_{k+1}\to Q^2\but\Delta_Q$
between $\psi_k|_{P_{k+1}}$ and $\phi|_{P_{k+1}}$ as follows.
Put $\Psi_t=g_0\o r_{3t}$ for $t\in [0,\frac{1}3]$,
$\Psi_t=g_{3t-1}\o r_1$ for $t\in [\frac{1}3,\frac{2}3]$, and
$\Psi_t=g_1\o r_{3-3t}$ for $t\in [\frac{2}3,1]$.
By the equivariant Borsuk Lemma $\Psi_t$ extends to an isovariant
$(4n_k)$-homotopy $\Psi_t\:X^2\to Q^2$ between $\psi_k$ and a map $\Psi_1$ such
that $\Psi_1=\phi$ on $P_{k+1}$ and $\Psi_1$ is $4n_k$-close to $\phi$.
The rest of the proof goes as in 8.2. \qed
\enddemo

\remark{Remark 8.4}
One can give an alternative proof of 8.3 using equivariant obstruction theory
(see \cite{CF}) as follows.
Since the pair $(P_{k+1},P_k)$ is equivariantly homotopy equivalent to
$(P_k\cup Y^{2n-1},P_k)$, the equivariant cohomology group
$H^{2n}\eq(P_{k+1},P_k,\pi)$ is zero for any group $\pi$.
Since $\psi_k$ and $\phi$ are close, and
$\pi_i(\R^{2n+1}\x\R^{2n+1}\but\Delta_{\R^{2n+1}})=0$ for $i<2n$, the
difference cocycle
$d(\psi_k,\phi)\in H^i\eq(P_{k+1},P_k,\pi_i(Q^2\but\Delta_Q))$ can be shown
to vanish for $i<2n$.
Hence $\psi_k$ can be isovariantly homotoped, keeping $P_k$ fixed, so as to
agree with $\phi$ on $P_{k+1}$.
\endremark

\enddocument
\end